\documentclass[11pt]{article}
\usepackage[english]{babel}
\usepackage{amssymb,amsmath,amsthm}
\newtheorem{theorem}{Theorem}[section]
\newtheorem{lemma}[theorem]{Lemma}
\newtheorem{proposition}[theorem]{Proposition}
\newtheorem{corollary}[theorem]{Corollary}
\newtheorem{question}[theorem]{Problem}

{\theoremstyle{definition}}
{\theoremstyle{definition}\newtheorem{example}[theorem]{Example}}
{\theoremstyle{definition}\newtheorem{definition}[theorem]{Definition}}
{\theoremstyle{definition}\newtheorem{remark}[theorem]{Remark}}

\newtheorem*{thmbg}{Theorem BG}
\newtheorem*{thmdu}{Theorem U}
\newtheorem*{thmd}{Proposition G}
\newtheorem*{thmsas}{Theorem B}
\newtheorem*{thmbp}{Theorem BP}
\newtheorem*{thmleon}{Theorem LM}
\newtheorem*{thmhc}{Theorem HC}
\newtheorem*{thmsc}{Theorem SC}
\newtheorem*{thmbm}{Theorem BM}

\numberwithin{equation}{section}

\font\Goth=eufm10 scaled 1440
\font\goth=eufm10 scaled 1200
\def\ccc{\hbox{{\Goth c}}}
\def\C{{\mathbb C}}
\def\Q{{\mathbb Q}}
\def\N{{\mathbb N}}
\def\Z{{\mathbb Z}}
\def\R{{\mathbb R}}
\def\T{{\mathbb T}}
\def\D{{\mathbb D}}

\def\M{\hbox{\goth M}}
\def\K{{\mathbb K}}
\def\F{{\mathcal F}}
\def\H{{\mathcal H}}
\def\E{{\mathcal E}}

\def\rr{{\mathcal R}}
\def\pp{{\mathcal P}}
\def\kk{{\bf k}}
\def\epsilon{\varepsilon}
\def\kappa{\varkappa}
\def\phi{\varphi}
\def\leq{\leqslant}
\def\geq{\geqslant}
\def\ilim{\mathop{\hbox{$\underline{\hbox{\rm lim}}$}}\limits}
\def\slim{\mathop{\hbox{$\overline{\hbox{\rm lim}}$}}\limits}

\def\dim{\hbox{\tt dim}\,}
\def\ker{\hbox{\tt ker}\,}
\def\Ker{\hbox{\tt ker}^\star\,}
\def\KER{\hbox{\tt ker}^\dagger\,}
\def\spann{\hbox{\tt span}\,}
\def\det{\hbox{\tt det}\,}
\def\deg{\hbox{\tt deg}\,}
\def\ssub#1#2{#1_{{}_{{\scriptstyle #2}}}}
\def\lbs{LB$_{\displaystyle\rm s}$}
\def\lbss{LB$_{\displaystyle s}$}
\def\NNN{{\mathcal M}^{\scriptscriptstyle \mathcal N}}

\title{Existence theorems in linear chaos}

\author{Stanislav Shkarin}

\date{}

\begin{document}

\maketitle

\begin{abstract} Chaotic linear dynamics deals primarily with
various topological ergodic properties of semigroups of continuous
linear operators acting on a topological vector space. We treat
questions of characterizing which of the spaces from a given class
support a semigroup of prescribed shape satisfying  a given
topological ergodic property.

In particular, we characterize countable inductive limits of
separable Banach spaces that admit a hypercyclic operator, show that
there is a non-mixing hypercyclic operator on a separable infinite
dimensional complex Fr\'echet space $X$ if and only if $X$ is
non-isomorphic to the space $\omega$ of all sequences with
coordinatewise convergence topology. It is also shown for any
$k\in\N$, any separable infinite dimensional Fr\'echet space $X$
non-isomorphic to $\omega$ admits a mixing uniformly continuous
group $\{T_t\}_{t\in \C^n}$ of continuous linear operators and that
there is no supercyclic strongly continuous operator semigroup
$\{T_t\}_{t\geq 0}$ on $\omega$. We specify a wide class of
Fr\'echet spaces $X$, including all infinite dimensional Banach
spaces with separable dual, such that there is a hypercyclic
operator $T$ on $X$ for which the dual operator $T'$ is also
hypercyclic. An extension of the Salas theorem on hypercyclicity of
a perturbation of the identity by adding a backward weighted shift
is presented and its various applications are outlined.
\end{abstract}


\normalsize

\section{Introduction \label{s1}}\rm

Unless stated otherwise, all vector spaces in this article are over
the field $\K$, being either the field $\C$ of complex numbers or
the field $\R$ of real numbers, all topological spaces {\it are
assumed to be Hausdorff} and all vector spaces {\it are assumed to
be non-trivial}. As usual, $\T=\{z\in\C:|z|=1\}$, $\Z$ is the set of
integers, $\Z_+$ is the set of non-negative integers, $\N$ is the
set of positive integers and $\R_+$ is the set of non-negative real
numbers. Symbol $L(X,Y)$ stands for the space of continuous linear
operators from a topological vector space $X$ to a topological
vector space $Y$. We write $L(X)$ instead of $L(X,X)$ and $X'$
instead of $L(X,\K)$. For each $T\in L(X)$, the dual operator
$T':X'\to X'$ is defined as usual: $(T'f)(x)=f(Tx)$ for $f\in X'$
and $x\in X$. By a {\it quotient} of a topological vector space $X$
we mean the space $X/Y$, where $Y$ is a closed linear subspace of
$X$. We start by recalling some definitions and facts.

\subsection{Notation and Definitions}

A topological vector space is called {\it locally convex} if it has
a base of neighborhoods of zero consisting of convex sets.
Equivalently, a topological vector space is locally convex if its
topology can be defined by a family of seminorms. For brevity, we
say {\it locally convex space} for a locally convex topological
vector space. A subset $B$ of a topological vector space $X$ is
called {\it bounded} if for any neighborhood $U$ of zero in $X$, a
scalar multiple of $U$ contains $B$. We say that $\tau$ is a {\it
locally convex topology} on a vector space $X$ if $(X,\tau)$ is a
locally convex space. If $X$ is a vector space and $Y$ is a linear
space of linear functionals on $X$ separating points of $X$, then
the weakest topology on $X$, with respect to which all functionals
from $Y$ are continuous, is denoted $\sigma(X,Y)$. The elements of
$X$ can be naturally interpreted as linear functionals on $Y$, which
allows one to consider the topology $\sigma(Y,X)$ as well. If $\cal
B$ is a family of bounded subsets of $(Y,\sigma(Y,X))$, whose union
is $Y$, then the seminorms
$$
p_B(x)=\sup_{f\in B}|f(x)|\quad \text{for}\quad B\in{\cal B}
$$
define the {\it topology on $X$ of uniform convergence on sets of
the family $\cal B$}. The topology on $X$ of uniform convergence on
all bounded subsets of $(Y,\sigma(Y,X))$ is called the {\it strong
topology} and denoted $\beta(X,Y)$. The topology of uniform
convergence on all compact convex subsets of $(Y,\sigma(Y,X))$ is
called the {\it Mackey topology} and is denoted $\tau(X,Y)$.
According to the Mackey--Arens theorem for a locally convex space
$(X,\tau)$, and a space $Y$ of linear functionals on $X$, the
equality $Y=X'$ holds if and only if $\sigma(X,Y)\subseteq
\tau\subseteq \tau(X,Y)$. We say that a locally convex space
$(X,\tau)$ {\it carries a weak topology} if $\tau$ coincides with
$\sigma(X,Y)$ for some space $Y$ of linear functionals on $X$,
separating points of $X$. If $X$ is a locally convex space, we write
$X_\beta$ for $(X,\beta(X,X'))$, $X_\tau$ for $(X,\tau(X,X'))$ and
$X_\sigma$ for $(X,\sigma(X,X'))$. Similarly we denote
$(X',\beta(X',X))$ by $X'_\beta$, $(X',\tau(X',X))$ by $X'_\tau$ and
$(X',\sigma(X',X))$ by $X'_\sigma$. An $\F$-space is a complete
metrizable topological vector space. A locally convex $\F$-space is
called a {\it Fr\'echet} space. If $\{X_\alpha:\alpha\in A\}$ is a
family of locally convex spaces, then their (locally convex) {\it
direct sum} is the algebraic direct sum
$X=\bigoplus\limits_{\alpha\in A}X_\alpha$ of the vector spaces
$X_\alpha$ endowed with the strongest locally convex topology, which
induces the original topology on each $X_\alpha$. Let
$\{X_n\}_{n\in\Z_+}$ be a sequence of vector spaces such that $X_n$
is a subspace of $X_{n+1}$ for each $n\in\Z_+$ and each $X_n$
carries its own locally convex topology $\tau_n$ such that $\tau_n$
is (maybe non-strictly) stronger than the topology
$\tau_{n+1}\bigr|_{X_n}$. Then the {\it inductive limit} of the
sequence $\{X_n\}$ is the space $X=\bigcup\limits_{n=0}^\infty X_n$
endowed with the strongest locally convex topology $\tau$ such that
$\tau\bigr|_{X_n}\subseteq\tau_n$ for each $n\in\Z_+$. In other
words, a convex set $U$ is a neighborhood of zero in $X$ if and only
if $U\cap X_n$ is a neighborhood of zero in $X_n$ for each
$n\in\Z_+$. An {\it LB-space} is an inductive limit of a sequence of
Banach spaces. An {\it \lbss-space} is an inductive limit of a
sequence of separable Banach spaces. We use symbol $\phi$ to denote
the locally convex direct sum of countably many copies of the
one-dimensional space $\K$ and the symbol $\omega$ to denote the
product of countably many copies of $\K$. Note that $\phi$ is a
space of countable algebraic dimension and carries the strongest
locally convex topology (=any seminorm on $\phi$ is continuous). We
can naturally interpret $\omega$ as the space $\K^{\N}$ of all
sequences with coordinatewise convergence topology. Clearly $\omega$
is a separable Fr\'echet space. Recall also that if $X$ is a locally
convex space and $A\subset X'$, then $A$ is called {\it uniformly
equicontinuous} if there exists a neighborhood $U$ of zero in $X$
such that $|f(x)|\leq 1$ for any $x\in U$ and $f\in A$.

Let $T$ be a continuous linear operator on a topological vector
space $X$. A vector $x\in X$ is called a {\it cyclic vector} for $T$
if the linear span of the orbit $O(T,x)=\{T^n:n\in\Z_+\}$ of $x$ is
dense in $X$. The operator $T$ is called {\it cyclic} if $T$ has a
cyclic vector. Recall also that for $n\in\N$, $T$ is called $n$-{\it
cyclic} if there are vectors $x_1,\dots,x_n\in X$ such that the
linear span of the set $\{T^nx_j:n\in\Z_+,\ 1\leq j\leq n\}$ is
dense in $X$. Obviously, 1-cyclicity coincides with cyclicity. We
say that $T$ is {\it multicyclic} if it is $n$-cyclic for some
$n\in\N$.

Let $X$ and $Y$ be topological spaces and $\{T_a:a\in A\}$ be a
family of continuous maps from $X$ to $Y$. An element $x\in X$ is
called {\it universal} for this family if the orbit $\{T_ax:a\in
A\}$ is dense in $Y$ and $\{T_a:a\in A\}$ is said to be {\it
universal} if it has a universal element. We say that a family
$\{T_n:n\in\Z_+\}$ is {\it hereditarily universal} if any its
infinite subfamily is universal. An {\it operator semigroup} on a
topological vector space $X$ is a family $\{T_t\}_{t\in A}$ of
elements of $L(X)$ labeled by elements of an abelian monoid $A$
(monoid is a semigroup with identity) and satisfying $T_0=I$,
$T_{s+t}=T_tT_s$ for any $t,s\in A$ (unless stated otherwise, we use
additive notation for the operation on $A$). A {\it norm} on $A$ is
a function $|\cdot|:A\to[0,\infty)$ satisfying $|na|=n|a|$ and
$|a+b|\leq |a|+|b|$ for any for any $n\in\Z_+$ and $a,b\in A$. An
abelian monoid equipped with a norm will be called a {\it normed
semigroup}. We will be mainly concerned with the case when $A$ is a
closed (additive) subsemigroup of $\R^k$ containing $0$ with the
norm $|a|$ being the Euclidean distance from $a$ to $0$. In the
latter case we consider $A$ to be equipped with topology inherited
from $\R^k$ and we say that an operator semigroup $\{T_t\}_{t\in A}$
is {\it strongly continuous} if the map $t\mapsto T_tx$ from $A$ to
$X$ is continuous for any $x\in X$. We say that an operator
semigroup $\{T_t\}_{t\in A}$ is {\it uniformly continuous} if there
exists a neighborhood $U$ of zero in $X$ such that for any sequence
$\{t_n\}_{n\in\Z_+}$ of elements of $A$ converging to $t\in A$,
$T_{t_n}x$ converges to $T_tx$ uniformly on $U$. Clearly, uniform
continuity is strictly stronger than strong continuity.

If $A$ is a normed semigroup and $\{T_t\}_{t\in A}$ is an operator
semigroup on a topological vector space $X$, then we say that
$\{T_t\}_{t\in A}$ is (topologically) {\it transitive} if for any
non-empty open subsets $U,V$ of $X$, the set $\{|t|:t\in A,\
T_t(U)\cap V\neq\varnothing\}$ is unbounded. We say that
$\{T_t\}_{t\in A}$ is (topologically) {\it mixing} if for any
non-empty open subsets $U,V$ of $X$, there is $r=r(U,V)>0$ such that
$T_t(U)\cap V\neq\varnothing$ provided $|t|>r$. We also say that
$\{T_t\}_{t\in A}$ is {\it hypercyclic} (respectively, {\it
supercyclic}) if the family $\{T_t:t\in A\}$ (respectively,
$\{zT_t:z\in\K,\ t\in A\}$) is universal. $\{T_t\}_{t\in A}$ is said
to be {\it hereditarily hypercyclic} (respectively, {\it
hereditarily supercyclic}) if for any sequence $\{t_n\}_{n\in\Z_+}$
of elements of $A$ such that $|t_n|\to\infty$, the family
$\{T_{t_n}:n\in\Z_+\}$ (respectively, $\{zT_{t_n}:z\in \K,\
n\in\Z_+\}$) is universal. A continuous linear operator $T$ acting
on a topological vector space $X$ is called {\it hypercyclic,
supercyclic, hereditarily hypercyclic, hereditarily supercyclic,
mixing {\rm or} transitive} if the semigroup $\{T^n\}_{n\in\Z_+}$
has the same property. It is worth noting that our definition of a
hereditarily hypercyclic operator follows Ansari \cite{ansa}, while
in the terminology of \cite{bp,gri}, the same property is called
'hereditarily hypercyclic with respect to the sequence $n_k=k$ of
all non-negative integers'. Hypercyclic and supercyclic operators
have been intensely studied during last few decades, see surveys
\cite{ge1,ge2,msa} and references therein. Clearly mixing implies
transitivity and hereditary hypercyclicity (respectively, hereditary
supercyclicity) implies hypercyclicity (respectively,
supercyclicity). Recall that a topological space $X$ is called a
{\it Baire space} if the intersection of countably many dense open
subsets of $X$ is dense in $X$. According to the classical Baire
theorem, complete metric spaces are Baire.

\begin{proposition} \label{untr2} Let $X$ be a
topological vector space, $A$ be a normed semigroup and ${\cal
S}=\{T_a\}_{a\in A}$ be an operator semigroup on $X$. Then
\begin{itemize}\itemsep=-2pt
\item[\rm(\ref{untr2}.1)]if ${\cal S}$ is hereditarily hypercyclic, then ${\cal S}$ is
mixing.
\end{itemize}
If additionally $X$ is Baire separable and metrizable, the converse
implication holds$:$
\begin{itemize}\itemsep=-2pt
\item[\rm(\ref{untr2}.2)]if ${\cal S}$ is mixing, then ${\cal S}$ is hereditarily hypercyclic.
\end{itemize}
\end{proposition}

The above proposition is a combination of well-known facts,
appearing in the literature in various modifications. It is worth
noting that a similar statement holds for hypercyclicity and
transitivity under certain natural additional assumptions. One can
also write down and prove a supercyclicity analogue of the above
proposition. In the next section we shall prove
Proposition~\ref{untr2} for sake of completeness. It is worth noting
that for any subsemigroup $A_0$ of $A$, not lying in the kernel of
the norm, $\{T_t\}_{t\in A_0}$ is mixing if $\{T_t\}_{t\in A}$ is
mixing. In particular, if $\{T_t\}_{t\in A}$ is mixing, then $T_t$
is mixing whenever $|t|>0$.

\subsection{Results}

The question of existence of supercyclic or hypercyclic operators or
semigroups on various types of topological vector spaces was
intensely studied. There are no hypercyclic operators on any finite
dimensional topological vector space and there are no supercyclic
operators on a finite dimensional topological vector space of real
dimension $>2$. These facts follow, for instance from the main
result of \cite{ww}. Herzog \cite{ger} demonstrated that there is a
supercyclic operator on any separable infinite dimensional Banach
space. Later Ansari \cite{ansa1} and Bernal-Gonz\'ales
\cite{bernal}, answering a question raised by Herrero, showed
independently that for any separable infinite dimensional Banach
space $X$ there is a hypercyclic operator $T\in L(X)$. Using the
same idea as in \cite{ansa1}, Bonet and Peris \cite{bonper} proved
that there is a hypercyclic operator on any separable infinite
dimensional Fr\'echet space and demonstrated that there is a
hypercyclic operator on an inductive limit $X$ of a sequence $X_n$
for $n\in\Z_+$ of separable Banach spaces provided there is
$n\in\Z_+$ for which $X_n$ is dense in $X$. Grivaux \cite{gri}
observed that hypercyclic operators $T$ provided in
\cite{ansa1,bernal,bonper} are in fact mixing and therefore
hereditarily hypercyclic. They actually come from the same source.
Namely, according to Salas \cite{sal} an operator of the shape
$I+T$, where $T$ is a backward weighted shift on $\ell_1$, is
hypercyclic. Virtually the same proof as in \cite{sal} demonstrates
that these operators are in fact mixing.  Moreover, all operators
constructed in the above papers, except for the ones acting on
$\omega$, are hypercyclic because of a quasisimilarity with one of
the operators of the shape identity plus a backward weighted shift.
The same quasisimilarity transfers the mixing property as
effectively as it transfers hypercyclicity. A similar idea was used
by Berm\'udez, Bonilla and Martin\'on \cite{ex1} and
Bernal-Gonz\'alez and Grosse-Erdmann \cite{ex2}, who have
demonstrated that any separable infinite dimensional Banach space
supports a hypercyclic strongly continuous semigroup $\{T_t\}_{t\in
\R_+}$. Berm\'udez, Bonilla, Conejero and Peris \cite{holo} have
shown for any separable infinite dimensional complex Banach space
$X$, there exists a mixing strongly continuous semigroup
$\{T_t\}_{t\in \Pi}$ with $\Pi=\{z\in\C:{\tt Re}\,z\geq 0\}$ such
that the map $t\mapsto T_t$ is holomorphic on the interior of $\Pi$.
As a matter of fact, one can easily see that the semigroup
constructed in \cite{holo} extends to a holomorphic mixing group
$\{T_t\}_{t\in\C}$. Finally, Conejero \cite{cone} proved that any
separable infinite dimensional Fr\'echet space non-isomorphic to
$\omega$ supports a uniformly continuous mixing operator semigroup
$\{T_t\}_{t\in\R_+}$.

The following theorem extracts the maximum of the method both in
terms of the class of spaces and semigroups. Although the general
idea remains the same, the proof requires dealing with a number of
technical details of various nature. For brevity we shall introduce
the following class of locally convex spaces.

\begin{definition}\label{MMM}We say that a sequence $\{x_n\}_{n\in\Z_+}$
of elements of a topological vector space $X$ is an $\ell_1$-{\it
sequence} if $x_n\to 0$ in $X$ and the series
$\sum\limits_{n=0}^\infty a_nx_n$ converges in $X$ for each
$a\in\ell_1$.

We say that a locally convex space $X$ belongs to the class $\M$ if
its topology is not weak and there exists an $\ell_1$-sequence in
$X$ with dense span.
\end{definition}

\begin{theorem}\label{saan} Let $X\in\M$. Then for any $k\in\N$,
there exists a hereditarily hypercyclic $(\!$and therefore mixing$)$
uniformly continuous operator group $\{T_t\}_{t\in\K^k}$. Moreover,
if $\,\K=\C$ the map $z\mapsto T_zx$ from $\C^k$ to $X$ is
holomorphic for each $x\in X$.
\end{theorem}

Since for any hereditarily hypercyclic semigroup
$\{T_t\}_{t\in\K^k}$ and any non-zero $t\in\K^k$, the operator $T_t$
is hereditarily hypercyclic, we have the following corollary.

\begin{corollary}\label{saan1} Let $X\in\M$. Then there is a
hereditarily hypercyclic $(\!$and therefore mixing$)$ operator $T\in
L(X)$.
\end{corollary}

\begin{remark}\label{Rem1} \ It is easy to see that if $X\in\M$, then
$X_\tau\in\M$. Indeed, if $\{x_n\}_{n\in\Z_+}$ is an
$\ell_1$-sequence in $X$ with dense span, then
$\{2^{-n}x_n\}_{n\in\Z_+}$ is an $\ell_1$-sequence in $X_\tau$ with
dense span. Of course, $X_\sigma$ never belongs to $\M$. On the
other hand, it is well-known that $L(X_\sigma)=L(X_\tau)$. Moreover,
since $\sigma(X,X')\subseteq\tau(X,X')$, then any strongly
continuous hereditarily hypercyclic operator semigroup
$\{T_t\}_{t\in\K^k}$ on $X_\tau$ is also strongly continuous and
hereditarily hypercyclic as an operator semigroup on $X_\sigma$.
Thus in the case $X_\tau\in\M$, Theorem~\ref{saan} implies that
there is a strongly continuous hereditarily hypercyclic operator
semigroup $\{T_t\}_{t\in\K^k}$ on $X_\sigma$. Unfortunately, the
nature of the weak topology does not allow to make such a semigroup
uniformly continuous.
\end{remark}

It is worth noting that any separable Fr\'echet space admits an
$\ell_1$-sequence with dense span. It is also well-known
\cite{shifer} that any Fr\'echet space carries the Mackey topology
and the topology on a Fr\'echet space $X$ differs from the weak
topology if and only if $X$ is infinite dimensional and is
non-isomorphic to $\omega$. That is, any separable infinite
dimensional Fr\'echet space non-isomorphic to $\omega$ belongs to
$\M$. Similarly, one can verify that an infinite dimensional
inductive limit $X$ of a sequence $X_n$ for $n\in\Z_+$ of separable
Banach spaces belongs to $\M$ provided there is $n\in\Z_+$ for which
$X_n$ is dense in $X$. Thus all the above mentioned existence
theorems are particular cases of Theorem~\ref{saan}.

Grivaux \cite{gri} raised a question whether each separable infinite
dimensional Banach space supports a hypercyclic non-mixing operator.
Since the class $\M$ contains separable infinite dimensional Banach
spaces, the following theorem provides an affirmative answer to this
question.

\begin{theorem}\label{nonmix} Let $X\in\M$.
Then there exists $T\in L(X)$ such that $T$ is hypercyclic and
non-mixing.
\end{theorem}

The simplest separable infinite dimensional locally convex space
space (and the only Fr\'echet space) outside $\M$ is $\omega$.
Curiously, the situation with $\omega$ is totally different.
Hypercyclic operators on the complex space $\omega$ have been
characterized by Herzog and Lemmert \cite{gele}. Namely, they proved
that a continuous linear operator $T$ on the complex Fr\'echet space
$\omega$ is hypercyclic if and only if the point spectrum
$\sigma_p(T')$ of $T'$ is empty. It also worth mentioning that B\'es
and Conejero \cite{beco} provided sufficient conditions for $T\in
L(\omega)$ to have an infinite dimensional closed linear subspace,
each non-zero vector of which is hypercyclic, and found common
hypercyclic vectors for some families of hypercyclic operators on
$\omega$. See also the related work \cite{peter1} by Petersson. The
following theorem extends the result of Herzog and Lemmert and
highlights the difference between $\omega$ and other Fr\'echet
spaces.

\begin{theorem}\label{omeg}
Let $T\in L(\omega)$ be such that $T'$ has no non-trivial finite
dimensional invariant subspaces and $\{p_l\}_{l\in\Z_+}$ be a
sequence of polynomials such that $\deg p_l\to\infty$ as
$l\to\infty$. Then the sequence $\{p_l(T):l\in\Z_+\}$ is universal.
Moreover, there is no strongly continuous supercyclic semigroup
$\{T_t\}_{t\in\R_+}$ on $\omega$.
\end{theorem}

Note that in the case $\K=\C$, $T'$ has no non-trivial finite
dimensional invariant subspaces if and only if
$\sigma_p(T')=\varnothing$. The first part of the above theorem
implies that any hypercyclic operator on $\omega$ is mixing. We
shall, in fact, verify the following more general statement.

\begin{theorem}\label{omegA} Let $X$ be a locally convex space
carrying weak topology and $T\in L(X)$. Then the following
conditions are equivalent
\begin{itemize}\itemsep=-2pt
\item[\rm(\ref{omegA}.1)]$T'$ has no non-trivial finite dimensional invariant subspaces$;$
\item[\rm(\ref{omegA}.2)]$T$ is transitive$;$
\item[\rm(\ref{omegA}.3)]$T$ is mixing$;$
\item[\rm(\ref{omegA}.4)]the semigroup $\{p(T)\}_{p\in\pp^*}$ is mixing, where
$\pp^*=\K[z]\setminus\{0\}$ is the multiplicative semigroup of
non-zero polynomials with the norm $|p|=\deg p$.
\end{itemize}
\end{theorem}

\begin{remark}\label{Rem2} Chan and Sanders \cite{chan} observed
that on the space $(\ell_2)_\sigma$, being the Hilbert space
$\ell_2$ with the weak topology, there is a transitive
non-hypercyclic operator. Theorem~\ref{omegA} provides a huge supply
of such operators. For instance, the backward shift $T$ on $\ell_2$
is mixing on $(\ell_2)_\sigma$ ($T'$ has no non-trivial finite
dimensional invariant subspaces) and $T$ is clearly non-hypercyclic
(each its orbit is bounded).
\end{remark}

Theorems~\ref{saan}, \ref{nonmix} and \ref{omeg} imply the following
curious corollary.

\begin{corollary}\label{COR} Let $X$ be a separable infinite
dimensional Fr\'echet space. Then the following are equivalent
\begin{itemize}\itemsep=-2pt
\item[\rm(\ref{COR}.1)]there is a hypercyclic non-mixing operator $T\in
L(X);$
\item[\rm(\ref{COR}.2)]there is a mixing uniformly continuous semigroup
$\{T_t\}_{t\in \R_+}$ on $X;$
\item[\rm(\ref{COR}.3)]there is a supercyclic strongly continuous semigroup
$\{T_t\}_{t\in \R_+}$ on $X;$
\item[\rm(\ref{COR}.4)]$X$ is non-isomorphic to $\omega$.
\end{itemize}
\end{corollary}

Another simple space outside $\M$ is $\phi$. Bonet and Peris
\cite{bonper} observed that there are no supercyclic operators on
$\phi$. On the other hand, Bonet, Frerick, Peris and Wengenroth
\cite{fre} constructed a hypercyclic operator on the locally convex
direct sum $X=\bigoplus\limits_{n=0}^\infty \ell_1$ of countably
many copies of the Banach space $\ell_1$. The space $X$ is clearly
an \lbs-space, is complete and non-metrizable. It is also easy to
see that $X\notin\M$ (there are no $\ell_1$-sequences in $X$ with
dense span). We find sufficient conditions of existence and of
non-existence of a hypercyclic operator on a locally convex space.
These conditions allow us to characterize the \lbs-spaces, which
admit a hypercyclic operator.

\begin{theorem}\label{lbs} Let $X$ be an \lbss-space, then the following conditions are equivalent:
\begin{itemize}\itemsep=-2pt
\item[{\rm (\ref{lbs}.1)}]$X$ admits no hypercyclic operator$;$
\item[{\rm (\ref{lbs}.2)}]$X$ admits no cyclic operator with dense
range$;$
\item[{\rm (\ref{lbs}.3)}]$X$ is isomorphic to $Y\times\phi$, where $Y$
is the inductive limit of a sequence $\{Y_n\}_{n\in\N}$ of separable
Banach spaces such that $Y_0$ is dense in $Y$.
\end{itemize}
\end{theorem}

The proof is based upon the following result, which is of
independent interest.

\begin{theorem}\label{timesphi} Let $X$ be a topological vector
space, which has no quotients isomorphic to $\phi$. Then there is no
cyclic operator with dense range on $X\times \phi$.
\end{theorem}

The following theorem provides another generalization of the
mentioned result of Bonet, Frerick, Peris and Wengenroth.

\begin{theorem}\label{sumfre} Let $\{X_n\}_{n\in\Z_+}$ be a sequence of
separable Fr\'echet spaces. Then there is a hypercyclic operator on
$X=\bigoplus\limits_{n=0}^\infty X_n$ if and only if the set
$\{n\in\Z_+:\text{$X_n$ is infinite dimensional}\}$ is infinite.
\end{theorem}

We derive the above theorem from the following result, concerning
more general spaces.

\begin{theorem}\label{sumM} Let $X_n\in\M$ for each $n\in\Z_+$ and
$X=\bigoplus\limits_{n=0}^\infty X_n$. Then there is a hypercyclic
operator on $X$.
\end{theorem}

The next issue, we discuss, are dual hypercyclic operators. Let $X$
be a locally convex space. Recall that $X'_\beta$ is the dual space
$X'$ endowed with the strong topology $\beta(X',X)$. It is worth
noting that if $X$ is a normed space, then the strong topology on
$X'$ coincides with the standard norm topology. Salas \cite{sal2}
has constructed an example of a hypercyclic operator $T$ on $\ell_2$
such that both $T$ and $T'$ are hypercyclic. This result motivated
Petersson \cite{peter} to introduce the following definition. We say
that a continuous linear operator $T$ on a locally convex space $X$
is {\it dual hypercyclic} if both $T$ and $T'$ are hypercyclic on
$X$ and $X'_\beta$ respectively. Using the construction of Salas,
Petersson proved that any infinite dimensional Banach space $X$ with
a monotonic and symmetric Schauder basis and with separable dual
admits a dual hypercyclic operator. He also raised the following
questions. Does there exist a dual hypercyclic operator on any
infinite dimensional Banach space with separable dual? Does there
exist a non-normable Fr\'echet space that admits a dual hypercyclic
operator? The first of these questions was recently answered
affirmatively by Salas \cite{saldual}. The following theorem
provides a sufficient condition for existence of a dual hypercyclic
operator on a locally convex space.

\begin{theorem}\label{dual} Let $X$ be an infinite dimensional
locally convex space admitting an $\ell_1$-sequence with dense span.
Assume also that there is an $\ell_1$-sequence $\{f_n\}_{n\in\Z_+}$
with dense span in $X'_\beta$ and at least one of the following
conditions is satisfied:
\begin{itemize}\itemsep=-2pt
\item[\rm(\ref{dual}.1)]the topology of $X$ coincides with
$\sigma(X,X');$
\item[\rm(\ref{dual}.2)]the topology of $X$ coincides with
$\tau(X,X');$
\item[\rm(\ref{dual}.3)]the set $\{f_n:n\in\Z_+\}$ is uniformly
equicontinuous.
\end{itemize}
Then $X$ admits a dual hypercyclic operator.
\end{theorem}

Since every separable Fr\'echet space admits an $\ell_1$-sequence
with dense span and every Fr\'echet space carries the Mackey
topology, the above theorem implies the following corollary.

\begin{corollary} \label{dual1}
Let $X$ be a separable infinite dimensional Fr\'echet space, such
that there is an $\ell_1$-sequence with dense span in $X'_\beta$.
Then there exists a dual hypercyclic operator $T\in L(X)$.
\end{corollary}

If $X$ is a Banach space, $X'_\beta$ is also a Banach space and
therefore has an $\ell_1$-sequence with dense span if and only if it
is separable. Thus Corollary~\ref{dual1} implies the next corollary,
which is the mentioned recent result of Salas.

\begin{corollary} \label{dual2}
Let $X$ be an infinite dimensional Banach space with separable dual.
Then there exists a dual hypercyclic operator $T\in L(X)$.
\end{corollary}

Corollary~\ref{dual1} also provides plenty of non-normable Fr\'echet
spaces admitting a dual hypercyclic operator, thus answering the
second of the above questions of Petersson. For instance, take the
complex Fr\'echet space $X$ of entire functions on one variable with
the topology of uniform convergence on compact sets. It is easy to
verify that the sequence of functionals $g_n(f)=(n!)^{-1}f^{(n)}(0)$
is an $\ell_1$-sequence with dense span in $X'_\beta$. Since $X$ is
also infinite dimensional and separable, Corollary~\ref{dual1}
implies that $X$ supports a dual hypercyclic operator.

The proofs of the above results are based upon the two main
ingredients. One of them are sufficient conditions of mixing and the
other is a criterion for a generic (in the Baire category sense)
operator from a given class to be hypercyclic. Our sufficient
conditions of mixing extend the result of Salas \cite{sal} on
hypercyclicity of perturbations of the identity by adding a backward
weighted shift. Apart from providing us with tools, these extensions
are of independent interest.

\begin{theorem}\label{main} Let $X$ be a topological vector space and
$T\in L(X)$ be such that the space
\begin{equation}\label{OmegaT}
\Lambda(T)=\spann\left(\bigcup_{n\in\N,\ |z|=1} ((T-zI)^n(X)\cap\ker
(T-zI)^{n})\right)
\end{equation}
is dense in $X$. Then $T$ is mixing. If additionally, $X$ is Baire,
separable and metrizable, then $T$ is hereditarily hypercyclic.
\end{theorem}

We shall see that the above theorem implies not only the mentioned
result of Salas, but also is applicable in many other situations.
For instance, we use the above theorem to prove the following
results.

\begin{theorem} \label{cococo} Let $X$ be a separable Banach
space and $\NNN$ be the operator norm closure in $L(X)$ of the set
of finite rank nilpotent operators. Then the set of $T\in\NNN$ for
which $T$ is supercyclic and $I+T$ is hypercyclic is a dense
$G_\delta$ subset of the complete metric space $\NNN$. If
additionally $X'$ is separable, then the set of $T\in\NNN$ for which
$T$ and $T'$ are supercyclic and $I+T$ and $I+T'$ are hypercyclic is
a dense $G_\delta$ subset of $\NNN$.
\end{theorem}

Note that if a Banach space $X$ has the approximation property
\cite{linden}, then  the set $\NNN$ from the above corollary is
exactly the set of compact quasinilpotent operators (in the case
$\K=\R$ by quasinilpotency of $T$ we mean quasinilpotency of the
complexification of $T$ or equivalently that $\|T^n\|^{1/n}\to 0$).
Thus we have the following corollary.

\begin{corollary} \label{cococo1} Let $X$ be a separable Banach
space with the approximation property and $\NNN\subset L(X)$ be the
set of compact quasinilpotent operators. Then the set of $T\in\NNN$
for which $T$ is supercyclic and $I+T$ is hypercyclic is a dense
$G_\delta$ subset of the complete metric space $\NNN$. If
additionally $X'$ is separable, then the set of $T\in\NNN$ for which
$T$ and $T'$ are supercyclic and $I+T$ and $I+T'$ are hypercyclic is
a dense $G_\delta$ subset of $\NNN$.
\end{corollary}

\begin{theorem} \label{cocococo} Let $X$ be a separable Banach
space and $\NNN$ be the set of nuclear quasinilpotent nilpotent
operators endowed with the nuclear norm metric. Then the set of
$T\in\NNN$ for which $T$ is supercyclic and $I+T$ is hypercyclic is
a dense $G_\delta$ subset of the complete metric space $\NNN$. If
additionally $X'$ is separable, then the set of $T\in\NNN$ for which
$T$ and $T'$ are supercyclic and $I+T$ and $I+T'$ are hypercyclic is
a dense $G_\delta$ subset of $\NNN$.
\end{theorem}

Theorems~\ref{cocococo} and~\ref{cococo} provide a large supply of
dual hypercyclic operators $T$ on any infinite dimensional Banach
space with separable dual.

\section{Extended backward shifts \label{ebsh}}

Godefroy and Shapiro \cite{gs} have introduced the notion of a
generalized backward shift. Namely, a continuous linear operator $T$
on a topological vector space $X$ is called a {\it generalized
backward shift} if its {\it generalized kernel}
$$
\Ker T=\bigcup_{n=1}^\infty \ker T^n
$$
is dense in $X$ and $\ker T$ is one-dimensional. We introduce a more
general concept. Namely, we say that $T$ is an {\it extended
backward shift} if
\begin{equation}\label{KER}
\KER T=\spann\left(\bigcup_{n=1}^\infty (T^n(X)\cap \ker
T^{n})\!\right)
\end{equation}
is dense in $X$. From the easy dimension argument \cite{gs} it
follows that if $T\in L(X)$ is a generalized backward shift, then
$\dim\ker T^n=n$ and $T(\ker T^{n+1})=\ker T^n$ for each $n\in\N$.
Hence $\ker T^n=T^n(\ker T^{2n})$ and therefore $\ker T^n=T^n(X)\cap
\ker T^{n}$ for any $n\in\N$. It follows that $\Ker T=\KER T$ for a
generalized backward shift. That is, any generalized backward shift
is an extended backward shift.

We also consider the following analog of the concept of an extended
backward shift for a $k$-tuple of operators. Let $T_1,\dots,T_k$ be
continuous linear operators on a topological vector space $X$. We
say that $T=(T_1,\dots,T_k)\in L(X)^k$ is a {\it EBS$_k$-tuple} if
$T_mT_j=T_jT_m$  for any $j,m\in\{1,\dots,k\}$ and
\begin{equation}\label{KERk}
\KER (T)=\spann\left(\bigcup_{n\in\N^k} \kappa(n,T)\!\!\right)\!,\ \
\text{where}\ \ \kappa(n,T)=T_1^{n_1}\!\dots
T_k^{n_k}\biggl(\bigcap_{j=1}^k \ker T_j^{2n_j}\!\biggr),
\end{equation}
is dense in $X$.

It is easy to see that in the case of one operator (that is, $k=1$
and $T=T_1\in L(X)$), $\kappa(n,T)=T^n(\ker T^{2n})=\ker T^n\cap
T^n(X)$ and therefore the last definition is a generalization of the
previous one. In order to study extended backward shifts we need to
establish some properties of the backward shift on the finite
dimensional space $\K^{2n}$.

\subsection{Backward shift on $\K^{2n}$ \label{tbs}}

The following lemma is a modification of a lemma from \cite{mrs}.

\begin{lemma}\label{detan}
For each $n\in\N$ and $z\in\C\setminus\{0\}$, the matrix
$A_{n,z}=\left\{\frac{z^{j+k-1}}{(j+k-1)!}\right\}_{j,k=1}^n$ is
invertible.
\end{lemma}

\begin{proof} For each $n,k\in\N$ consider the matrix
$M_{n,k}=\Bigl\{\frac{(k+n-l)!}{(k+n-l+j-1)!}\Bigr\}_{j,l=1}^n$.
First, we demonstrate that the determinants of $M_{n,k}$ satisfy the
recurrent formula
\begin{equation}
\det M_{n,k}=\frac{(n-1)!k!(k+1)!}{(k+n-1)!(k+n)!}\,\det
M_{n-1,k+2}\ \ \text{for $n\geq 2$.} \label{dere}
\end{equation}
The equality (\ref{dere}) for $n=2$ is trivial. Suppose now that
$n\geq 3$. Subtracting the previous column from each column of
$M_{n,k}$ except the first one, we see that $\det M_{n,k}=\det
N_{n,k}$ where $N_{n,k}=\Bigl\{j\frac{(k+n-l-1)!}{(k+n-l+j)!}
\Bigr\}_{j,l=1}^{n-1}$. Dividing the $j$-th row of $N_{n,k}$ by $j$
and multiplying the $j$-th column by $a_j=(k+n-j+1)(k+n-j)$ for
$1\leq j\leq n$, we arrive at the matrix $M_{n-1,k+2}$. Hence
$$
\det M_{n,k}=\det M_{n-1,k+2}\prod_{j=1}^{n-1} ja_j^{-1}
=\frac{(n-1)!k!(k+1)!}{(k+n-1)!(k+n)!} \,\det M_{n-1,k+2},
$$
which proves (\ref{dere}). Since $\det M_{1,k}=1$ for each $k\in\N$,
from (\ref{dere}), it follows that $\det M_{n,k}\neq 0$ for any
$n,k\in\N$. Let now $B_n$ be the matrix obtained from $A_{n,1}$ by
putting the columns of $A_{n,1}$ in the reverse order. Clearly $\det
A_{n,1}=(-1)^{n-1}\det B_n$. On the other hand, multiplying the
$j$-th column of $B_n$ by $(n-j+1)!$ for $1\leq j\leq n$, we get the
matrix $M_{n,1}$. Hence,
$$
\det A_{n,1}=(-1)^{n-1}\det M_{n,1} \prod\limits_{j=1}^n (j!)^{-1}\
\ \text{for each $n\in\N$.}
$$
Since $\det M_{n,1}\neq 0$, we see that $\det A_{n,1}\neq 0$ and
therefore $A_{n,1}$ is invertible. Finally, for any $z\in\C$
consider the diagonal $n\times n$ matrix $D_{n,z}$ with the entries
$(1,z,\dots,z^{n-1})$ on the main diagonal. It is straightforward to
verify that
\begin{equation}\label{dia}
A_{n,z}=zD_{n,z}A_{n,1}D_{n,z} \ \ \ \text{for any $z\in\C$}.
\end{equation}
Since $A_{n,1}$ and $D_{n,z}$ for $z\neq 0$ are invertible, we see
that $A_{n,z}$ is invertible for each $n\in\N$ and
$z\in\C\setminus\{0\}$.
\end{proof}

\begin{lemma}\label{jordan}
Let $n\in\N$ and $e_1,\dots,e_{2n}$ be the canonical basis of
$\K^{2n}$ and $S\in L(\K^{2n})$ be the backward shift defined by
$Se_1=0$ and $Se_k=e_{k-1}$ for $2\leq k\leq 2n$ and $P$ the linear
projection on $\K^{2n}$ onto the subspace
$E=\spann\{e_1,\dots,e_n\}$ along
$F=\spann\{e_{n+1},\dots,e_{2n}\}$. Then for any
$z\in\K\setminus\{0\}$ and $u,v\in E$, there exists a unique
$x^z=x^z(u,v)\in \K^{2n}$ such that
\begin{equation}\label{ab}
Px^z=u\quad\text{and}\quad Pe^{zS}x^z=v.
\end{equation}
Moreover, for any bounded subset $B$ of $E$ and any $\epsilon>0$,
there is $c=c(\epsilon,B)>0$ such that
\begin{align}
\sup_{u,v\in B}|(x^z(u,v))_{n+j}|&\leq c|z|^{-j}\quad\text{for
$1\leq j\leq n$ and $|z|\geq\epsilon;$} \label{Ox1}
\\
\sup_{u,v\in B}|(e^{zS}x^z(u,v))_{n+j}|&\leq
c|z|^{-j}\quad\quad\text{for $1\leq j\leq n$ and $|z|\geq\epsilon$.}
\label{Ox2}
\end{align}
In particular, $x^z(u,v)\to u$ and $e^{zS}x^z(u,v)\to v$ as $|z|\to
\infty$ uniformly for $u,v$ from any bounded subset of $E$.
\end{lemma}

\begin{proof} Let $u,v\in E$ and $z\in\K\setminus\{0\}$.
For $y\in \K^{2n}$ we denote
$$
\overline{y}=(y_{n+1},\dots,y_{2n})\in\K^n.
$$
One easily sees that (\ref{ab}) is equivalent to the vector equation
\begin{equation}
A_{n,z}\overline{x}^z=w^z, \label{vector}
\end{equation}
where $A_{n,z}$ is the matrix from Lemma~\ref{detan} and
$w^z=w^z(u,v)\in\K^n$ is defined as
\begin{equation} \label{wm}
w^z_j=v_{n-j+1}-\sum_{k=n-j+1}^n\frac{z^{k+j-n-1}u_k}{(k+j-n-1)!}
\quad\text{for $1\leq j\leq n$},
\end{equation}
provided we set $x_{j}=u_{j}$ for $1\leq j\leq n$. According to
Lemma \ref{detan}, the matrix $A_{n,z}$ is invertible for any
$z\in\K\setminus\{0\}$ and therefore (\ref{vector}) is uniquely
solvable. Thus there exists a unique $x^z=x(z,u,v)\in \K^{2n}$
satisfying (\ref{ab}). It remains to verify the estimates
(\ref{Ox1}) and (\ref{Ox2}). From (\ref{wm}) it follows that for any
bounded subset $B$ of $E$ and any $\epsilon>0$, there is
$a=a(\epsilon,B)$ such that
\begin{equation}
|(w^z(u,v))_j|\leq a|z|^{j-1}\quad\text{if $u,v\in B$,
$|z|\geq\epsilon$ and $1\leq j\leq n$}. \label{Ow}
\end{equation}
Recall that $D_{n,z}$ is the diagonal $n\times n$ matrix with the
entries $(1,z,\dots,z^{n-1})$ on the diagonal. Equalities
(\ref{vector}) and (\ref{dia}) imply
$$
\overline{x}^z=A_{n,z}^{-1}w^z=
z^{-1}D_{n,z}^{-1}A_{n,1}^{-1}D_{n,z}^{-1}w^z,
$$
where we use invertibility of $A_{n,1}$ provided by
Lemma~\ref{detan}. According to (\ref{Ow}), the set
$\{D_{n,z}^{-1}w^z(u,v):|z|\geq \epsilon,\ u,v\in B\}$ is bounded in
$\K^n$. Hence the set $Q=\{A_{n,1}^{-1}D_{n,z}^{-1}w^z(u,v):|z|\geq
\epsilon,\ u,v\in B\}$ is bounded in $\K^n$. From the last display
we see that
$$
(x^z(u,v))_{n+j}=\overline{x}^z_j\subseteq \{z^{-1}(D_{n,z}^{-1}
y)_j:y\in Q\}\ \ \text{if $|z|\geq \epsilon$,\ $u,v\in B$}.
$$
Boundedness of $Q$ implies now that (\ref{Ox1}) is satisfied with
some $c=c_1(\epsilon,B)$. Finally, since
$$
(e^{zS}x^z)_{n+j}=\sum_{l=n+j}^{2n} \frac{z^{l-n-j}x^z_l}{(l-n-j)!}
\quad\text{for}\ \ 1\leq j\leq n,
$$
there exists $c=c_2(\epsilon,B)$ for which (\ref{Ox2}) is satisfied.
Hence both (\ref{Ox2}) and (\ref{Ox1}) are satisfied with
$c=\max\{c_1(\epsilon,B),c_2(\epsilon,B)\}$.
\end{proof}

The next corollary follows immediately from Lemma~\ref{jordan}.

\begin{corollary}\label{jordan1}
Let $n\in\N$, $E\subseteq \K^{2n}$ and $S\in L(\K^{2n})$ be as in
Lemma~$\ref{jordan}$. Then for any $u,v\in E$ and any sequence
$\{z_j\}_{j\in\Z_+}$ in $\K$ satisfying $|z_j|\to\infty$, there
exists a sequence $\{x_j\}_{j\in\Z_+}$ of elements of $\K^{2n}$ such
that $x_j\to u$ and $e^{z_jS}x_j\to v$ as $j\to\infty$.
\end{corollary}

\begin{lemma}\label{jordan11}
Let $n\in\N$, $E\subseteq \K^{2n}$ and $S\in L(\K^{2n})$ be as in
Lemma~$\ref{jordan}$. Then for any bounded sequences
$\{u_j\}_{j\in\Z_+}$ and $\{v_j\}_{j\in\Z_+}$ of elements of $E$,
there exists a sequence $\{x_j\}_{j\in\Z_+}$ of elements of
$\K^{2n}$ such that $x_j-u_j\to 0$ and $(I+S)^jx_j-v_j\to 0$ as
$j\to\infty$.
\end{lemma}

\begin{proof}It is easy to see that there is $J\in L(\K^{2n})$ such
that $J$ has an upper triangular matrix, $J$ is invertible and
$S=J^{-1}(e^S-I)J$. Indeed, $S$ and $e^S-I$ are similar since they
are nilpotent of maximal rank $2n-1$. Moreover, since $S$ and
$e^S-I$ are upper triangular, the similarity operator can be chosen
upper triangular and therefore $J(E)\subseteq E$.

Let now $x_j=J^{-1}x^j(Ju_j,Jv_j)$, where $x^z(u,v)$ is defined in
Lemma~\ref{jordan}. Since the set
$$
\{Ju_j:j\in\Z_+\}\cup\{Jv_j:j\in\Z_+\}
$$
is bounded and is contained in $E$ because $J(E)\subseteq E$, from
Lemma~\ref{jordan} it follows that $x^j(Ju_j,Jv_j)-Ju_j\to 0$ and
$e^{jS}x^j(Ju_j,Jv_j)-Jv_j\to 0$ as $j\to\infty$. Multiplying by
$J^{-1}$, we obtain $x_j-u_j\to 0$ and
$J^{-1}e^{jS}Jx_j=(I+S)^jx_j-v_j\to 0$ as $j\to\infty$.
\end{proof}

In order to construct multi-parameter mixing semigroups, we need the
following multi-operator version of Corollary~\ref{jordan1}.

\begin{lemma}\label{jordan2} Let $k\in\N$, $n_1,\dots,n_k\in \N$,
for each $j\in\{1,\dots,k\}$ let $e^j_1,\dots,e^j_{2n_j}$ be the
canonical basis in $\K^{2n_j}$,
$E_j=\spann\{e^j_1,\dots,e^j_{n_j}\}$ and $S_j\in L(\K^{2n_j})$ be
the backward shift: $S_je^j_1=0$ and $S_je^j_l=e^j_{l-1}$ for $2\leq
l\leq 2n_j$. Let also $X=\K^{2n_1}\otimes {\dots}\otimes \K^{2n_k}$,
$E=E_1\otimes {\dots}\otimes E_k$ and for $1\leq j\leq k$,
$$
T_j\in L(X),\quad T_j=I\otimes{\dots}I\otimes S_j\otimes I\otimes
{\dots}\otimes I,
$$
where $S_j$ sits in the $j^{\rm th}$ place. Finally, let
$\{z_m\}_{m\in\Z_+}$ be a sequence of elements of $\K^k$ with
$|z_m|\to\infty$. Then for any $u,v\in E$, there exists a sequence
$\{x_m\}_{m\in\Z_+}$ of elements of $X$ such that $x_m\to u$ and
$e^{\langle z_m,T\rangle}x_m\to v$ as $m\to\infty$, where $\langle
s,T\rangle=s_1T_1+{\dots}+s_kT_k$.
\end{lemma}

\begin{proof} If the statement of the lemma is false, then there
are $u,v\in E$ and a subsequence $\{z'_{m}\}_{m\in\Z_+}$ of
$\{z_m\}$ such that $(u,v)$ does not belong to the closure of the
set $\{(x,e^{\langle z'_{l},T\rangle}x):x\in X,\ m\in\Z_+\}$. Let
$\overline{\K}=\K\cup\{\infty\}$ be the one-point compactification
of $\K$. Since $\overline{\K}^k$ is compact and metrizable, we can
pick a convergent in $\overline{\K}^k$ subsequence $\{z''_m\}$ of
$\{z'_{m}\}$. Clearly the statement of the lemma remains false with
$\{z_m\}$ replaced by $\{z''_m\}$. That is, it suffices to consider
the case when $\{z_m\}$ converges in $\overline{\K}^k$.

Thus without loss of generality, we can assume that $\{z_m\}$
converges to $w\in \overline{\K}^k$. Let $C=\{j:w_j=\infty\}$. Since
$|z_m|\to\infty$, the set $C$ is non-empty. Without loss of
generality, we may also assume that $C=\{1,\dots,r\}$ with $1\leq
r\leq k$.

Denote by $\Sigma$ the set of $(u,v)\in X\times X$ such that there
exists a sequence $\{x_m\}_{m\in\Z_+}$ of elements of $X$ for which
$x_m\to u$ and $e^{\langle z_m,T\rangle}x_m\to v$ as $m\to\infty$.
We have to demonstrate that $E\times E\subseteq \Sigma$. Let $u_j\in
E_j$ for $1\leq j\leq k$ and $u=u_1\otimes{\dots}\otimes u_k$. By
Corollary~\ref{jordan1}, for $1\leq j\leq r$, there exist sequences
$\{x_{j,m}\}_{m\in\Z_+}$ and $\{y_{j,m}\}_{m\in\Z_+}$ of elements of
$\K^{2n_j}$ such that
$$
\text{$x_{j,m}\to 0$, $e^{(z_m)_j S_j}x_{j,m}\to u_j$, $y_{j,m}\to
u_j$ and $e^{(z_m)_j S_j}y_{j,m}\to 0$ as $m\to\infty$ for $1\leq
j\leq r$}.
$$
Now we put $x_{j,m}=e^{-w_jS_j}u_j$ and $y_{j,m}=u_j$ for $r<j\leq
k$ and $m\in\Z_+$. Consider the sequences $\{x_m\}_{m\in\Z_+}$ and
$\{y_m\}_{m\in\Z_+}$ of elements of $X$ defined by the formula
$x_m=x_{1,m}\otimes{\dots}\otimes x_{k,m}$ and
$y_m=y_{1,m}\otimes{\dots}\otimes y_{k,m}$. According to the
definition of $x_m$ and $y_m$ and the above display, $x_m\to 0$ and
$y_m\to u$. Indeed, $x_m\to 0$ because for any $j$, the sequence
$x_{j,m}$ is bounded and $x_{1,m}\to 0$. Similarly, taking into
account that $(z_m)_j\to w_j$ for $j>r$, we see that $e^{\langle
z_m,T\rangle}x_m\to u$ and $e^{\langle z_m,T\rangle}y_m\to 0$. Hence
$(u,0)\in \Sigma$ and $(0,u)\in\Sigma$. Thus
$$
(\{0\}\times E_0)\cup(E_0\times \{0\})\subseteq \Sigma,
$$
where $E_0=\{u_1\otimes{\dots}\otimes u_k:u_j\in E_j,\ 1\leq j\leq
k\}$. On the other hand, it is easy to see that the linear span of
the set $(\{0\}\times E_0)\cup(E_0\times \{0\})$ is exactly $E\times
E$. Since $\Sigma$ is a linear space, the above display implies that
$E\times E\subseteq \Sigma$.
\end{proof}

For applications it is more convenient to reformulate the above
lemma in the coordinate form.

\begin{corollary}\label{jordan3} Let $k\in\N$, $n_1,\dots,n_k\in \N$,
for each $j\in\{1,\dots,k\}$ let $N_j=\{1,\dots,2n_j\}$,
$Q_j=\{1,\dots,n_j\}$. Consider the sets $M=N_1\times {\dots}\times
N_k$, $M_0=Q_1\times {\dots}\times Q_k$ and let $\{e_m:m\in M\}$ be
the canonical basis of the finite dimensional vector space $X=\K^M$.
Let $E=\spann\{e_m:m\in M_0\}$ and for $1\leq j\leq k$, $T_j\in
L(X)$ be the operator acting on the canonical basis in the following
way: $T_je_m=0$ if $m_j=1$, $T_je_m=e_{m'}$ if $m_j>1$, where
$m'_l=m_l$ if $l\neq j$, $m'_j=m_j-1$. Then for any sequence
$\{z_m\}_{m\in\Z_+}$ of elements of $\K^k$ with $|z_m|\to\infty$ and
any $u,v\in E$, there exists a sequence $\{x_m\}_{m\in\Z_+}$ of
elements of $X$ such that $x_m\to u$ and $e^{\langle
z_m,T\rangle}x_m\to v$ as $m\to\infty$.
\end{corollary}

\subsection{The key lemmas \label{key}}

\begin{lemma}\label{kerim0}
Let $X$ be a Hausdorff topological vector space, $k\in\N$,
$n=(n_1,\dots,n_k)\in\N^k$ and $A=(A_1,\dots,A_k)\in L(X)^k$ be such
that $A_jA_l=A_lA_j$ for any $l,j\in \{1,\dots,k\}$. Then for each
$x$ from the space $\kappa(n,A)$ defined by $(\ref{KERk})$, there
exists a common finite dimensional invariant subspace $Y$ for
$A_1,\dots,A_k$ such that for any sequence $\{z_m\}_{m\in\Z_+}$ of
elements of $\K^k$ with $|z_m|\to\infty$ there exist sequences
$\{x_m\}_{m\in\Z_+}$, $\{y_m\}_{m\in\Z_+}$ of elements of $Y$ for
which
\begin{equation}\label{converg0}
x_m\to0,\ \ e^{\langle z_m,A\rangle}x_m\to x,\ \ y_m\to x,\ \
e^{\langle z_m,A\rangle}y_m\to 0 \qquad \text{as \ $m\to\infty$},
\end{equation}
where $\langle s,A\rangle=(s_1A_1+{\dots}+s_kA_k)\bigr|_Y$.
\end{lemma}

\begin{proof}
Since $x\in \kappa(n,T)$, there exists $y\in X$ such that
$x=A^{n_1}\!\dots A^{n_k}y$ and $A^{2n_j}y=0$ for $1\leq j\leq k$.
For each $j\in\{1,\dots,k\}$ let $N_j=\{1,\dots,2n_j\}$ and
$Q_j=\{1,\dots,n_j\}$. Denote $M=N_1\times {\dots}\times N_{k}$,
$M_0=Q_1\times {\dots}\times Q_{k}$. For any $l\in M$, let
$h_l=A^{2n_1-l_1}\!\dots A^{2n_{k}-l_{k}}y$ and let
$Y=\spann\{h_l:l\in M\}$. Clearly $Y$ is finite dimensional. It is
also straightforward to verify that $A_j h_l=0$ if $l_j=1$ and $A_j
h_l=h_{l'}$ if $l_j>1$, where $l'_r=l_r$ for $r\neq j$,
$l'_j=l_j-1$, $1\leq j\leq k_0$. It follows that $Y$ is invariant
for $A_j$ for $1\leq j\leq k$. Consider the linear operator
$J:\K^M\to Y$ defined on the canonical basis by the formulas
$Je_l=h_l$ for $l\in M$. Let also $E=\spann\{e_l:l\in M_0\}$ and
$T_j\in L(\K^M)$ be the operators from Corollary~\ref{jordan3}.
Taking into account the definition of $T_j$ and the action of $A_j$
on $h_l$, we see that $A_jJ=JT_j$ for $1\leq j\leq k$. Clearly
$n=(n_1,\dots,n_{k})\in M_0$ and therefore $e_n\in E$. Since
$x=A^{n_1}\!\dots A^{n_k}y$, we also see that $x=h_n$. According to
Corollary~\ref{jordan3}, there exist sequences $\{u_m\}_{m\in\Z_+}$
and $\{v_m\}_{m\in\Z_+}$ of elements of $\K^M$ such that $u_m\to
e_n$, $e^{\langle z_m,T\rangle}u_m\to 0$, $v_m\to 0$ and $e^{\langle
z_m,T\rangle}u_m\to e_n$ as $m\to\infty$. Now let $y_m=Ju_m$ and
$x_m=Jv_m$ for $m\in\Z_+$. Then $\{x_m\}$ and $\{y_m\}$ are
sequences of elements of $Y$. From the intertwining relations
$A_jJ=JT_j$ and the fact that $\K^M$ and $Y$ are finite dimensional,
it follows that $x_m\to J0=0$, $y_m\to Je_n=x$, $e^{\langle
z_m,A\rangle}x_m\to Je_n=x$, $e^{\langle z_m,T\rangle}y_m\to J0=0$.
Thus (\ref{converg0}) is satisfied.
\end{proof}

Let $X$ be a topological vector space and $T=(T_1,\dots,T_k)\in
L(X)^k$. We write $T\in \E(k,X)$ if for any $z\in \K^k$, the series
$\sum\limits_{n=0}^\infty \frac{1}{n!} \langle z,T\rangle^n x$
converges in $X$ and the map $(z,x)\to e^{\langle z,T\rangle}x$ from
$\K^k\times X$ to $X$ is separately continuous, where $\langle
z,T\rangle=z_1T_1+{\dots}+z_kT_k$ and
$e^{S}x=\sum\limits_{n=0}^\infty \frac{1}{n!} S^n x$. The following
proposition is an elementary exercise. We leave its proof for the
reader.

\begin{proposition}\label{expo} Let $k\in\N$, $X$ be a locally
convex space and $T\in \E(k,X)$. Assume also that $T_jT_m=T_mT_j$
for any $m,j\in\{1,\dots,k\}$. Then $\{e^{\langle
z,T\rangle}\}_{z\in\K^k}$ is a strongly continuous operator group.
Moreover, if $\K=\C$, then the map $z\mapsto e^{\langle
z,T\rangle}x$ from $\C^k$ to $X$ is holomorphic for each $x\in X$.
\end{proposition}

\begin{remark}\label{commo} It is worth noting that the semigroup
property $e^{\langle z+w,T\rangle}=e^{\langle z,T\rangle}e^{\langle
w,T\rangle}$ fails if the operators $T_j$ are not pairwise
commuting.
\end{remark}

\begin{corollary}\label{kerim1}
Let $X$ be a locally convex space, $k\in\N$, $n_1,\dots,n_k\in\N$
and $A=(A_1,\dots,A_k)\in \E(k,X)$ be such that $A_jA_l=A_lA_j$ for
any $l,j\in \{1,\dots,k\}$. Then for each $x,y$ from the space
$\KER(A)$ defined by $(\ref{KERk})$ and any sequence
$\{z_m\}_{m\in\Z_+}$ of elements of $\K^k$ with $|z_m|\to\infty$,
there exist a sequence $\{u_m\}_{m\in\Z_+}$ in $X$ such that $u_m\to
x$ and $e^{\langle z_m,A\rangle}u_m\to y$ as $m\to\infty$.
\end{corollary}

\begin{proof} Fix a sequence $\{z_m\}_{m\in\Z_+}$ of
elements of $\K^k$ with $|z_m|\to\infty$. Let $\Sigma$ be the set of
$(x,y)\in X\times X$ for which there exists a sequence
$\{u_m\}_{n\in\Z_+}$ in $X$ such that $u_m\to x$ and $e^{\langle
z_m,A\rangle}u_m\to y$ as $m\to\infty$. According to
Lemma~\ref{kerim0}, $\kappa(n,A)\times \{0\}\subseteq \Sigma$  and
$\{0\}\times \kappa(n,A)\subseteq \Sigma$ for any
$n=(n_1,\dots,n_k)\in\N^k$, where the space $\kappa(n,A)$ is defined
by (\ref{KERk}). On the other hand, from the definition of $\Sigma$
it is clear that $\Sigma$ is a linear subspace of $X\times X$. Thus
$$
\KER(A)\times\KER(A)=\spann\left(\bigcup_{n\in\N^k}
\bigl((\kappa(n,A)\times \{0\})\cup (\{0\}\times
\kappa(n,A))\bigr)\right)\subseteq\Sigma.
$$
Hence $(x,y)\in\Sigma$ for any $x,y\in \KER(A)$.
\end{proof}

\begin{lemma}\label{kerim}
Let $X$ be a topological vector space, $z\in\K$, $|z|=1$, $A\in
L(X)$, $m\in\N$ and $x\in A^m(X)\cap \ker A^m$. There exist
sequences $\{u_k\}_{k\in\Z_+}$ and $\{v_k\}_{k\in\Z_+}$ of elements
of $X$ such that
\begin{equation}\label{converg}
u_k\to0,\ \ z^k(I+A)^ku_k\to x,\ \ v_k\to x,\ \ z^k(I+A)^kv_k\to 0
\qquad \text{as $k\to\infty$}.
\end{equation}
\end{lemma}

\begin{proof} If $x=0$, we can take $u_k=v_k=0$, so we may assume that
$x\neq0$. Let $n$ be the smallest positive integer for which
$A^nx=0$. Since $A^mx=0$, we have $n\leq m$. Hence $x\in
A^m(X)\subseteq A^n(X)$. Thus we can pick $w\in X$ such that
$A^nw=x$. Denote
$$
h_{j}=A^{2n-j}w\quad \text{for}\quad 1\leq j\leq 2n \quad
\text{and}\quad  Y=\spann\{h_1,\dots,h_{2n}\}.
$$
Clearly $Ah_j=h_{j-1}$ for $2\leq j\leq 2n$, $h_n=A^nh_{2n}=A^nw=x$
and $Ah_1=A^{2n}h_{2n}=A^nx=0$. By definition of $n$ we have
$h_1=A^{2n-1}h_{2n}=A^{n-1}x\neq 0$. In particular $A(Y)\subseteq
Y$, $(A\bigr|_{Y})^{2n}=0$ and $(A\bigr|_Y)^{2n-1}\neq 0$. Since the
order of nilpotency of a nilpotent operator on $Y$ cannot exceed the
dimension of $Y$, we have $\dim Y\geq 2n$. On the other hand $Y$ is
the span of the $2n$-elements set $\{h_1,\dots,h_{2n}\}$. Hence
$\{h_1,\dots,h_{2n}\}$ is a linear basis of $Y$. Thus there exists a
unique linear isomorphism $J:\K^{2n}\to Y$ such that $Je_k=h_k$ for
$1\leq j\leq 2n$. Since $Ah_j=h_{j-1}$ for $2\leq j\leq 2n$ and
$Ah_1=0$, we have $A\bigr|_Y=JSJ^{-1}$, where $S\in L(\K^{2n})$ is
the backward shift operator from Lemma~\ref{jordan}. Applying
Lemma~\ref{jordan11} with $u^m=(0,\dots,0,1)$ and
$v^m=(0,\dots,0,0)$, we find that there exists a sequence
$\{g_k\}_{k\in\Z_+}$ of vectors in $\K^{2n}$ such that $g_k\to e_n$
and $(I+S)^kg_k\to 0$ as $k\to\infty$. Applying Lemma~\ref{jordan11}
with $u^m=(0,\dots,0,0)$ and $v^m=(0,\dots,0,z^{-m})$,  there exists
also a sequence $\{f_k\}_{k\in\Z_+}$ of vectors in $\K^{2n}$ such
that $f_k\to 0$ and $z^k(I+S)^kf_k\to e_n$ as $k\to\infty$. Define
now $u_k=Jf_k$ and $v_k=Jg_k$ for $k\in\Z_+$. Since $Y$ and
$\K^{2n}$ are finite dimensional, we see that
\begin{align*}
&v_k\to Je_n=x,\ \  u_k\to 0,\ \
(I+A)^kv_k=J(I+S)^kJ^{-1}Jg_k=J(I+S)^kg_k\to 0,
\\
&z^k(I+A)^ku_k=z^kJ(I+S)^kJ^{-1}Jf_k=z^kJ(I+S)^kf_k\to Je_n=x\ \
\text{as $k\to\infty$.}
\end{align*}
Thus the sequences $\{u_k\}$ and $\{v_k\}$ satisfy the desired
conditions.
\end{proof}

The following corollary of Lemma~\ref{kerim} seems to be of
independent interest.

\begin{theorem}\label{omt}
Let $X$ be a topological vector space, $T\in L(X)$ and $\Lambda(T)$
be the set defined in $(\ref{OmegaT})$. Then for any
$u,v\in\Lambda(T)$, there exists a sequence $\{x_k\}_{k\in\Z_+}$ in
$X$ such that $x_k\to u$ and $T^kx_k\to v$ as $k\to\infty$.
\end{theorem}

\begin{proof}
Let $\Sigma$ be the set of pairs $(x,y)\in X\times X$ for which
there exists a sequence $\{x_k\}_{k\in\Z_+}$ of elements of $X$ such
that $x_k\to x$ and $T^kx_k\to y$. Pick $z\in\K$ with $|z|=1$ and
$n\in\N$ and let
$$
x\in\ker(T-zI)^n\cap (T-zI)^n(X)=\ker A^n\cap A^n(X),\ \
\text{where}\ \ A=z^{-1}T-I.
$$
By Lemma \ref{kerim}, there exist sequences $\{u_k\}_{k\in\Z_+}$ and
$\{v_k\}_{k\in\Z_+}$ in $X$ satisfying (\ref{converg}). Since
$I+A=z^{-1}T$, (\ref{converg}) can be rewritten in the following
way:
$$
v_k\to x,\ \ T^kv_k\to 0,\ \ u_k\to 0,\ \ T^ku_k\to x\quad\text{as\
\ $k\to\infty$}.
$$
This shows that for any $n\in\N$ and any $z\in\K$ with $|z|=1$,
$$
\bigl(\ker(T-zI)^n\cap (T-zI)^n(X)\bigr)\times \{0\}\subseteq
\Sigma\quad\text{and}\quad \{0\}\times \bigl(\ker (T-zI)^n\cap
(T-zI)^n(X)\bigr)\subseteq \Sigma.
$$
On the other hand, the fact that $\Sigma$ is a linear subspace of
$X\times X$, the definition of $\Lambda(T)$ and the above display
imply that $\Lambda(T)\times \Lambda(T)\subseteq \Sigma$.
\end{proof}

\subsection{Mixing semigroups and extended backward shifts}

We start by proving Proposition~\ref{untr2}. The next observation is
Proposition~1 in \cite{ge1}.

\begin{thmd} Let $X$ be a topological space and ${\cal
F}=\{T_\alpha:\alpha\in A\}$ be a family of continuous maps from $X$
to $X$ such that $T_\alpha T_\beta=T_\beta T_\alpha$ for any
$\alpha,\beta\in A$ and $T_\alpha(X)$ is dense in $X$ for any
$\alpha\in A$. Then the set of universal elements for $\cal F$ is
either empty or dense in $X$.
\end{thmd}

The following general theorem can be found in
\cite[pp.~348--349]{ge1}.

\begin{thmdu} Let $X$ be a Baire topological space, $Y$ a second countable
topological space and $\{T_a:a\in A\}$ a family of continuous maps
from $X$ into $Y$. Then the following assertions are equivalent:

\noindent{\rm (U1)}\ \ The set of universal elements for $\{T_a:a\in
A\}$ is dense in $X;$

\noindent{\rm (U2)}\ \ The set of universal elements for $\{T_a:a\in
A\}$ is a dense $G_\delta$-subset of $X;$

\noindent{\rm (U3)}\ \ The set $\{(x,T_ax):x\in X,\ a\in A\}$ is
dense in $X\times Y.$
\end{thmdu}

\begin{proof}[Proof of Proposition~$\ref{untr2}$]
Assume that $\{T_t\}_{t\in A}$ is hereditarily hypercyclic. That is,
$\{T_{t_n}:n\in\Z_+\}$ is universal for any sequence
$\{t_n\}_{n\in\Z_+}$ of elements of $A$ such that $|t_n|\to\infty$.
Applying this property to the sequence $t_n=nt$ with $t\in A$,
$|t|>0$, we see that $T_t$ is hypercyclic. Since any hypercyclic
operator has dense range \cite{ge1}, we get that $T_t(X)$ is dense
in $X$ for any $t\in A$ with $|t|>0$. We proceed by reasoning ad
absurdum. Assume that $\{T_t\}_{t\in A}$ is non-mixing. Then there
exist non-empty open subsets $U$ and $V$ of $X$ and a sequence
$\{t_n\}_{n\in\Z_+}$ of elements of $A$ such that $|t_n|\to\infty$
and $|t_n|>0$, $T_{t_n}(U)\cap V=\varnothing$ for each $n\in\Z_+$.
Since each $T_{t_n}$ has dense range and $T_{t_n}$ commute with each
other, Proposition~G implies that the set $W$ of universal elements
of $\{T_{t_n}:n\in\Z_+\}$ is either empty or dense in $X$. Since
$\{T_{t_n}:n\in\Z_+\}$ is universal, $W$ is non-empty and therefore
dense in $X$. Hence we can pick $x\in W\cap U$. Since $x$ is
universal for $\{T_{t_n}:n\in\Z_+\}$, there is $n\in\Z_+$ for which
$T_{t_n}x\in V$. Hence $T_{t_n}x\in T_{t_n}(U)\cap V=\varnothing$.
This contradiction completes the proof of (\ref{untr2}.1).

Next, assume that $X$ is Baire separable and metrizable,
$\{T_t\}_{t\in A}$ is mixing and $\{t_n\}_{n\in\Z_+}$ is a sequence
of elements of $A$ such that $|t_n|\to\infty$. From the definition
of mixing it follows that for any non-empty open subsets $U$ and $V$
of $X$, $T_{t_n}(U)\cap V\neq\varnothing$ for all sufficiently large
$n\in\Z_+$. Hence $\{(x,T_{t_n}x):x\in X,\ n\in\Z_+\}$ is dense in
$X\times X$. By Theorem~U, $\{T_{t_n}:n\in\Z_+\}$ is universal.
\end{proof}

\begin{theorem}\label{mi0} Let $X$ be a
topological vector space, $k\in\N$ and $A=(A_1,\dots,A_k)\in
\E(k,X)$ be a EBS$_k$-tuple. Then the strongly continuous group
$\{e^{\langle z,A\rangle}\}_{z\in\K^k}$ is mixing. If additionally
$X$ is Baire separable and metrizable, $\{e^{\langle
z,A\rangle}\}_{z\in\K^k}$ is hereditarily hypercyclic
\end{theorem}

\begin{proof} Assume that $\{e^{\langle z,A\rangle}\}_{z\in\K^k}$ is non-mixing.
Then we can find non-empty open subsets $U$ and $V$ of $X$ and a
sequence $\{z_m\}_{m\in\Z_+}$ in $\K^k$ such that $|z_m|\to\infty$
as $m\to \infty$ and $e^{\langle z_m,A\rangle}(U)\cap V=\varnothing$
for each $m\in\Z_+$. Let $\Sigma$ be the set of pairs $(x,y)\in
X\times X$ for which there exists a sequence $\{x_m\}_{m\in\Z_+}$ of
elements of $X$ such that $x_m\to x$ and $e^{\langle
z_m,A\rangle}x_m\to y$. According to Corollary~\ref{kerim1},
$\KER(A)\times \KER(A)\subseteq \Sigma$. Since $A=(A_1,\dots,A_k)$
is a EBS$_k$-tuple, $\KER(A)$ is dense in $X$ and therefore $\Sigma$
is dense in $X\times X$. In particular, $\Sigma$ meets $U\times V$,
which is not possible since $e^{\langle z_m,A\rangle}(U)\cap
V=\varnothing$ for any $m\in\Z_+$. This contradiction shows that
$\{e^{\langle z,A\rangle}\}_{z\in\K^k}$ is mixing. If $X$ is Baire
separable and metrizable, Proposition~\ref{untr2} implies that
$\{e^{\langle z,A\rangle}\}_{z\in\K^k}$ is hereditarily hypercyclic.
\end{proof}

It is easy to see that if $X$ is a Banach space and $k\in\N$, then
$L(X)^k=\E(k,X)$. Moreover, each operator group of the shape
$\{e^{\langle z,A\rangle}\}_{z\in\K^k}$ is uniformly continuous.
Hence, we get the following corollary of Theorem~\ref{mi0}.

\begin{corollary}\label{mi00} Let $X$ be a separable Banach space and
$(A_1,\dots,A_k)\in L(X)^k$ be a EBS$_k$-tuple. Then $\{e^{\langle
z,A\rangle}\}_{z\in\K^k}$ is a hereditarily hypercyclic uniformly
continuous group.
\end{corollary}

\section{$\ell_1$-sequences, equicontinuous sets and the class $\M$}

\begin{lemma}\label{equi1}Let $Y_0$ and $Y_1$ be closed linear subspaces
of a locally convex space $Y$ such that $Y_0\subset Y_1$ and the
topology of $Y_1/Y_0$ is not weak. Then there is a sequence
$\{f_n:n\in\Z_+\}$ in $Y'$ such that
\begin{itemize}\itemsep=-2pt
\item[\rm (\ref{equi1}.1)] $\phi\subseteq \bigl\{\{f_n(y)\}_{n\in\Z_+}:
y\in Y_1\bigr\};$
\item[\rm (\ref{equi1}.2)] $f_n\bigr|_{Y_0}=0$ for each $n\in\Z_+;$
\item[\rm (\ref{equi1}.3)] $\{f_n:n\in\Z_+\}$ is uniformly
equicontinuous.
\end{itemize}
\end{lemma}

\begin{proof} Since the topology of $Y_1/Y_0$ is not weak, there
exists a continuous seminorm $\widetilde p$ on $Y_1/Y_0$ such that
the closed linear space $\ker \widetilde p=\widetilde p^{-1}(0)$ has
infinite codimension in $Y_1/Y_0$. Clearly the seminorm $p$ on $Y_1$
defined by the formula $p(y)=\widetilde p(y+Y_0)$ is also continuous
and $\ker p$ has infinite codimension and contains $Y_0$. In
particular the space $Y_p=Y_1/\ker p$ endowed with the norm
$\|x+\ker p\|=p(x)$ is an infinite dimensional normed space. Hence
we can choose sequences $\{y_n\}_{n\in\Z_+}$ in $Y_1$ and
$\{g_n\}_{n\in\Z_+}$ in $Y'_p$ such that $\|g_n\|\leq 1$ for each
$n\in\Z_+$ and $g_n(y_k+\ker p)=\delta_{n,k}$ for $n,k\in\Z_+$,
where $\delta_{n,k}$ is the Kronecker $\delta$ (every infinite
dimensional normed space admits a biorthogonal sequence). Now let
the functionals $h_n:Y_1\to \K$ be defined by the formula
$h_n(y)=g_n(y+\ker p)$. The above properties of the functionals
$g_n$ can be rewritten in terms of $h_n$ in the following way
$$
|h_n(y)|\leq p(y)\ \ \text{and}\ \ h_n(y_k)=\delta_{n,k}\ \
\text{for any $n,k\in\Z_+$ and $y\in Y_1$}.
$$
Since any continuous seminorm on a subspace of a locally convex
space extends to a continuous seminorm on the entire space
\cite{shifer,bonet}, we can find a continuous seminorm $q$ on $Y$
such that $q\bigr|_{Y_1}=p$. Applying the Hahn--Banach theorem, we
can find $f_n\in Y'$ for $n\in\Z_+$ such that
$$
f_n\bigr|_{Y_1}=h_n \ \ \text{and}\ \ |f_n(y)|\leq q(y)\ \ \text{for
any $n\in\Z_+$ and $y\in Y$.}
$$
From the last two displays we have $f_n(y_k)=\delta_{n,k}$, which
implies (\ref{equi1}.1) since $y_k\in Y_1$. From the inequality in
the above display it follows that each $|f_n|$ is bounded by $1$ on
the unit ball $W$ of the seminorm $q$. Since $W$ is a neighborhood
of zero in $Y$, condition (\ref{equi1}.3) is satisfied. Since
$Y_0\subseteq \ker p\subseteq \ker q$, from the inequality
$|f_n(y)|\leq q(y)$ it follows that each $f_n$ vanishes on $Y_0$.
That is, (\ref{equi1}.2) is satisfied.
\end{proof}

Applying Lemma~\ref{equi1} with $Y_0=0$ and $Y_1=Y$, we obtain the
following corollary.

\begin{corollary}\label{equi2}Let $Y$ be a locally convex space,
whose topology is not weak. Then there exists a linearly independent
sequence $\{f_n:n\in\Z_+\}$ in $Y'$ such that $\{f_n:n\in\Z_+\}$ is
uniformly equicontinuous and $\phi\subseteq
\bigl\{\{f_n(y)\}_{n\in\Z_+}: y\in Y\bigr\}$.
\end{corollary}

Recall that a subset $D$ of a locally convex space $X$ is called a
{\it disk} if $D$ is bounded, convex and balanced (=is stable under
multiplication by any $\lambda\in\K$ with $|\lambda|\leq1$). The
symbol $X_D$ stands for the space $\spann(D)$ endowed with the norm
being the Minkowskii functional of the set $D$. Boundedness of $D$
implies that the topology of $X_D$ is stronger than the one
inherited from $X$. A disk $D$ in $X$ is called a {\it Banach disk}
if the normed space $X_D$ is complete. It is well-known that a
sequentially complete disk is a Banach disk, see, for instance,
\cite{bonet}. In particular, a compact or a sequentially compact
disk is a Banach disk. We say that $D$ is a {\it Banach $s$-disk} in
$X$ if $D$ is a Banach disk and the Banach space $X_D$ is separable.

\begin{lemma}\label{l11} Let $\{x_n\}_{n\in\Z_+}$ be an
$\ell_1$-sequence in a locally convex space $X$. Then the set
$$
K=\biggl\{\sum_{n=0}^\infty a_nx_n:a\in\ell_1,\ \|a\|_1\leq
1\biggr\}
$$
is a compact and metrizable disk. Moreover, $K$ is a Banach $s$-disk
and $E=\spann\{x_n:n\in \Z_+\}$ is dense in the Banach space $X_K$.
\end{lemma}

\begin{proof} Let $Q=\{a\in \ell_1:\|a\|_1\leq 1\}$ be endowed with
the coordinatewise convergence topology. It is easy to see that $Q$
is a metrizable compact topological space as a closed subspace of
$\D^{\Z_+}$, where $\D=\{z\in\K:|z|\leq 1\}$. Obviously, the map
$\Phi:Q\to K$,  $\Phi(a)=\sum\limits_{n=0}^\infty a_nx_n$ is onto.
It is also easy to see that $\Phi$ is continuous. Indeed, let $p$ be
a continuous seminorm on $X$, $a\in Q$ and $\epsilon>0$. Since
$x_n\to 0$, there is $m\in\Z_+$ such that $p(x_n)\leq \epsilon$ for
$n>m$. Let $\delta=\epsilon(1+p(x_0)+{\dots}+p(x_m))^{-1}$ and
$W=\{b\in Q:|a_j-b_j|<\delta\ \ \text{for}\ \ 0\leq j\leq m\}$. Then
$W$ is a neighborhood of $a$ in $Q$ and for each $b\in W$ we have
$$
p(\Phi(b)-\Phi(a))=p\biggl(\sum_{n=0}^\infty (b_n-a_n)x_n\biggr)\leq
\sum_{n=0}^\infty |b_n-a_n|p(x_n).
$$
Taking into account that $p(x_n)<\epsilon$ for $n>m$ and
$|a_n-b_n|<\delta$ for $n\leq m$, we obtain
$$
p(\Phi(b)-\Phi(a))\leq \delta\sum_{n=0}^m
p(x_m)+\epsilon\sum_{n=m+1}^\infty |b_n-a_n|.
$$
Using the definition of $\delta$ and inequalities $\|a\|_1\leq 1$,
$\|b\|_1\leq 1$, we see that $p(\Phi(b)-\Phi(a))\leq 3\epsilon$.
Since $a$, $p$ and $\epsilon$ are arbitrary, the map $\Phi$ is
continuous. Hence $K$ is compact and metrizable as a continuous
image of a compact metrizable space. Obviously $K$ is convex and
balanced. Hence $K$ is a Banach disk (any compact disk is a Banach
disk). Let us show that $E$ is dense in $X_K$. Take $u\in X_K$. Then
there is $a\in \ell_1$ such that $u=\sum\limits_{k=0}^\infty
a_kx_k$. Clearly, for any $n\in\Z_+$, $u_n=\sum\limits_{k=0}^n
a_kx_k\in E$. Let $\|\cdot\|$ be the norm of the Banach space $X_K$.
Then for any $n\in\Z_+$,
$$
\|u-u_n\|=\biggl\|\sum_{k=n+1}^\infty
a_kx_k\biggr\|\leq\sum_{k=n+1}^\infty |a_k|\to 0\ \ \text{as
$n\to\infty$.}
$$
Hence $E$ is dense in $X_K$ and therefore $X_K$ is separable and $K$
is a Banach $s$-disk.
\end{proof}

\begin{lemma}\label{bp1} Let $X$ be a separable metrizable
topological vector space and $\{f_n\}_{n\in\Z_+}$ be a linearly
independent sequence in $X'$. Then there exist sequences
$\{x_n\}_{n\in\Z_+}$ in $X$ and $\{\alpha_{k,j}\}_{k,j\in\Z_+,\
j<k}$ in $\K$ such that $\spann\{x_k:k\in\Z_+\}$ is dense in $X$,
$g_n(x_k)=0$ for $n\neq k$ and $g_n(x_n)\neq 0$ for $n\in\Z_+$,
where $g_n=f_n+\sum\limits_{j<n}\alpha_{n,j}f_j$.
\end{lemma}

\begin{proof} Let $\{U_n\}_{n\in\Z_+}$ be a base of topology of $X$.
First, we construct inductively sequences
$\{\alpha_{k,j}\}_{k,j\in\Z_+,\ j<k}$ in $\K$ and
$\{y_n\}_{n\in\Z_+}$ in $X$ such that for any $k\in\Z_+$,
\begin{itemize}\itemsep=-2pt
\item[(b1)]$y_k\in U_k$;
\item[(b2)]$g_k(y_k)\neq 0$, where $g_k=f_k+\sum\limits_{j<k}\alpha_{k,j}f_j$;
\item[(b3)]$g_k(y_m)=0$ if $m<k$.
\end{itemize}
Let $g_0=f_0$. Since $f_0\neq 0$, there is $y_0\in U_0$ such that
$f_0(y_0)=g_0(y_0)\neq 0$. This provides us with the base of
induction. Assume now that $n\in\N$ and $y_k$, $\alpha_{k,j}$ with
$j<k<n$ satisfying (b1--b3) are already constructed. According to
(b2) and (b3), we can find $\alpha_{n,0},\dots,\alpha_{n,n-1}\in \K$
such that $g_n(y_m)=0$ for $m<n$, where
$g_n=f_n+\sum\limits_{j<n}\alpha_{n,j}f_j$. Next since $f_j$ are
linearly independent, $g_n\neq 0$ and therefore there is $y_n\in
U_n$ such that $g_n(y_n)\neq 0$. This concludes the description of
the inductive procedure of constructing sequences
$\{\alpha_{k,j}\}_{k,j\in\Z_+,\ j<k}$ and $\{y_n\}_{n\in\Z_+}$
satisfying (b1--b3) for each $k\in\Z_+$.

Using (b2) and (b3), one can easily demonstrate that there is a
sequence $\{\beta_{k,j}\}_{k,j\in\Z_+,\ j<k}$ in $\K$ such that
$g_n(x_k)=0$ for $k\neq m$, where
$x_k=y_k+\sum\limits_{j<k}\beta_{k,j}y_j$. From (b2) and (b3) it
also follows that $g_n(x_n)\neq 0$ for each $n\in\Z_+$. It remains
to notice that according to (b1), $\{y_n:n\in\Z_+\}$ is dense in
$X$. Since $\{y_n:n\in\Z_+\}\subseteq
\spann\{y_n:n\in\Z_+\}=\spann\{x_n:n\in\Z_+\}$, we see that
$\spann\{x_n:n\in\Z_+\}$ is dense in $X$.
\end{proof}

\begin{lemma}\label{ma1} Let $X\in \M$. Then there exist sequences
$\{x_n\}_{n\in\Z_+}$ and $\{f_k\}_{k\in\Z_+}$ in $X$ and $X'$
respectively, such that
\begin{itemize}\itemsep=-2pt
\item[{\rm (\ref{ma1}.1)}]$\{x_n\}_{n\in\Z_+}$ is an $\ell_1$-sequence in $X;$
\item[{\rm (\ref{ma1}.2)}]the space $E=\spann\{x_n:n\in\Z_+\}$ is dense in $X;$
\item[{\rm (\ref{ma1}.3)}]$f_k(x_n)=0$ if $k\neq n$ and $f_k(x_k)\neq0$ for each $k\in\Z_+;$
\item[{\rm (\ref{ma1}.4)}]$\{f_k:k\in\Z_+\}$ is uniformly equicontinuous.
\end{itemize}
Moreover, $\{f_k\}$ can be chosen from the linear span of any
linearly independent uniformly equicontinuous sequence in $X'$.
\end{lemma}

\begin{proof}According to Lemma~\ref{l11}, there exists a Banach
$s$-disk $K$ in $X$ such that $X_K$ is dense in $X$. By
Corollary~\ref{equi2}, there is a linearly independent sequence
$\{g_n\}_{n\in\N}$ in $X'$ such that $\{g_n:n\in\N\}$ is uniformly
equicontinuous. Since $X_K$ is dense in $X$, the functionals
$g_n\bigr|_{X_K}$ on $X_K$ are linearly independent. Applying
Lemma~\ref{bp1} to the sequence $\{g_n\bigr|_{X_K}\}$, we find that
there exist sequences $\{y_n\}_{n\in\Z_+}$ in $X_K$ and
$\{\alpha_{k,j}\}_{k,j\in\Z_+,\ j<k}$ in $\K$ such that
$E=\spann\{y_k:k\in\Z_+\}$ is dense in $X_K$, $h_n(y_k)=0$ for
$n\neq k$ and $h_n(y_n)\neq 0$ for $n\in\Z_+$, where
$h_n=g_n+\sum\limits_{j<n}\alpha_{n,j}g_j$. Let $q$ be the norm of
the Banach space $X_K$. Consider $f_n=c_nh_n$, where
$c_n=\left(1+\sum\limits_{j<n}|\alpha_{n,j}|\right)^{-1}$. Since
$\{g_n:n\in\N\}$ is uniformly equicontinuous and
$h_n=g_n+\sum\limits_{j<n}\alpha_{n,j}g_j$, we immediately see that
$\{f_n:n\in\N\}$ is uniformly equicontinuous. Next, let
$x_n=b_ny_n$, where $b_n=2^{-n}q(x_n)^{-1}$. Since $x_n$ converges
to 0 in the Banach space $X_K$, $\{x_n\}_{n\in\N}$ is an
$\ell_1$-sequence in $X_K$. Since the topology of $X_K$ is stronger
than the one inherited from $X$, $\{x_n\}_{n\in\N}$ is an
$\ell_1$-sequence in $X$. Since
$\spann\{x_n:n\in\Z_+\}=\spann\{y_n:n\in\Z_+\}$ is dense in $X_K$,
it is also dense in $X_K$ in the topology inherited from $X$. Since
$X_K$ is dense in $X$, $\spann\{x_n:n\in\Z_+\}$ is dense in $X$.
Finally since $f_n(x_k)=c_nb_kh_n(y_k)$, we see that $f_n(x_k)=0$ if
$n\neq k$ and $f_n(x_n)\neq 0$ for any $n\in\Z_+$. Thus conditions
(\ref{ma1}.1--\ref{ma1}.4) are satisfied.
\end{proof}

\begin{lemma}\label{dera} Let $X$ be a locally convex space, whose topology is
not weak and $Y$ be a locally convex space admitting an
$\ell_1$-sequence with dense span. Then there is $T\in L(Y,X)$ such
that $T(X)$ is dense in $Y$.
\end{lemma}

\begin{proof} Let $\{y_n\}_{n\in\Z_+}$ be an $\ell_1$-sequence in
$Y$ with dense span. By Corollary~\ref{equi2}, there is a uniformly
equicontinuous sequence $\{f_n\}_{n\in\Z_+}$ in $X'$ such that
$\phi\subseteq \bigl\{\{f_n(x)\}_{n\in\Z_+}: x\in X\bigr\}$.
Consider the linear operator $T:X\to Y$ defined by the formula
$Tx=\sum\limits_{n=0}^\infty 2^{-n}f_n(x)y_n$. Since $\{f_n\}$ is
uniformly equicontinuous, there is a continuous seminorm $p$ on $X$
such that $|f_n(x)|\leq p(x)$ for any $x\in X$ and $n\in\Z_+$. Since
$\{y_n\}$ is an $\ell_1$-sequence, Lemma~\ref{l11} implies that the
closed convex balanced hull $Q$ of $\{y_n:n\in\Z_+\}$ is compact and
metrizable. Let $q$ be the Minkowskii functional of $Q$ (=the norm
on $Y_Q$). It is easy to see that $q(Tx)\leq p(x)$ for each $x\in
X$. Hence $T$ is continuous as an operator from $X$ to the Banach
space $Y_Q$. Since the latter carries the topology stronger than the
one inherited from $Y$, $T\in L(X,Y)$. Next, the inclusion
$\phi\subseteq \bigl\{\{f_n(x)\}_{n\in\Z_+}: x\in x\bigr\}$ implies
that $T(X)$ contains the linear span of $\{y_n:n\in\Z_+\}$, which is
dense in $Y$. Thus $T$ has dense range.
\end{proof}

\begin{lemma}\label{ma2} Let $X$ be an infinite dimensional locally convex space
and assume that there exist $\ell_1$-sequences with dense span in
both $X$ and $X'_\beta$. Then there exist sequences
$\{x_n\}_{n\in\Z_+}$ and $\{f_k\}_{k\in\Z_+}$ in $X$ and $X'$
respectively, such that
\begin{itemize}\itemsep=-2pt
\item[{\rm (\ref{ma2}.1)}]$\{x_n\}_{n\in\Z_+}$ is an $\ell_1$-sequence in $X$ and
$\{f_n\}_{n\in\Z_+}$ is an $\ell_1$-sequence in $X'_\beta;$
\item[{\rm (\ref{ma2}.2)}]$E=\spann\{x_n:n\in\Z_+\}$ is dense in $X$ and
$F=\spann\{f_n:n\in\Z_+\}$ is dense in $X'_\beta;$
\item[{\rm (\ref{ma2}.3)}]$f_k(x_n)=0$ if $k\neq n$ and $f_k(x_k)\neq 0$ for each
$k\in\Z_+$.
\end{itemize}
If the original $\ell_1$-sequence in $X'_\beta$ is uniformly
equicontinuous, than we can also ensure that
\begin{itemize}\itemsep=-2pt
\item[{\rm (\ref{ma2}.4)}]$\{f_n:n\in\Z_+\}$ is uniformly equicontinuous.
\end{itemize}
\end{lemma}

\begin{proof} By Lemma~\ref{l11}, there exists a Banach
$s$-disk $K$ in $X$ such that $X_K$ is dense in $X$. Fix an
$\ell_1$-sequence $\{g_n\}_{n\in\N}$ in $X'_\beta$ with dense span.
Since any sequence of elements of a linear space with infinite
dimensional span has a linearly independent subsequence with the
same span, we, passing to a subsequence, if necessary, can assume
that $g_n$ are linearly independent. Since $X_K$ is dense in $X$,
the functionals $g_n\bigr|_{X_K}$ on $X_K$ are linearly independent.
Let $D$ be the closed absolutely convex hull of the set
$\{g_n:n\in\Z_+\}$ in $X'_\beta$. By Lemma~\ref{l11}, $D$ is a
Banach disk in $X'_\beta$. Applying Lemma~\ref{bp1} to the sequence
$\{g_n\bigr|_{X_K}\}$, we find that there exist sequences
$\{y_n\}_{n\in\Z_+}$ in $X_K$ and $\{\alpha_{k,j}\}_{k,j\in\Z_+,\
j<k}$ in $\K$ such that $E=\spann\{y_k:k\in\Z_+\}$ is dense in
$X_K$, $h_n(y_k)=0$ for $n\neq k$ and $h_n(y_n)\neq 0$ for
$n\in\Z_+$, where $h_n=g_n+\sum\limits_{j<n}\alpha_{n,j}g_j$. Let
$q$ be the norm of the Banach space $X_K$ and $p$ be the norm of the
Banach space $X'_D$. Consider $f_n=c_nh_n$, where
$c_n=2^{-n}\left(1+\sum\limits_{j<n}|\alpha_{n,j}|\right)^{-1}$.
Since $p(g_n)\leq 1$ for each $n\in\Z_+$ and
$h_n=g_n+\sum\limits_{j<n}\alpha_{n,j}g_j$, we immediately see that
$p(f_n)\leq 2^{-n}$ for any $n\in\Z_+$. Hence $\{f_n\}_{n\in\Z_+}$
is an $\ell_1$-sequence in the Banach space $X'_D$. Since the
topology of $X'_D$ is stronger than the one inherited from
$X'_\beta$, $\{f_n\}_{n\in\Z_+}$ is an $\ell_1$-sequence in
$X'_\beta$. Similarly, $\{f_n:n\in\Z_+\}$ is uniformly
equicontinuous if $\{g_n:n\in\Z_+\}$ is. Since the sequences
$\{g_n\}_{n\in\Z_+}$ and $\{f_n\}_{n\in\Z_+}$ have the same spans,
the span of $\{f_n:n\in\Z_+\}$ is dense in $X'_\beta$. Next, let
$x_n=b_ny_n$, where $b_n=2^{-n}q(x_n)^{-1}$. Since $x_n$ converges
to 0 in the Banach space $X_K$, $\{x_n\}_{n\in\N}$ is an
$\ell_1$-sequence in $X_K$. Since the topology of $X_K$ is stronger
than the one inherited from $X$, $\{x_n\}_{n\in\N}$ is an
$\ell_1$-sequence in $X$. Since
$\spann\{x_n:n\in\Z_+\}=\spann\{y_n:n\in\Z_+\}$ is dense in $X_K$,
it is also dense in $X_K$ in the topology inherited from $X$. Since
$X_K$ is dense in $X$, $\spann\{x_n:n\in\Z_+\}$ is dense in $X$.
Finally since $f_n(x_k)=c_nb_kh_n(y_k)$, we see that $f_n(x_k)=0$ if
$n\neq k$ and $f_n(x_n)\neq 0$ for any $n\in\Z_+$. Thus conditions
(\ref{ma2}.1--\ref{ma1}.3) (and (\ref{ma1}.4) if $\{g_n:n\in\Z_+\}$
is uniformly equicontinuous) are satisfied.
\end{proof}

\begin{lemma}\label{clM} The class $\M$ is stable under finite or
countable products. Moreover, if $X\in\M$, then $X\times\K\in\M$.
\end{lemma}

\begin{proof} It is clear that if the topology of one of the
locally convex spaces $X_j$ is not weak, then so is the topology of
their product $\prod X_j$. Thus the only thing we have to worry
about is the existence of an $\ell_1$-sequence with dense span. Let
$J$ be a finite or countable infinite set, $X_j\in \M$ for each
$j\in J$ and $X=\prod\limits_{j\in J}X_j$. Let for any $j\in J$,
$\{x_{j,n}\}_{n\in\Z_+}$ be an $\ell_1$-sequence in $X_j$ with dense
span. For $(j,n)\in J\times \Z_+$ let $u_{j,n}\in X$ be such that
$j^{\rm th}$ component of $u_{j,n}$ is $x_{j,n}$, while all other
components are $0$. Clearly $\{u_{j,n}:j\in J,\ n\in\Z_+\}$ is a
countable subset of $X$. It is easy to see that enumerating this
subset by elements of $\Z_+$, we get an $\ell_1$-sequence in $X$
with dense span. The proof of the second part of the lemma is even
easier: one have just to add one vector $(0,1)$ to an
$\ell_1$-sequence with dense span in $X=X\times\{0\}$, to obtain an
$\ell_1$-sequence with dense span in $X\times\K$.
\end{proof}

\section{Proof of Theorem~\ref{saan}}

Let $X\in\M$. By Lemma~\ref{ma1}, there exist sequences
$\{x_n\}_{n\in\Z_+}$ and $\{f_k\}_{k\in\Z_+}$ in $X$ and $X'$
respectively satisfying (\ref{ma1}.1--\ref{ma1}.4). Uniform
equicontinuity of $\{f_n\}$ is equivalent to the existence of a
continuous seminorm $p$ on $X$ such that each $|f_n|$ is bounded by
$1$ on the unit ball $\{x\in X:p(x)\leq 1\}$ of $p$. Since $\{x_n\}$
is an $\ell_1$-sequence in $X$, Lemma~\ref{l11}, absolutely convex
closed hull $K$ of $\{x_n:n\in\Z_+\}$ is compact and is a Banach
disk in $X$. Let $q$ be the norm of the Banach space $X_K$. Then
$q(x_n)\leq 1$ for each $n\in\Z_+$. Let $c=\sup\{p(x):x\in K\}$.
Compactness of $K$ implies that $c$ is finite. Clearly $p(x)\leq
cq(x)$ for any $x\in Y$.

\begin{lemma}\label{ex} Let $\alpha,\beta:\Z_+\to\Z_+$ be any maps
and $a=\{a_n\}_{n\in\Z_+}\in \ell_1$. Then the formula
\begin{equation}\label{T}
Tx=\sum_{n\in\Z_+} a_n f_{\alpha(n)}(x)x_{\beta(n)}
\end{equation}
defines a linear operator on $X$. Moreover, the series
$\sum\limits_{n=0}^\infty\frac{z^n}{n!}T^nx=e^{zT}x$ converges in
$X$ for any $x\in X$ and $z\in\K$ and
\begin{equation}\label{esti}
q(Tx)\leq \|a\|p(x),\quad q(e^{zT}x-x)\leq c^{-1}e^{c\|a\||z|}p(x)\
\ \text{for each $x\in X$ and $z\in\K$},
\end{equation}
where $\|a\|$ is the $\ell_1$-norm of $a$.
\end{lemma}

\begin{proof} Condition (\ref{ma1}.4) implies that the sequence
$\{f_{\alpha(n)}(x)\}$ is bounded for any $x\in X$. Since $\{x_n\}$
is an $\ell_1$-sequence and $a\in\ell_1$, we see that the series in
(\ref{T}) converges for any $x\in X$ and therefore defines a linear
operator on $X$. Moreover, if $p(x)\leq 1$, then $|f_k(x)|\leq 1$
for each $k\in\Z_+$ and since $q(x_m)\leq 1$ for $m\in\Z_+$, formula
(\ref{T}) implies that $q(Tx)\leq \|a\|$. Hence $q(Tx)\leq
\|a\|p(x)$ for each $x\in X$. Then  $q(T^2x)\leq \|a\|p(Tx)\leq
c\|a\|q(Tx)\leq c\|a\|^2p(x)$. Iterating this argument, we see that
$$
q(T^nx)\leq c^{n-1}\|a\|^n p(x)\ \ \text{for each $n\in\N$.}
$$
It follows that for any $x\in X$, the series
$\sum\limits_{n=1}^\infty \frac{z^n}{k!}T^kx$ is absolutely
convergent in the Banach space $X_K$ (and therefore in $X$) for any
$z\in\K$. Naturally we denote its sum as $e^{zT}x$. Moreover, using
the above display, we obtain
$$
q(e^{zT}x-x)=q\biggl(\sum_{n=1}^\infty \frac{z^n}{n!}T^nx\biggr)\leq
c^{-1}p(x)\sum_{n=1}^\infty
{(c\|a\||z|)^n}{n!}<c^{-1}e^{c\|a\||z|}p(x).
$$
Thus (\ref{esti}) is satisfied.
\end{proof}

Now fix a bijection $\phi:\Z_+^k\to\Z_+$. By symbol $e_j$ we denote
the element of $\Z_+^k$ defined by $(e_j)_l=\delta_{j,l}$, where
$\delta_{j,l}$ is the Kronecker delta. For $n\in\Z_+^k$, we also
write $|n|=n_1+{\dots}+n_k$. For each $m\in\Z_+$ let
$$
\epsilon_m=\min_{n\in\Z_+^k\atop |n|=m+1}
\bigl|f_{\phi(n)}(x_{\phi(n)})\bigr|.
$$
According to (\ref{ma1}.1), (\ref{ma1}.3) and (\ref{ma1}.4),
$\{\epsilon_m\}$ is a bounded sequence positive numbers. Next, we
pick any sequence $\{\alpha_m\}_{m\in\Z_+}$ of positive numbers such
that
\begin{equation}\label{al}
\alpha_{m+1}\geq 2^m\alpha_m\epsilon_m^{-1}\quad\text{for any
$m\in\Z_+$}
\end{equation}
and consider the operators $A_j\in L(X)$ for $1\leq j\leq k$ defined
by the formula
$$
A_jx=\sum_{n\in\Z_+^k} \frac{\alpha_{|n|}f_{\phi(n+e_j)}(x)}
{\alpha_{|n|+1}f_{\phi(n+e_j)}(x_{\phi(n+e_j)})} x_{\phi(n)}.
$$
From the estimates (\ref{al}) it follows that the series defining
$A_j$ can be written as
\begin{equation*}
A_jx=\sum_{n\in\Z_+^k} c_{j,n}f_{\phi(n+e_j)}(x) x_{\phi(n)}\ \
\text{with $0<|c_{j,n}|<2^{-|n|}$}.
\end{equation*}
Clearly $A_j$ have shape (\ref{T}) with $\|a\|\leq
C=\sum\limits_{n\in\Z_+^k}2^{-|n|}$. By Lemma~\ref{ex}, $A_j$ are
linear operators on $X$ satisfying $q(Tx)\leq Cp(x)$ for any $x\in
X$. Hence $A_j$ are continuous as operators from $X$ to the Banach
space $X_K$. Since $X_K$ carries a stronger topology than the one
inherited from $X$, $A_j\in L(X)$ for $1\leq j\leq k$. Next, using
the definition of $A_j$, it is easy to verify that
$A_jA_lx_n=A_lA_jx_n$ for any $1\leq j<l\leq k$ and $n\in\Z_+$.
Indeed, for any $n\in\Z_+$, there is a unique $m\in\Z_+^k$ such that
$n=\phi(m)$. If either $m_j=0$ or $m_l=0$, from the definition of
$A_j$ and $A_l$ it follows that $A_jA_lx_n=A_lA_jx_n=0$. If $m_j\geq
1$ and $m_l\geq 1$ from the same definition we obtain
$A_jA_lx_n=A_lA_jx_n=\frac{\alpha_{|m|-2}}{\alpha_{|m|}}x_{\phi(m-e_j-e_l)}$.
From (\ref{ma1}.2) it follows now that the operators $A_1,\dots,A_n$
are pairwise commuting. Next, for $z\in \K^k$, the operator $\langle
z,A\rangle$ has shape (\ref{T}) with $\|a\|\leq C\|z\|$, where
$\|z\|$ is the $\ell_1$-norm of $z$. By Lemma~\ref{ex}, the series
$\sum\limits_{n=0}^\infty \frac1{n!}\langle z,A\rangle^n$ converges
pointwisely to a linear operator $e^{\langle z,A\rangle}$ and
\begin{equation}\label{esti1}
q(e^{\langle z,A\rangle}x-x)\leq c^{-1}e^{cC\|z\|}p(x)\ \ \text{for
any $x\in X$ and $z\in \K^k$}.
\end{equation}
Exactly as for the operators $A_j$, we see that $e^{\langle
z,A\rangle}\in L(X)$ for any $z\in \K^k$. From (\ref{esti1}), if
$z_n\to z$ in $\K^k$, then $e^{\langle z_n,A\rangle}x\to e^{\langle
z,A\rangle}x$ uniformly on $\{x\in X:p(x)<1\}$. By
Proposition~\ref{expo}, $\{e^{\langle z,A\rangle}\}_{z\in\K^k}$ is a
uniformly continuous group and the map $z\mapsto e^{\langle
z,A\rangle}x$ is holomorphic for any $x\in X$ provided $\K=\C$. The
proof will be complete if we show that the group $\{e^{\langle
z,A\rangle}\}_{z\in\K^k}$ is hereditarily hypercyclic. To this end
consider the restrictions $B_j$ of $A_j$ to $X_K$ as bounded linear
operators on the Banach space $X_K$. Then $B_j$ commute with each
other as the restrictions of commuting operators. Let us show that
the operator group $\{e^{\langle z,B\rangle}\}_{z\in\K^k}$ on $X_K$
is hereditarily hypercyclic. In view of Corollary~\ref{mi00}, it
suffices to verify that $B=(B_1,\dots,B_k)$ is a EBS$_k$-tuple. We
already know that $B_j$ are pairwise commuting. From the definition
of $B_j$ it is easy to see that $\ker B_j^m$ contains
$E_{j,m}=\spann\{x_{\phi(n)}:n\in\Z_+^k,\ n_j\leq m-1\}$. Using this
fact it is easy to see that for any $m\in\N^k$, the set
$\kappa(m,B)$ defined by (\ref{KERk}) contains
$$
E_m=\bigcap_{j=1}^k E_{j,m_j}=\spann\{x_{\phi(n)}:n\in\Z_+^k,\
n_j\leq m_j-1,\ 1\leq j\leq k\}.
$$
Hence the space $\KER B$ defined in (\ref{KERk}) contains
$\bigcup\limits_{m\in\N^k}E_m=\spann\{x_n:n\in\Z_+\}$. By
Lemma~\ref{l11}, $\KER(B)$ is dense in $X_K$. Hence $B$ is an
EBS$_k$-tuple. Thus according to Corollary~\ref{mi00}, $\{e^{\langle
z,B\rangle}\}_{z\in\K^k}$ is hereditarily hypercyclic. Since $X_K$
is dense in $X$ and carries a topology stronger than the one
inherited from $X$, $\{e^{\langle z,A\rangle}\}_{z\in\K^k}$ is
hereditarily hypercyclic. By Proposition~\ref{untr2}, $\{e^{\langle
z,A\rangle}\}_{z\in\K^k}$ is mixing.

\section{Theorem~\ref{main}: proof and applications}

\begin{proof}[Proof of Theorem~$\ref{main}$] Assume that $T$ is non-mixing.
Then we can choose non-empty open subsets $U$ and $V$ of $X$ and a
strictly increasing sequence $\{n_k\}_{k\in\Z_+}$ of positive
integers such that $T^{n_k}(U)\cap V=\varnothing$ for any
$k\in\Z_+$. Let $\Sigma$ be the set of $(x,y)\in X\times X$ for
which there exists a sequence $\{x_k\}_{n\in\Z_+}$ of elements of
$X$ such that $x_k\to x$ and $T^{n_k}x_k\to y$. According to
Theorem~\ref{omt}, $\Sigma$ contains $\Lambda(T)\times \Lambda(T)$.
Since $\Lambda(T)$ is dense in $X$, we see that $\Sigma$ is dense in
$X\times X$. Hence $\Sigma$ intersects $U\times V$, which is
impossible since $T^{n_k}(U)\cap V=\varnothing$ for any $k\in\Z_+$.
This contradiction shows that $T$ is mixing. If $X$ is Baire
separable and metrizable, then by Proposition~\ref{untr2}, $T$ is
hereditarily hypercyclic.
\end{proof}

\subsection{An extension of the Salas theorem}

It is easy to see that $\KER T\subseteq \Lambda(I+T)$ for any linear
operator $T$. Thus Theorem~\ref{main} implies the following
corollary.

\begin{corollary}\label{main1} Let $T$ be an extended backward
shift on a topological vector space $X$. Then $I+T$ is mixing. If
additionally $X$ is Baire, separable and metrizable, then $I+T$ is
hereditarily hypercyclic.
\end{corollary}

Recall that a backward weighted shift on $\ell_p=\ell_p(\Z_+)$ for
$1\leq p<\infty$ is the operator $T$ acting on the canonical basis
$\{e_n\}_{n=0}^\infty$ of $\ell_p$ as follows: $Te_0=0$ and
$Te_n=w_ne_{n-1}$ for $n\geq 1$, where $w=\{w_n\}_{n\in\N}$ is a
bounded sequence of non-zero numbers in $\K$. Clearly any backward
weighted shift is a generalized backward shift and therefore is an
extended backward shift. Hence Corollary~\ref{main1} contains the
Salas theorem on hypercyclicity of the operators $I+T$ with $T$
being a backward weighted shift as a particular case. It is also
easy to see that if $X$ is a topological vector space and $T\in
L(X)$ is surjective, then $\KER T=\Ker T$. Thus we obtain the
following corollary.

\begin{corollary}\label{main2} Let $X$ be a topological vector
space and $T\in L(X)$ be such that $T(X)=X$ and $\Ker T$ is dense in
$X$. Then $I+T$ is mixing. If additionally $X$ is Baire, separable
and metrizable, then $I+T$ is hereditarily hypercyclic.
\end{corollary}

\subsection{An extension of the Hilden--Wallen theorem}

\begin{lemma}\label{uni} Let
$X$ be a Baire topological space, $Y$ be a second countable
topological space and $\{T_n\}_{n\in \Z_+}$ be a sequence of
continuous maps from $X$ to $Y$. Let also $\Sigma$ be the set of
$(x,y)\in X\times Y$ for which there exists a sequence
$\{x_n\}_{n\in\Z_+}$ of elements of $X$ such that $x_n\to x$ and
$T_nx_n\to y$ as $n\to\infty$. If $\Sigma$ is dense in $X\times Y$,
then $\{T_n:n\in\Z_+\}$ is hereditarily universal.
\end{lemma}

\begin{proof} The density of $\Sigma$ in
$X\times Y$ implies that for any infinite set $A\subseteq \Z_+$,
condition (U3) from Theorem~U is satisfied for the family
$\{T_n:n\in A\}$ and therefore this family is universal.
\end{proof}

Hilden and Wallen \cite{hw} demonstrated that any backward weighted
shift on $\ell_p$ for $1\leq p<\infty$ is supercyclic. Many
particular cases of the following proposition are known, see for
instance \cite{sal1}. We include it here in its full generality for
the sake of completeness.

\begin{proposition}\label{super} Let $X$ be a Baire separable metrizable
topological vector space $T\in L(X)$. Suppose also that $T$ has
dense range and dense generalized kernel. Then $T$ is hereditarily
supercyclic.
\end{proposition}

\begin{proof} Since $T$ has dense range, we have that $T^k(X)$ is dense in $X$ for
each $k\in\Z_+$. Thus for any $x\in X$ there exists a sequence
$\{u^x_k\}_{k\in\Z_+}$ in $X$ such that $T^ku^x_k\to x$ as
$k\to\infty$. Fix a dense countable set $B\subseteq X$. Since $X$ is
metrizable and $B$ is countable, we can choose a sequence
$\{\lambda_k\}_{k\in\Z_+}$ of positive numbers such that
$\lambda_k^{-1}u^x_k\to 0$ as $k\to\infty$ for each $x\in B$.

Let now $T_k=\lambda_k T^k$ for each $k\in\Z_+$ and $\Sigma$ be the
set of $(x,y)\in X\times X$ for which there exists a sequence
$\{x_k\}_{k\in\Z_+}$ of elements of $X$ such that $x_k\to x$ and
$T_kx_k\to y$ as $k\to\infty$. For any $x\in B$, let
$x_k=\lambda_k^{-1}u^x_k$. Then $x_k\to 0$ as $k\to\infty$ and
$T_kx_k=x$ for each $k\in\Z_+$. Hence $\{0\}\times
B\subseteq\Sigma$. On the other hand, for any $x\in\Ker T$, we have
$T_k x=0$ for all sufficiently large $k$ and therefore $T_k x\to 0$
as $k\to\infty$. Considering the constant sequence $x_k=x$ for
$k\in\Z_+$, we see that $\Ker T\times\{0\}\subseteq \Sigma$.

Finally, observe that $\Sigma$ is a linear subspace of $X\times X$.
Hence $\Ker T\times B\subseteq \Sigma$. Since both $B$ and $\Ker T$
are dense in $X$, we see that $\Sigma$ is dense in $X\times X$. By
Lemma~\ref{uni}, for each infinite set $A\subseteq\Z_+$, the family
$\{T_n:n\in A\}=\{\lambda_nT^n:n\in A\}$ is universal. Hence $T$ is
hereditarily supercyclic.
\end{proof}

\subsection{Remarks on Theorem~\ref{main}}

Theorem~\ref{main} is reminiscent of the following criterion of
hypercyclicity of Bayart and Grivaux \cite{bg} in terms of the
unimodular point spectrum.

\begin{thmbg}\it Let $X$ be a complex separable infinite dimensional Banach
space, $T\in L(X)$ and assume that there exists a continuous Borel
probability measure $\mu$ on the unit circle $\T$ such that for each
Borel set $A\subseteq \T$ with $\mu(A)=1$, the space
$$
\spann\left(\bigcup_{z\in A}\ker(T-zI)\right)
$$
is dense in $X$. Then $T$ is hypercyclic.
\end{thmbg}\rm

It is worth noting that in the case of operators on Banach spaces,
neither Theorem~\ref{main} implies the result of Bayart and Grivaux,
nor their result implies Theorem~\ref{main}, see Examples \ref{ex1},
\ref{ex2} and \ref{tap} below. Theorem~\ref{main} is also strictly
stronger than Proposition~$2.2$ in the article \cite{herw} by
Herrero and Wang, which is a key tool in the proof of the main
result in \cite{herw} that any element of the operator norm closure
of the set of hypercyclic operators on $\ell_2$ is a compact
perturbation of a hypercyclic operator. Namely, they assume that the
span taken in (\ref{OmegaT}) is dense taking only into account the
$z$'s for which $T-zI$ has closed range, and this allows them to use
the Kitai criterion. We would like also to mention the following
fact.

\begin{proposition}\label{griii} Let $X$ be a locally convex space
and $Y$ be a closed linear subspace of $X$ such that $X$ admits an
$\ell_1$-sequence with dense span and the topology of $X/Y$ is not
weak. Then there is a hereditarily hypercyclic operator $T\in L(X)$
such that $Ty=y$ for any $y\in Y$.
\end{proposition}

\begin{proof} Since the topology of $X/Y$ is not weak,
Lemma~\ref{equi1} implies that there is a linearly independent
uniformly equicontinuous sequence $\{g_n\}_{n\in\Z_+}$ such that
$Y\subseteq \ker g_n$ for each $n\in\Z_+$. By Lemma~\ref{ma1}, we
can find sequences $\{x_n\}_{n\in\Z_+}$ and $\{f_k\}_{k\in\Z_+}$ in
$X$ and $X'$ respectively, such that conditions
(\ref{ma1}.1--\ref{ma1}.4) are satisfied and each $f_k$ belongs to
$\spann\{g_m:m\in\Z_+\}$. Hence $Y\subseteq \ker f_k$ for any
$k\in\Z_+$. Uniform equicontinuity of $\{f_k\}$ and the fact that
$\{x_k\}$ is an $\ell_1$-sequence implies that the formula
$$
Tx=\sum_{n=0}^\infty 2^{-n}f_{n+1}(x)x_{n}
$$
defines a continuous linear operator on $X$, which also acts
continuously on the Banach space $X_K$, where $K$ is the Banach disk
being the closed convex balanced hull of $\{x_k:k\in\Z_+\}$. By
Lemma~\ref{l11}, $X_K$ is separable and $E=\spann\{x_k:k\in\Z_+\}$
is dense in $X_K$. It is also straightforward to verify that
$\KER(T_K)$ contains $E$, where $T_K\in L(X_K)$ is the restriction
of $T$ to $X_K$. Hence $T_K$ is an extended backward shift on the
separable Banach space $X_K$. By Corollary~\ref{main1}, $I+T_K$ is
hereditarily hypercyclic. Since $X_K$ is dense in $X$ and carries
stronger topology, $I+T\in L(X)$ is hereditarily hypercyclic. Since
$Y\subseteq\ker f_k$ for each $k\in\Z_+$ from the definition of $T$
it follows that $(I+T)y=y$ for each $y\in Y$. Thus $I+T$ satisfies
all required conditions.
\end{proof}

Grivaux \cite{gri1} proved that if $Y$ is a closed linear subspace
of a separable Banach space $X$ such that $X/Y$ is infinite
dimensional, then there exists a hypercyclic $T\in L(X)$ such that
$Ty=y$ for any $y\in Y$. This result is an immediate corollary of
Proposition~\ref{griii}. If we apply Proposition~\ref{griii} in the
case when $X$ is a separable Fr\'echet space, we obtain that
whenever $Y$ is a closed linear subspace of $X$ such that $X/Y$ is
infinite dimensional and non-isomorphic to $\omega$, there is a
hereditarily hypercyclic $T\in L(X)$ such that $Ty=y$ for any $y\in
Y$. It is also worth noting that any non-normable Fr\'echet space
has a quotient isomorphic to $\omega$, see \cite{bonet}. If $X$ is a
non-normable separable Fr\'echet space non-isomorphic to $\omega$
and $Y$ is a closed linear subspace of $X$ such that $X/Y$ is
isomorphic to $\omega$, then there is no hypercyclic operator $T\in
L(X)$ such that $Ty=y$ for each $y\in Y$. Indeed, assume that such a
$T$ does exist. Then $S:X/Y\to X$, $S(x+Y)=Tx-x$ is a continuous
linear operator with dense range from $X/Y$ to $X$. Since $X/Y$ is
isomorphic to $\omega$ and $\omega$ carries the minimal locally
convex topology, $S$ is onto and therefore $X$ is isomorphic to
$\omega$, which is a contradiction.

\begin{example}\label{ex1}
Let $\H$ be the complex Sobolev space $W^{-1,2}(\T)$ which consists
of the distributions $f$ on $\T$ such that
$\sum\limits_{n=-\infty}^\infty |\widehat
f_n|^2(1+n^2)^{-1}<\infty$, where the $\widehat f_n$ denotes the
$n^{\rm th}$ Fourier coefficient of $f$. Then the operator $T$
acting on $\H$ as $Tf(z)=zf(z)$ satisfies the conditions of
Theorem~{\rm BG} (and therefore is hypercyclic) and does not satisfy
the conditions of Theorem~$\ref{main}$.
\end{example}

\begin{proof}For each $z\in\T$, $\ker(T-zI)=\Ker(T-zI)$ is the
one-dimensional space spanned by the Dirac $\delta$-function
$\delta_z$, which does not belong to the range of $T-zI$. Thus
$\Lambda(T)=\{0\}$ and $T$ does not satisfy conditions of
Theorem~\ref{main}. On the other hand, $\spann\{\delta_z:z\in A\}$
is dense in $\H$ for any set $A\subseteq\T$ which is dense in $\T$.
Since any subset of $\T$ of full Lebesgue measure is dense, the
conditions Theorem~BG are satisfied with $\mu$ being the normalized
Lebesgue measure.
\end{proof}

\begin{example}\label{ex2}
If $\H$ is the complex Sobolev space $W^{-2,2}(\T)$ which consists
of the distributions $f$ on $\T$ such that
$\sum\limits_{n=-\infty}^\infty |\widehat
f_n|^2(1+n^2)^{-2}<\infty$, then the operator $Tf(z)=zf(z)$ acting
on $\H$ satisfies the conditions of both Theorem~{\rm BG} and
Theorem~$\ref{main}$.
\end{example}

\begin{proof}For each $z\in\T$, $\delta_z\in \ker(T-zI)\cap (T-zI)(\H)$
and therefore the dense linear span of the set $\{\delta_z:z\in\T\}$
is contained in $\Lambda(T)$. Hence $T$ satisfies the conditions of
Theorem~\ref{main}. As in the above example, $\spann\{\delta_z:z\in
A\}$ is dense in $\H$ for any set $A\subseteq\T$ which is dense in
$\T$ and therefore $T$ satisfies conditions of Theorem~BG.
\end{proof}

\begin{example}\label{ex4}
If $T$ is a quasinilpotent generalized backward shift on a separable
complex Banach space, then $I+T$ satisfies the conditions of
Theorem~\ref{main} and does not satisfy the conditions of
Theorem~BG.
\end{example}

\begin{proof} Since $\sigma_p(I+T)=\{1\}$, conditions of Theorem~BG
are not satisfied for the operator $I+T$. On the other hand, $\KER
T\subseteq \Lambda(I+T)$ is dense and therefore $I+T$ satisfies
conditions of Theorem~\ref{main}.
\end{proof}

Examples of chaotic operators on a complex Hilbert space $H$ which
are not mixing are constructed in \cite{BaGr}. For such an operator,
the linear span of the union of $\ker(T-zI)$ for $z\in\T$ satisfying
$z^n=1$ for some $n\in\N$ is dense in $H$. This shows that the
assumption of Theorem~\ref{main} cannot be relaxed into a weaker
assumption like density of the linear span of the union of
$\ker(T-zI)^{n}$ for $z\in\T$ and $n\in\N$. Note also that the class
of operators $T$ for which $\Lambda(T)$ is dense is closed under
finite direct sums. Moreover, this class for operators acting on
Banach spaces is closed under infinite $c_0$-sums and $\ell_p$-sums
for $1\leq p<\infty$. In particular, $I+T$ is mixing, when $T$ is a
finite or countable $c_0$-sum or $\ell_p$-sum with $1\leq p<\infty$
of (possibly different) backward weighted shifts.

Now we describe another class of operators to which
Theorem~\ref{main} applies. Let $\alpha\in L_\infty[0,1]$ and
$\phi:[0,1]\to [0,1]$ be a Borel measurable map. Consider the
integral operator $ T_{\alpha,\phi}$ defined by
$$
T_{\alpha,\phi}f(x)=\int_0^{\phi(x)}\alpha(t)f(t)\,dt.
$$
It is straightforward to verify that $T_{\alpha,\phi}$ is a compact
linear operator on $L_p[0,1]$ for $1\leq p\leq \infty$.

\begin{example} \label{tap} Assume that $\psi$ is continuous, strictly
increasing, $\psi(x)<x$ for $0<x\leq 1$ and that $\alpha(x)\neq 0$
almost everywhere on $[0,1]$. Then $T=T_{\alpha,\psi}$ acting on
$L_p[0,1]$ for $1\leq p<\infty$ is hereditarily supercyclic and
$I+T$ is mixing.
\end{example}

\begin{proof} Consider the sequence $a_1=\psi (1)$ and
$a_{n+1}=\psi (a_n)$ for $n\in\N$. Clearly $\{a_n\}$ is strictly
decreasing and tends to zero as $n\to\infty$. One can easily verify
that
$$
\ker T^n=\bigl\{f: f\bigr|_{[0,a_n]}=0\bigr\}.
$$
Using this equality it is straightforward to check that $T(\ker
T^{n+1})$ is dense in $\ker T^n$ for each $n\in\N$ and that $\Ker T$
is dense in the entire space. Thus $\KER T$ is dense in $L_p[0,1]$.
That is, $T$ is an extended backward shift. It remains to apply
Proposition~\ref{super} and Corollary~\ref{main1}.
\end{proof}

\section{Universality of generic families \label{hysugener}}

\begin{theorem}\label{toogeneric} Let $X$ be a separable metrizable
Baire topological space, $\Omega$ be a Baire topological space, $A$
be a set, and for each $a\in A$ let $\Psi_a$ be a map from $\Omega$
into the set $C(X,X)$ of continuous maps from $X$ to $X$ such that
\begin{align}
&\text{for any $a\in A$, the map \ $\Phi_a:\Omega\times X\to X$,
$\Phi_a(\alpha,x)=\Psi_a(\alpha)x$ \ is continuous} \label{gen11}
\\
&\text{and $\{(\alpha,x,\Psi_a(\alpha)x)\,\,:\,\, \alpha\in \Omega,\
x\in X,\ a\in A\}$ is dense in $\Omega\times X\times
X$.}\label{gen12}
\end{align}
Then
\begin{equation*}
U=\bigl\{\alpha\in\Omega\,\,:\,\,\text{$\{\Psi_a(\alpha):a\in A\}$
is universal}\bigr\} \ \ \text{is a dense $G_\delta$-subset of
$\Omega$}.
\end{equation*}
\end{theorem}

\begin{proof} Let $\{U_n\}_{n\in\Z_+}$ be a sequence of non-empty open
subsets of $X$, which form a basis of the topology of $X$. By
Theorem~U, $\{\Psi_a(\alpha):a\in A\}$ is universal if and only if
for each $k,m\in\Z_+$ there exists $a\in A$ for which
$\Phi_a(\alpha)(U_k)\cap U_m\neq\varnothing$. Thus we have
\begin{equation*}
U=\bigcap_{k,m=0}^\infty\,\,\bigcup_{x\in U_k}\,\,\bigcup_{a\in A}
\{\alpha\in \Omega:\Phi_a(\alpha,x)\in U_m\}.
\end{equation*}
According to (\ref{gen11}), the sets $\{\alpha\in
\Omega:\Phi_a(\alpha,x)\in U_m\}$ are open in $\Omega$. Hence, the
above display implies that $U$ is a $G_\delta$-subset of $\Omega$.
It remains to show that $U$ is dense in $\Omega$. Let
\begin{equation*}
U_0=\bigl\{(\alpha,x)\in\Omega\times X\,\,:\,\,\text{$x$ is a
universal element for $\{\Psi_a(\alpha):a\in A\}$}\bigr\}.
\end{equation*}
Clearly $U$ is  the projection of $U_0$ onto $\Omega$. On the other
hand, $U_0$ is the set of universal elements of the family
$\{\Phi_a:a\in A\}$. Since the product of two Baire spaces, one of
which is second countable, is Baire (see \cite{oxtoby}),
$\Omega\times X$ is Baire. Applying Theorem~U, we see that $U_0$ is
dense in $\Omega\times X$. Since the projection onto $\Omega$ of a
dense subset of $\Omega\times X$ is dense, we get that $U$ is dense
in $\Omega$.
\end{proof}

We apply the above general result to two types of universality:
hypercyclicity and supercyclicity.

\begin{theorem}\label{generic}
Let $X$ be a separable metrizable Baire topological vector space,
$\Omega$ be a Baire topological space and $\alpha\mapsto T_\alpha$
be a map from $\Omega$ into $L(X)$ such that
\begin{equation}
\text{for any $n\in\Z_+$, the map \ $\Phi_n:\Omega\times X\to X$,
$\Phi_n(\alpha,x)=T_\alpha^nx$ \ is continuous.} \label{gen1}
\end{equation}
Then the following conditions are equivalent$:$
\begin{align}\label{gen2}
&\text{$\{(\alpha,x,T_\alpha^nx)\,\,:\,\, \alpha\in \Omega,\ x\in
X,\ n\in\Z_+\}$ is dense in $\Omega\times X\times X$,}
\\
\label{gen2a} &H=\{\alpha\in\Omega\,\,:\,\,\text{$T_\alpha$ is
hypercyclic}\} \ \ \text{is a dense $G_\delta$-subset of $\Omega$}.
\end{align}
\end{theorem}

\begin{proof}Let $A=\Z_+$ and $\Psi_n:\Omega\to L(X)$ for $n\in A$ be defined as
$\Psi_n(\alpha)=T_\alpha^n$. Applying Theorem~\ref{toogeneric}, we
see that (\ref{gen2}) implies (\ref{gen2a}). Since the set of
hypercyclic vectors of any hypercyclic operator is dense,
(\ref{gen2a}) implies (\ref{gen2}).
\end{proof}

\begin{theorem}\label{genericsup}
Let $X$ be a separable metrizable Baire topological vector space,
$\Omega$ be a Baire topological space and $\alpha\mapsto T_\alpha$
be a map from $\Omega$ to $L(X)$ such that $(\ref{gen1})$ is
satisfied. Then the following conditions are equivalent:
\begin{align}\label{gen3}
&\text{$\{(\alpha,x,\lambda T_\alpha^nx)\,\,:\,\, \alpha\in \Omega,\
x\in X,\ n\in\Z_+,\ \lambda\in\C\}$ is dense in $\Omega\times
X\times X,$}
\\
\label{gen3a} &S=\{\alpha\in\Omega\,\,:\,\,\text{$T_\alpha$ is
supercyclic}\} \ \ \text{is a dense $G_\delta$-subset of $\Omega$.}
\end{align}
\end{theorem}

\begin{proof} Let $A=\Z_+\times\C$ and $\Psi_{n,\lambda}:\Omega\to L(X)$ for
$(n,\lambda)\in A$ be defined as $\Psi_{n,\lambda}(\alpha)=\lambda
T_\alpha^n$. Applying Theorem~\ref{toogeneric}, we see that
(\ref{gen3}) implies (\ref{gen3a}). Since the set of supercyclic
vectors of any supercyclic operator is dense, (\ref{gen3a}) implies
(\ref{gen3}).
\end{proof}

\begin{corollary}\label{generic1} Let $X$ be a separable Baire metrizable
topological vector space, $\Omega$ be a Baire topological space and
$\alpha\mapsto T_\alpha$ be a map from $\Omega$ to $L(X)$ satisfying
$(\ref{gen1})$. Suppose also that for any non-empty open subset $W$
of $\Omega$ and any nonempty open subsets $U$, $V$ of $X$, there
exist $\alpha\in W$, $x\in U$ and $y\in V$ such that
$x,y\in\Lambda(T_\alpha)$, where $\Lambda(T_\alpha)$ is defined in
$(\ref{OmegaT})$. Then the set of $\alpha\in\Omega$ for which
$T_\alpha$ is hypercyclic is a dense $G_\delta$-subset of $\Omega$.
\end{corollary}

\begin{proof}Let $(\alpha,x,y)\in \Omega\times X\times X$ be such that
$x,y\in\Lambda(T_\alpha)$. By Theorem~\ref{omt} $(x,y)$ is in the
closure of the set $\{(u,T_\alpha^nu):u\in X,\ n\in\Z_+\}$. Hence
$(\alpha,x,y)$ is in the closure of the set
$\{(\alpha,u,T_\alpha^nu)\,\,:\,\, \alpha\in \Omega,\ u\in X,\
n\in\Z_+\}$. By the assumptions of the corollary it follows that
this last set is dense in $\Omega\times X\times X$. It remains to
apply Theorem~\ref{generic}.
\end{proof}

The following supercyclicity analog of the above corollary turns out
to be much easier.

\begin{proposition}\label{generic3} Let $X$ be a Baire separable metrizable
topological vector space, $\Omega$ be a Baire topological space and
$\alpha\mapsto T_\alpha$ be a map from $\Omega$ to $L(X)$ satisfying
$(\ref{gen1})$. Suppose also that for any non-empty open subset $W$
of $\Omega$ and any nonempty open subsets $U$ and $V$ of $X$, there
exist $\alpha\in W$, $x\in U$, $y\in V$ and $n\in\N$ such that
$x\in\ker T_\alpha^n$ and $y\in T_\alpha^n(X)$. Then the set of
$\alpha\in\Omega$ for which $T_\alpha$ is supercyclic is a dense
$G_\delta$-subset of $\Omega$.
\end{proposition}

\begin{proof}Let $(\alpha,x,y)\in \Omega\times X\times X$ be such
that there exists $n\in\N$ for which $x\in\ker T_\alpha^n$ and $y\in
T_\alpha^n(X)$. Let $v\in X$ be such that $T_\alpha^nv=y$. For each
$m\in\N$ consider $u_m=x+m^{-1}v$. Then $u_m\to x$ as $m\to\infty$.
Moreover, $mT^nu_m=y$ for any $m\in\N$. Thus for all sufficiently
large $m$, $(\alpha,u_m,mT^nu_m)\in W\times U\times V$. Since $W$,
$U$ and $V$ were arbitrary, the set
$\{(\alpha,u,zT^nu):\alpha\in\Omega,\ n\in\Z_+,\ z\in\K\}$  is dense
in $\Omega\times X\times X$. It remains to apply
Theorem~\ref{genericsup}.
\end{proof}

The obvious inclusion $\KER (T-I)\subset \Lambda(T)$ and the fact
that a map $\alpha\to T_\alpha$ satisfies (\ref{gen1}) if and only
if the map $\alpha\to I+T_\alpha$ satisfies (\ref{gen1}) imply the
following corollary of Corollary~\ref{generic1} and
Proposition~\ref{generic3}.

\begin{corollary}\label{generic2} Let $X$ be a Baire separable metrizable
topological vector space, $\Omega$ be a Baire topological space and
$\alpha\mapsto T_\alpha$ be a map from $\Omega$ to $L(X)$ satisfying
$(\ref{gen1})$. Suppose also that for any non-empty open subset $W$
of $\Omega$ and any nonempty open subsets $U$ and $V$ of $X$, there
exist $\alpha\in W$, $x\in U$ and $y\in V$ such that $x,y\in\KER
T_\alpha$. Then the set of $\alpha\in\Omega$ for which $T_\alpha$ is
supercyclic and $I+T_\alpha$ is hypercyclic is a dense
$G_\delta$-subset of $\Omega$.
\end{corollary}

Now we apply the above general results in order to show, in
particular, that for a generic (in the Baire category sense) compact
quasinilpotent operator $T$ on a separable Hilbert space, $T$ is
supercyclic and $I+T$ is hypercyclic. Throughout the remaining part
of the section $X$ and $Y$ are two infinite dimensional Banach
spaces, $b:X\times Y\to \K$ is a continuous bilinear form separating
points of $X$ and of $Y$. That is, for each non-zero $x\in X$, there
is $y\in Y$ satisfying $b(x,y)\neq 0$ and for each non-zero $y\in
Y$, there is $x\in X$ such that $b(x,y)\neq 0$. In particular, $b$
is a dual pairing between $X$ and $Y$. Recall that the injective
norm on the tensor product $X\otimes Y$ is defined by the formula
$$
\epsilon(\xi)=\sup\biggl\{\biggl|\sum_{j=1}^n
\phi(x_j)\psi(y_j)\biggr|: \phi\in X',\ \psi\in Y',\ \|\phi\|\leq
1,\ \|\psi\|\leq 1\biggr\}\; \text{ for}\ \ \xi=\sum_{j=1}^n
x_j\otimes y_j\in X\otimes Y.
$$
The completion of $X\otimes Y$ with respect to this norm is called
the injective tensor product of $X$ and is $Y$ and denoted
$X\widehat{\otimes}_\epsilon Y$. It is again a Banach space. For
each $\xi\in X\otimes Y$ we consider the linear operators $T_\xi\in
L(X)$ and $S_\xi\in L(Y)$ defined by the formulae
\begin{equation}\label{txsx}
T_\xi x=\sum_{j=1}^n b(x,y_j)x_j\ \ \text{and}\ \ S_\xi
y=\sum_{j=1}^n b(x_j,y)y_j,\ \text{where}\ \ \xi=\sum_{j=1}^n
x_j\otimes y_j.
\end{equation}
Clearly $\xi\mapsto T_\xi$ and $\xi\mapsto S_\xi$ are bounded linear
operators from $X\otimes Y$ endowed with the injective norm into the
Banach spaces $L(X)$ and $L(Y)$ respectively. Hence they admit
unique continuous extensions to $X\widehat{\otimes}_\epsilon Y$,
which we again denote by $\xi\mapsto T_\xi$ and $\xi\mapsto S_\xi$.
Note that $S_\xi$ is the dual operator of $T_\xi$ with respect to
the dual pairing $b$:
\begin{equation}
\label{dupa} b(T_\xi x,y)=b(x,S_\xi y)\ \ \text{for any $x\in X$ and
$y\in Y$.}
\end{equation}

Let now ${\mathcal M}$ be a linear subspace of
$X\widehat{\otimes}_\epsilon Y$ which we suppose to be endowed with
its own $\F$-space topology stronger than the one inherited from
$X\widehat{\otimes}_\epsilon Y$. Suppose also that $X\otimes Y$ is a
dense linear subspace of $\mathcal M$ and  let $\NNN$ be the closure
in $\mathcal M$ of the set $\{\xi\in X\otimes Y:T_\xi\ \ \text{is
nilpotent}\}$. Remark that (\ref{dupa}) implies that $T_\xi^n=0$ if
and only if $S_\xi^n=0$. Thus nilpotency of $T_\xi$ is equivalent to
nilpotency of $S_\xi$ and this implies that $\NNN$ coincides with
the closure in $\mathcal M$ of the set $\{\xi\in X\otimes Y:S_\xi\ \
\text{is nilpotent}\}$.

\begin{proposition}\label{ten}If the Banach space $X$ is separable,
then the set
$$
\{\xi\in \NNN:T_\xi\ \ \text{is supercyclic and}\ \ I+T_\xi\ \
\text{is hypercyclic}\}
$$
is a dense $G_\delta$ subset of $\NNN$. If $X$ and $Y$ are both
separable, then the set
$$
\{\xi\in \NNN:\text{$T_\xi$ and $S_\xi$ are supercyclic and}\ \
\text{$I+T_\xi$ and $I+S_\xi$ are hypercyclic}\}
$$
is a dense $G_\delta$ subset of $\NNN$.
\end{proposition}

\begin{proof} Let $W$ be a non-empty
open subset of $\NNN$ and $U_1,U_2$ be non-empty open subsets of
$X$. Pick $x_1\in U_1$ and $x_2\in U_2$. By the definition of
$\NNN$, there exists $\xi\in X\otimes Y$ such that $\xi\in W$ and
$T_\xi$ is nilpotent. Let $n\in\N$ be such that $T_\xi^n=0$. Clearly
the space
$$
L=\{y\in Y:S_\xi y=0,\ b(x_1,y)=b(x_2,y)=0\}
$$
has finite codimension in $Y$, while $\ker T_\xi$ has finite
codimension in $X$. Since $b$ is a dual pairing of $X$ and $Y$, we
can find $f_j\in L$ and $u_j\in \ker T_\xi$ for $1\leq j\leq 2n$
such that $b(u_k,f_j)=\delta_{k,j}$ for $1\leq j,k\leq 2n$. Consider
the element
$$
\eta=x_1\otimes f_1+x_2\otimes f_{n+1}+\sum_{j=2}^n(u_{j-1}\otimes
f_j+u_{n+j-1}\otimes f_{n+j})\in X\otimes Y
$$
and let $\xi_t=\xi+t\eta$ for $t\in \K$. From the above properties
of $f_j$ and $u_j$ it immediately follows that the range of
$T^n_{\xi_t}$ is contained in the range $Q$ of $T_\xi$ and that the
restrictions of $T_{\xi_t}$ and $T_\xi$ to $Q$ coincide. Hence
$T^{2n}_{\xi_t}=0$ for any $t\in\K$. Thus all $T_{\xi_t}$ are
nilpotent and therefore all $\xi_t$ belong to $\NNN$. Since $W$ is
open in $\NNN$ we can pick $s\in\K\setminus \{0\}$ close enough  to
zero to ensure that $\xi_s\in W$. Now from the definition of $\xi_s$
it immediately follows that $T_{\xi_s}^nu_n=s^nx_1$,
$T_{\xi_s}^nu_{2n}=s^nx_2$, $T_{\xi_s}^n x_1=T_\xi^nx_1=0$ and
$T_{\xi_s}^n x_2=T_\xi^n x_2=0$. Thus
$$
x_1,x_2\in \ker T_{\xi_s}^n\cap T_{\xi_s}^n(X)\subset \KER
T_{\xi_s}.
$$
Summarizing the above, we have found $x_1\in U_1$, $x_2\in U_2$ and
$\xi_s\in W$ such that $x_1,x_2\in \KER T_{\xi_s}$. Using the fact
that $\NNN$ carries a topology stronger than the one defined by the
injective norm, we see that the map $\xi\mapsto T_\xi$ from $\NNN$
to $L(X)$ satisfies (\ref{gen1}). By Corollary~\ref{generic2}, the
set of $\xi\in{\NNN}$ for which $T_\xi$ is supercyclic and $I+T_\xi$
is hypercyclic is a dense $G_\delta$ subset of $\NNN$. The proof of
the first part of the proposition is complete.

Applying the first part of the proposition to $(Y,X)$ instead of
$(X,Y)$ with $b(x,y)$ replaced by $b(y,x)$, we obtain that if $Y$ is
separable, then the set of $\xi\in\NNN$ for which $S_\xi$ is
supercyclic and $I+S_\xi$ is hypercyclic is a dense $G_\delta$
subset of $\NNN$. The second part now follows from the first part
and the fact that the intersection of two dense $G_\delta$ subsets
of a complete metric space is again a dense $G_\delta$ set.
\end{proof}

\subsection{Proof of Theorems~\ref{cococo} and~\ref{cocococo}}

Let us consider the case where $Y=X'$, $b(x,y)=y(x)$ and ${\mathcal
M}=X\widehat{\otimes}_\epsilon X'$. In this case the map $\xi\mapsto
T_\xi$ is an isometry and $\{T_\xi:\xi\in \NNN\}$ is exactly the
operator norm closure of the set of finite rank nilpotent operators.
Moreover, taking into account that $T'_\xi=S_\xi$ for any
$\xi\in{\mathcal M}$, we see that Theorem~\ref{cococo} is an
immediate corollary of Proposition~\ref{ten} for this specific
choice of $\cal M$ and $b$.

Recall that the projective norm on the tensor product $X\otimes Y$
of the  Banach spaces $X$ and $Y$ is defined by the formula
$$
\pi(\xi)=\inf\biggl\{\sum \|x_j\|\|y_j\|\biggl\},
$$
where the infimum is taken over all possible representations of
$\xi$ as a finite sum $\sum x_j\otimes y_j$. If we consider the case
$Y=X'$, $b(x,y)=y(x)$ and ${\mathcal M}=X\widehat{\otimes}_\pi X'$
the completion of $X\otimes X'$ with respect to the projective norm,
then $\{T_\xi:\xi\in {\mathcal N}\}$ is exactly the set of nuclear
quasinilpotent operators, and we see that Theorem~\ref{cocococo} is
an immediate corollary of Proposition~\ref{ten} for this choice of
$\cal M$ and $b$.

\subsection{Proof of Theorem~\ref{nonmix}}

Let $X\in\M$. By Lemma~\ref{ma1}, there exist sequences
$\{x_n\}_{n\in\Z_+}$ and $\{f_n\}_{n\in\Z_+}$ in $X$ and $X'$
respectively such that conditions (\ref{ma1}.1--\ref{ma1}.4) are
satisfied. Since $\{f_n\}$ is uniformly equicontinuous, we can pick
a non-zero continuous seminorm $p$ on $X$ such that $|f_n(x)|\leq 1$
for any $n\in\Z_+$ whenever $p(x)\leq 1$. By Lemma~\ref{l11}, the
closed balanced convex hull $K$ of $\{x_n:n\in\Z_+\}$ a Banach
$s$-disk. That is, the Banach space $X_K$ is separable. It is also
clear from Lemma~\ref{l11} that any $x\in X_K$ has shape
$x=\sum\limits_{n=0}^\infty a_nx_n$ for some $a\in\ell_1$ and
$\|x\|\leq \|a\|_1$, where $\|x\|$ is the norm of $X$ in the Banach
space $X_K$. Consider the bilinear form on $X\times \ell_1$ defined
by the formula
$$
b(x,a)=\sum_{n=0}^\infty a_nf_n(x).
$$
It is easy to see that $b$ is well-defined, continuous and
$b(x,a)\leq p(x)\|a\|_1$ for any $x\in X$ and $a\in\ell_1$.
Moreover, $b$ separates points of $X_K$ and $\ell_1$. Indeed, let
$x\in X_K$, $x\neq 0$. Pick $\alpha\in \ell_1$ such that
$x=\sum\limits_{n=0}^\infty \alpha_nx_n$. Since $x\neq 0$, there is
$m\in\Z_+$ such that $\alpha_m\neq 0$. Using (\ref{ma1}.3), we see
that $b(x,e_m)=f_m(x)=\alpha_m f_m(x_m)\neq 0$, where
$\{e_j\}_{j\in\Z_+}$ is the standard basis in $\ell_1$. Similarly,
let $a\in \ell_1$, $a\neq 0$. Then there is $m\in\Z_+$ for which
$a_m\neq 0$ and therefore, according to (\ref{ma1}.3),
$b(x_m,a)=a_mf_m(x_m)\neq 0$. Thus $b$ separates points of $X_K$ and
$\ell_1$. Let ${\mathcal M}=X_K\widehat{\otimes}_\pi \ell_1$ be the
projective tensor product of the Banach spaces $X_K$ and $\ell_1$
and $\NNN$ be the closure in $\mathcal M$ of the set of $\xi\in
X_K\otimes \ell_1$, for which the operator $T_\xi$ defined in
(\ref{txsx}) is nilpotent. By Proposition~\ref{ten}, the set
$\{\xi\in\NNN:I+T_\xi\ \ \text{and}\ \ I+S_\xi\ \ \text{are both
hypercyclic}\}$ is a dense $G_\delta$ subset of $\NNN$. In
particular, we can pick $\xi\in \NNN\subset X_D\widehat{\otimes}_\pi
\ell_1$ such that $I+T_\xi$ and $I+S_\xi$ are hypercyclic. Using the
theorem characterizing the shape of elements of the projective
tensor product \cite{shifer}, we see that there exist bounded
sequences $\{y_n\}_{n\in\Z_+}$ and $\{w_n\}_{n\in\Z_+}$ in $X_K$ and
$\ell_1$ respectively and $\lambda\in\ell_1$ such that
$\xi=\sum\limits_{n=0}^\infty \lambda_ny_n\otimes w_n$. Then the
operators $T_\xi$ and $S_\xi$ act according to the following
formulae on $X_K$ and $\ell_1$ respectively:
$$
T_\xi x=\sum_{n=0}^\infty \lambda_n b(x,w_n)y_n,\ \ \text{and}\ \
S_\xi a=\sum_{n=0}^\infty \lambda_n b(y_n,a)w_n.
$$
Using boundedness of $\{y_n\}$ in $X_K$ and $\{w_n\}$ in $\ell_1$,
summability of $\{|\lambda_n|\}$ and the definition of $b$, we see
that the right-hand side of the first equality in the above display
defines a continuous linear operators $T:X\to X$, taking values in
$X_K$. Since the restriction $I+T_\xi$ of $I+T$ to $X_K$ is
hypercyclic on $X_K$, $X_K$ is dense in $X$ and $X_K$ carries the
topology stronger than the one inherited from $X$, we see that $I+T$
is hypercyclic (any hypercyclic vector for $I+T_\xi$ is also
hypercyclic for $I+T$). Pick a hypercyclic vector $a$ for $I+S_\xi$.
Since $a\neq 0$ and $b$ separates points of $\ell_1$, the functional
$b(\cdot,a)$ is non-zero. Since $a$ is hypercyclic for $I+S_\xi$, we
can pick a strictly increasing sequence $\{n_k\}_{k\in\Z_+}$ of
positive integers such that $\|(I+S_\xi)^{n_k}a\|_1\leq 1$ for any
$k\in\Z_+$. Then
$$
|b((I+T)^{n_k}x,a)|=|b(x,(I+S_\xi)^{n_k}a)|\leq
p(x)\|(+S_\xi)^{n_k}a\|_1\leq p(x)\ \ \text{for any $x\in X$ and
$k\in\Z_+$}.
$$
Let $U=\{x\in X:p(x)<1\}$ and $V=\{x\in X:|b(x,a)|>1\}$. Clearly $U$
and $V$ are non-empty open subsets of $X$ ($V\neq \varnothing$ since
the functional $b(\cdot,a)$ is non-zero). Moreover, from the above
display it follows that $(I+T)^{n_k}(U)\cap V=\varnothing$ for each
$k\in\Z_+$. Hence $I+T$ is non-mixing. Since $I+T$ is hypercyclic,
the proof of Theorem~\ref{nonmix} is complete.

\subsection{Proof of Theorem~\ref{dual}}

Let $X$ be an infinite dimensional locally convex space, such that
both $X$ and $X'_\beta$ admit $\ell_1$-sequences with dense span. By
Lemma~\ref{ma2}, there exist sequences $\{x_n\}_{n\in\Z_+}$ and
$\{f_n\}_{n\in\Z_+}$ in $X$ and $X'$ respectively, satisfying
(\ref{ma2}.1--\ref{ma2}.3). By Lemma~\ref{l11} the closed balanced
convex hulls $K$ and $D$ of $\{x_n:n\in\Z\}$ and $\{f_n:n\in\Z\}$
are Banach $s$-disks in $X$ and $X'_\beta$ respectively. Moreover,
$X_K$ is dense in $X$ and $X'_D$ is dense in $X'_\beta$. Since $D$
is $\beta$-compact, it is also $\sigma(X',X)$-compact. Hence the
seminorm $p(x)=\sup\{|f(x)|:f\in D\}$ on $X$ is continuous with
respect to the Mackey topology $\tau=\tau(X,X')$. Clearly each
$|f_n|$ is bounded by $1$ on $\{x\in X:p(x)\leq1\}$ and therefore
$\{f_n:n\in\Z\}$ is uniformly equicontinuous. By Lemma~\ref{ma2}, we
can assume that the same holds true if the original
$\ell_1$-sequence in $X'_\beta$ is uniformly equicontinuous. Assume
for time being that either $X$ carries the Mackey topology $\tau$ or
the original $\ell_1$-sequence in $X'_\beta$ is uniformly
equicontinuous. Then $\{f_n:n\in\Z\}$ is uniformly equicontinuous.

Consider the bilinear form on $X\times X'$ defined by the formula
$\beta(x,f)=f(x)$. Clearly $\beta$ separates points of $X$ and $X'$
and $\beta$ is separately continuous on $X\times X'_\sigma$. Since
$X_K$ is dense in $X$ and $X'_D$ is dense in $X'_\sigma$, the
bilinear form $b:X_K\times X'_D\to\K$, being the restriction of
$\beta$ to $X_K\times X'_D$, separates points of $X_K$ and $X'_D$.
Moreover separate continuity of $\beta$ implies separate continuity
of $b$ and therefore continuity of $b$ on $X_K\times X'_D$ by means
of the uniform boundedness principle (every separately continuous
bilinear form on product of Banach spaces is continuous). Let
${\mathcal M}=X_K\widehat{\otimes}_\pi X'_D$ be the projective
tensor product of the Banach spaces $X_K$ and $X'_D$ and $\NNN$ be
the closure in $\mathcal M$ of the set of $\xi\in X_K\otimes X'_D$,
for which the operator $T_\xi$ defined in (\ref{txsx}) is nilpotent.
By Proposition~\ref{ten}, the set $\{\xi\in\NNN:I+T_\xi\ \
\text{and}\ \ I+S_\xi\ \ \text{are both hypercyclic}\}$ is a dense
$G_\delta$ subset of $\NNN$. In particular, we can pick $\xi\in
\NNN\subset X_D\widehat{\otimes}_\pi \ell_1$ such that $I+T_\xi$ and
$I+S_\xi$ are hypercyclic. Using once again the theorem
characterizing the shape of elements of the projective tensor
product, we see that there exist bounded sequences
$\{y_n\}_{n\in\Z_+}$ and $\{g_n\}_{n\in\Z_+}$ in $X_K$ and $X'_D$
respectively and $\lambda\in\ell_1$ such that
$\xi=\sum\limits_{n=0}^\infty \lambda_ny_n\otimes g_n$. Then the
operators $T_\xi$ and $S_\xi$ act according to the following
formulae on $X_K$ and $X'_D$ respectively:
$$
T_\xi x=\sum_{n=0}^\infty \lambda_n g_n(x)y_n,\ \ \text{and}\ \
S_\xi f=\sum_{n=0}^\infty \lambda_n f(y_n)g_n,
$$
where we used the specific shape of our bilinear form. Using
boundedness of $\{y_n\}$ in $X_K$ and $\{g_n\}$ in $X'_D$ and
summability of $\{|\lambda_n|\}$, we see that the right-hand sides
of the equalities in the above display define linear operators $T$
and $S$ on $X$ and $X'$, taking values in $X_K$ and $X'_D$
respectively. Since $\{f_n:n\in\Z_+\}$ is uniformly equicontinuous
and $\{g_n\}$ is bounded in $X'_D$, $\{g_n:n\in\Z_+\}$ is also
uniformly equicontinuous. It follows that $T$ is continuous as an
operator from $X$ to $X_K$ and therefore $T\in L(X)$. It is also
easy to verify that $S=T'$ and therefore $S\in L(X'_\beta)$. Since
the restriction $I+T_\xi$ of $I+T$ to $X_K$ is hypercyclic on $X_K$,
$X_K$ is dense in $X$ and $X_K$ carries the topology stronger than
the one inherited from $X$, we see that $I+T$ is hypercyclic (any
hypercyclic vector for $I+T_\xi$ is also hypercyclic for $I+T$).
Similarly $I+S=(I+T)'$ is hypercyclic on $X'_\beta$ (any hypercyclic
vector for $I+S_\xi$ is also hypercyclic for $I+S$). Hence $I+T$ is
a dual hypercyclic operator. In order to complete the proof of
Theorem~\ref{dual}, it remains to consider the case when $X$ carries
the weak topology $\sigma=\sigma(X,X')$. By the already proven part
of Theorem~\ref{dual}, there is a dual hypercyclic operator $R$ on
$X_\tau$. Since $L(X_\tau)=L(X_\sigma)$, $R\in L(X)$. Since $\tau$
is stronger than $\sigma$, $R$ is hypercyclic on $X=X_\sigma$. Hence
$R$ is dual hypercyclic on $X$. The proof of Theorem~\ref{dual} is
complete.

\subsection{Generic bilateral weighted shifts}

For each $w\in\ell_\infty(\Z)$, $T_w$ stands for the bounded linear
operator acting on $\ell_p(\Z)$, $1\leq p<\infty$ or $c_0(\Z)$,
defined on the canonical basis $\{e_n\}_{n\in\Z}$ by
$$
T_we_n=w_ne_{n-1}\qquad \text{for\ \ }n\in\Z.
$$
If additionally $w_n\neq 0$ for each $n\in\Z$, the operator $T_w$ is
called the {\it bilateral weighted shift} with  weight sequence $w$.
Cyclic properties of bilateral weighted shifts have been intensely
studied. Hypercyclicity and supercyclicity of bilateral weighted
shifts were characterized by Salas \cite{sal,sal1} in terms of the
weight sequences. It was observed in \cite[Proposition~5.1]{shk1}
that the Salas conditions admit a simpler equivalent form:

\begin{thmsas}
Let $T$ be a bilateral weighted shift with weight sequence $w$
acting on $\ell_p(\Z)$ with $1\leq p<\infty$ or $c_0(\Z)$. Then $T$
is hypercyclic if and only if for any $m\in\Z_+$,
\begin{equation}
\ilim\limits_{n\to\infty} \max\{\widetilde w(m-n+1,m), (\widetilde
w(m+1,m+n))^{-1}\}=0 \label{sal3}
\end{equation}
and $T$ is supercyclic if and only if for any $m\in\Z_+$,
\begin{equation}
\ilim\limits_{n\to+\infty}
\widetilde{w}(m-n+1,m)\widetilde{w}(m+1,m+n)^{-1}=0, \label{sal4}
\end{equation}
where $\widetilde{w}(a,b)=\prod\limits_{j=a}^b\,|w_j|$\ \ for
$a,b\in\Z$ with $a\leq b$.
\end{thmsas}

We address the issue of hypercyclicity and supercyclicity of generic
bilateral weighted shifts in the Baire category sense. For each
$c>0$ let
\begin{equation*}
B_c=\{w\in \ell_\infty(\Z):\|w\|_\infty\leq c\}.
\end{equation*}
 Clearly
$B_c$ endowed with the coordinatewise convergence topology is a
compact metrizable topological space.

\begin{theorem}\label{genshi}Let $1\leq p<\infty$. For each $c>1$ the set of
$w\in B_c$ for which $T_w$ acting on $\ell_p(\Z)$ is hypercyclic is
a dense $G_\delta$-subset of $B_c$. For each $c>0$ the set of $w\in
B_c$ for which $T_w$ is supercyclic and $I+T_w$ is hypercyclic is a
dense $G_\delta$-subset of $B_c$.
\end{theorem}

It is worth noting that if $w_n=0$ for some $n\in\Z$ then the range
of the operator $T_w$ is not dense and therefore $T_w$ can not be
supercyclic. Thus any $T_w$ for $w$ from the dense $G_\delta$-sets
in the above theorem are indeed bilateral weighted shifts. Recall
that $T_w$ is compact if and only if $w\in c_0(\Z)$. For compact
bilateral weighted shifts we can replace the coordinatewise
convergence topology on the space of weights by stronger topologies.

\begin{theorem}\label{genshi1}Let $1\leq p<\infty$. Let $E$ be a linear
subspace of $c_0(\Z)$ carrying its own $\F$-space topology stronger
than the one inherited from $c_0(\Z)$ and such that the space
$\phi(\Z)$ of sequences with finite support is densely contained in
$E$. Then the set of $w\in E$ for which $T_w$ acting on $\ell_p(\Z)$
is supercyclic and $I+T_w$ is hypercyclic is a dense
$G_\delta$-subset of $E$.
\end{theorem}

\begin{remark}
As a corollary of the above theorem we obtain that the set of
weights $w\in c_0(\Z)$ for which $T_w$ acting on $\ell_p(\Z)$ is
supercyclic and $I+T_w$ is hypercyclic is a dense $G_\delta$-subset
of the Banach space $c_0(\Z)$.
\end{remark}

It is easy to see that the dual of $T_w$ acting on $\ell_p(\Z)$ for
$1<p<\infty$ acts on $\ell_q(\Z)$ with $\frac1p+\frac1q=1$ according
to the formula
$$
T'_we_n=w_{n+1}e_{n+1}\qquad \text{for}\ \ n\in\Z.
$$
Considering the isometry $U$ on $\ell_q(\Z)$  defined by
$Ue_n=e_{-n}$ for each $n\in\Z$, we see that $U^{-1}T'_wU=T_{w'}$,
where $w'_n=w_{1-n}$ for any $n\in\Z$. Thus $T'_w$ is hypercyclic or
supercyclic if and only if $T_{w'}$ (acting on $\ell_q(\Z)$) is. In
the case $p=2$, the Hilbert space adjoint $T_w^\star$ acts on
$\ell_2(\Z)$ in a similar way $T_w^\star
e_n=\overline{w_{n+1}}e_{n+1}$ and is unitarily similar to $T'_w$
with a diagonal unitary operator providing the similarity. Thus the
cyclicity properties of $T'_w$ and $T_w^\star$ are the same.

\smallskip

Taking into account the fact that the map $w\mapsto w'$ is a
homeomorphism of $B_c$ onto itself for each $c>0$, we immediately
obtain the following corollary of Theorem~\ref{genshi}.

\begin{corollary}\label{c-genshi}For each $c>1$ the set of
$w\in B_c$ for which both $T_w$ and $T'_w$ acting on $\ell_2(\Z)$
are hypercyclic is a dense $G_\delta$-subset of $B_c$. For each
$c>0$ the set of $w\in B_c$ for which both $T_w$ and $T'_w$ are
supercyclic and both $I+T_w$ and $I+T'_w$ are hypercyclic is a dense
$G_\delta$-subset of of $B_c$.
\end{corollary}

Similarly the next corollary follows from Theorem~\ref{genshi1}.

\begin{corollary}\label{c-genshi1}Suppose that the space $E$ from
Theorem~$\ref{genshi1}$ satisfies the additional symmetry condition
that $(Jx)_n=x_{k-n}$ is an invertible continuous linear operator on
$E$ for some $k\in\Z$. Then set of $w\in E$ for which both $T_w$ and
$T'_w$ acting on $\ell_2(\Z)$ are supercyclic and both $I+T_w$ and
$I+T'_w$ are hypercyclic is a dense $G_\delta$-subset of $E$.
\end{corollary}

Applying the above corollary to  weighted $c_0$-spaces with
symmetric weight sequence yields the following result.

\begin{corollary}\label{coco}Let $\{a_n\}_{n\in\Z_+}$
be any sequence of positive numbers. Then there exists $w\in
c_0(\Z)$ such that $|w_n|\leq a_{|n|}$ for each $n\in\Z$, $T_w$ and
$T'_w$ acting on $\ell_2(\Z)$ are supercyclic and $I+T_w$ and
$I+T'_w$ are hypercyclic.
\end{corollary}

We conclude this section by proving Theorems~\ref{genshi}
and~\ref{genshi1}.

\begin{proof}[Proof of Theorem~$\ref{genshi}$]
It is straightforward to verify that the maps $w\mapsto T_w$ and
$w\mapsto I+T_w$ from $\Omega=B_c$ into $L(\ell_p(\Z))$ satisfy
(\ref{gen1}). Pick a non-empty open subset $U$ of $B_c$ and
non-empty open subsets $V$ and $W$ of $\ell_p(\Z)$.
\par\smallskip
{\bf Case} $c>1$: By definition of the topology of $B_c$, there
exist $w\in U$ and a positive integer $m$ such that $w_k=c$ for
$k>m$, $w_k=c^{-1}$ for $k<-m$ and $w_k\neq 0$ for $-m\leq k\leq m$.
According to Theorem~B, the bilateral weighted shift $T_w$ is
hypercyclic. Hence we can choose $x\in V$ and $n\in\Z_+$ such that
$T_w^nx\in W$. Thus $(w,x,T_w^nx)\in U\times V\times W$ and
therefore (\ref{gen2}) is satisfied. By Theorem~\ref{generic}, the
set of $w\in B_c$ for which $T_w$ is hypercyclic is a dense
$G_\delta$-subset of $B_c$.
\par\smallskip
{\bf Case} $c>0$: As above, there exist $w\in U$ and a positive
integer $m$ such that $w_k=c$ for $k>m$, $w_k=c/2$ for $k<-m$ and
$w_k\neq 0$ for $-m\leq k\leq m$. According to Theorem~B, the
bilateral weighted shift $T_w$ is supercyclic. Hence we can choose
$x\in V$, $n\in\Z_+$ and $\lambda\in\K$ such that $\lambda T_w^nx\in
W$. Thus $(w,x,\lambda T_w^nx)\in U\times V\times W$ and therefore
(\ref{gen3}) is satisfied. By Theorem~\ref{genericsup}, the set of
$w\in B_c$ for which $T_w$ is supercyclic is a dense
$G_\delta$-subset of $B_c$.
\par\smallskip
Finally, we can pick $m\in\Z$, $w\in U$, $x\in V$ and $y\in W$ such
that $w_m=0$, $w_k\neq 0$ for $k>m$, $x$ and $y$ have finite
supports and $x_k=y_k=0$ for $k<m$. It is straightforward to check
that
$$
x,y\in\spann\{e_m,e_{m+1},\dots\}\subseteq \KER T_w.
$$
By Corollary~\ref{generic2}, the set of $w\in B_c$ for which $I+T_w$
is hypercyclic is a dense $G_\delta$-subset of $B_c$.
\end{proof}

\begin{proof}[Proof of Theorem~$\ref{genshi1}$]
We use the same notation as in the proof of Theorem~\ref{genshi}.
Since $\phi(\Z)$ is dense in $\ell_p(\Z)$ and in $E$, we can choose
$w'\in U$, $x'\in V$, $y'\in W$ and a positive integer $m$ such that
$x'_k=y'_k=w'_k=0$ if $|k|>m$. Choosing positive numbers
$\epsilon_n$ for $n>-m$ small enough, we can ensure that the series
$\sum\limits_{n=m+1}^\infty \epsilon_ne_n$ converges in $E$ to
$w''\in E$, with $w=w'+w''\in U$ and $w_k=w'_k+w''_k\neq 0$ for
$k\geq -m$. Clearly $F=\{u\in\ell_p(\Z):u_k=0\ \text{for}\ k<-m\}$
is a closed invariant subspace of $T_w$ and $x,y\in F$. Moreover,
the restriction $R$ of $T_w$ to $F$ is isometrically similar to the
backward weighted  shift on $\ell_p$ with weight sequence
$\{w_{-m},w_{-m+1},\dots\}$ By the Hilden and Wallen theorem
\cite{hw} and the Salas theorem \cite{sal}, $R$ is supercyclic and
$I+R$ is hypercyclic. Since $x'\in V\cap F$ and $y'\in W\cap F$ we
see that the open subsets $V\cap F$ and $W\cap F$ of $F$ are
non-empty. Since $R$ is supercyclic, we can choose $x\in V\cap F$,
$n\in\Z_+$ and $\lambda\in\K$ such that $\lambda R^nx=\lambda
T_w^nx\in W\cap F$. Thus $(w,x,\lambda T_w^nx)\in U\times V\times W$
and therefore (\ref{gen3}) is satisfied. By
Theorem~\ref{genericsup}, the set of $w\in E$ for which $T_w$ is
supercyclic is a dense $G_\delta$-subset of $E$. Finally, since
$I+R$ is hypercyclic, we can choose $x\in V\cap F$ and $n\in\Z_+$
such that $(I+R)^nx=(I+T_w)^nx\in W\cap F$. Thus $(w,x,
(I+T_w)^nx)\in U\times V\times W$ and therefore (\ref{gen2}) is
satisfied for the map $w\mapsto I+T_w$. By Theorem~\ref{generic},
the set of $w\in E$ for which $I+T_w$ is hypercyclic is a dense
$G_\delta$-subset of $E$.
\end{proof}

\section{Mixing operators on spaces with weak topology \label{wea}}

In this section we shall prove Theorem~\ref{omegA}. In order to do
so we need a characterization of linear maps with no non-trivial
finite dimensional invariant subspaces. The underlying field plays
no role in this linear algebraic statement, so, for sake of
generality we formulate and prove it for linear maps on linear
spaces over an arbitrary field.

\subsection{Linear maps without finite dimensional invariant
subspaces}

Throughout this subsection $\kk$ is a field, $X$ is a linear space
over $\kk$, $T:X\to X$ is a $\kk$-linear map, $\pp =\kk[z]$ is the
space of polynomials on one variable, $\rr=\kk(z)$ is the space of
rational functions and $M$ is the operator on $\rr$ of
multiplication by the argument:
$$
M:\rr\to\rr,\quad Mf(z)=zf(z).
$$
We denote $\pp^*=\pp\setminus\{0\}$ and consider the degree function
$\deg:\rr\to\Z\cup\{-\infty\}$, extending the conventional degree of
a polynomial. We set $\deg(0)=-\infty$ and let $\deg(p/q)=\deg
p-\deg q$, where $p$ and $q$ are non-zero polynomials and the
degrees in the right hand side are the conventional degrees of
polynomials. Clearly this function is well-defined and is a grading
on $\rr$, that is, it satisfies the properties:
\begin{itemize}\itemsep=-2pt
\item[(d1)]$\deg(r_1r_2)=\deg(r_1)+\deg(r_2)$ for any
$r_1,r_2\in\rr$;
\item[(d2)]$\deg(r_1+r_2)\leq \max\{\deg r_1,\deg r_2\}$ for
any $r_1,r_2\in\rr$;
\item[(d3)]if $r_1,r_2\in\rr$ and $\deg r_1\neq \deg r_2$, then
$\deg(r_1+r_2)=\max\{\deg r_1,\deg r_2\}$.
\end{itemize}
Note that if $p$ is a non-zero polynomial, then $\deg p$ is the
usual degree of $p$.

As usual, a linear subspace $E$ of $X$ is called $T$-{\it invariant}
if $T(E)\subseteq E$ and it is called $T$-{\it biinvariant} if
$T(E)\subseteq E$ and $T^{-1}(E)\subseteq E$. The following lemma is
a key ingredient in the proof of Theorem~\ref{omegA}.

\begin{lemma} \label{empty1}
Let $T$ be a linear operator on a linear space $X$ with no
non-trivial finite dimensional invariant subspaces and let $L$ be a
finite dimensional subspace of $X$. Then there exists $n_0=n_0(L)\in
\N$ such that $p(T)(L)\cap L=\{0\}$ for any $p\in\pp$ with $\deg
p\geq n_0$.
\end{lemma}

In order to prove the above lemma, we need some preparation.

\begin{lemma} \label{inj} Let $T$ be a linear operator on a linear space
$X$. Then $T$ has no non-trivial finite dimensional invariant
subspaces if and only if $p(T)$ is injective for any non-zero
polynomial $p$.
\end{lemma}

\begin{proof} If $p$ is a non-zero polynomial and $p(T)$ is
non-injective, then there is non-zero $x\in X$ such that $p(T)x=0$.
Let $k=\deg p$. It is straightforward to verify that
$E=\spann\{x,Tx,\dots,T^{k-1}x\}$ is a non-trivial finite
dimensional invariant subspace for $T$. Assume now that $T$ has a
non-trivial finite dimensional invariant subspace $L$. Let $p$ be
the characteristic polynomial of the restriction of $T$ to $L$. By
the Hamilton--Cayley theorem $p(T)$ vanishes on $L$. Hence $p(T)$ is
non-injective.
\end{proof}

\begin{definition}\label{Tind}
For a linear operator $T$ on a linear space $X$ we say that vectors
$x_1,\dots,x_n$ in $X$ are $T$-{\it independent} if for any
polynomials $p_1,\dots,p_n$, the equality $\sum\limits_{j=1}^n
p_j(T)x_j=0$ implies $p_j=0$ for $1\leq j\leq n$. Otherwise, we say
that $x_1,\dots,x_n$ are $T$-{\it dependent}. A set $A\subset X$ is
called $T$-{\it independent} if any pairwise different vectors
$x_1,\dots,x_n\in A$ are $T$-{\it independent}.

For any non-zero $x\in X$, we define
\begin{equation}\label{xa}
F(T,x)=\{y\in X:\text{there are $p\in \pp^*$ and $q\in\pp$ such that
$p(T)y=q(T)x$}\}.
\end{equation}
\end{definition}

\begin{lemma} \label{empty00} Let $T$ be a linear operator on a linear space
$X$, $x\in X\setminus\{0\}$ and $F(T,x)$ be the space defined in
$(\ref{xa})$. Then $F(T,x)$ is a $T$-biinvariant linear subspace of
$X$.
\end{lemma}

\begin{proof} Let $y,u\in F(T,x)$ and $t,s\in \kk$. Then we can pick
$p_1,p_2\in\pp^*$ and $q_1,q_2\in \pp$ such that $p_1(T)y=q_1(T)x$
and $p_2(T)u=q_2(T)x$. Hence
$(p_1p_2)(T)(ty+su)=(tq_1p_2+sq_2p_1)(T)x$. Since $p_1$ and $p_2$
are non-zero, the polynomial $p_1p_2$ is also non-zero and we have
$ty+su\in F(T,x)$ and therefore $F(T,x)$ is a linear subspace of
$X$. Clearly $p_1(T)(Ty)=\widetilde q_1(T)x$, where $\widetilde
q_1(z)=zq(z)$. Hence $Ty\in F(T,x)$, which proves the $T$-invariance
of $F(T,x)$. Assume now that $w\in X$ and $Tw\in F(T,x)$. Thus we
can pick $p_3\in \pp^*$ and $q_3\in \pp$ such that
$p_3(T)Tw=q_3(T)x$. Hence $\widetilde p_3(T)w=q_3(T)x$, where
$\widetilde p_3(z)=zp_3(z)$, and therefore $w\in F(T,x)$. That is,
$F(T,x)$ is $T$-biinvariant.
\end{proof}

By the above lemma, we can consider linear operators, being
restrictions of $T$ to the invariant subspaces $F(T,x)$. The
following lemma describes these restrictions in the case when $T$
has no non-trivial finite dimensional invariant subspaces.

\begin{lemma} \label{empty0} Let $T$ be a linear operator on a linear space
$X$ with no non-trivial finite dimensional invariant subspaces. Let
also $x\in X\setminus\{0\}$ and $F(T,x)$ be defined in $(\ref{xa})$.
For each $y\in F(T,x)$ and any $p\in\pp^*$, $q\in\pp$ satisfying
$p(T)y=q(T)x$, we write $r_{x,y}=q/p$. Then the rational function
$r_{x,y}$ does not depend on the choice of $p\in\pp^*$ and $q\in\pp$
satisfying $p(T)y=q(T)x$, the map $S_xy=r_{x,y}$ from $F(T,x)$ to
$\rr$ is linear and $S_xTy=MS_xy$ for any $y\in F(T,x)$.
\end{lemma}

\begin{proof} Let $y\in F(T,x)$ and $p,p_1\in\pp^*$, $q,q_1\in\pp$
be such that $p(T)y=q(T)x$ and $p_1(T)y=q_1(T)x$. Hence
$(p_1p)(T)y=(q_1p)(T)x=(qp_1)(T)x$ and therefore
$(q_1p-qp_1)(T)x=0$. Since $x\neq 0$, Lemma~\ref{inj} implies that
$q_1p=qp_1$, or equivalently, $q/p=q_1/p_1$. Thus $q/p$ does not
depend on the choice of $p\in\pp^*$ and $q\in\pp$ satisfying
$p(T)y=q(T)x$. In particular, the map $y\mapsto r_{x,y}$ from
$F(T,x)$ to $\rr$ is well defined. Next, let $y,u\in F(T,x)$ and
$t,s\in \kk$. Then we can pick $p_1,p_2\in\pp^*$ and $q_1,q_2\in
\pp$ such that $p_1(T)y=q_1(T)x$ and $p_2(T)u=q_2(T)x$. Hence
$(p_1p_2)(T)(ty+su)=(tq_1p_2+sq_2p_1)(T)x$. It follows that
$$
r_{x,ty+su}=\frac{tq_1p_2+sq_2p_1}{p_1p_2}=t\frac{q_1}{p_1}+s\frac{q_2}{p_2}=
tr_{x,y}+sr_{x,u}
$$
and the linearity of the map $y\mapsto r_{x,y}$ is also verified. It
remains to show that $S_xTy=MS_xy$. Clearly $S_xy=q_1/p_1$. Since
$p_1(T)(Ty)=Tq_1(T)x$, we have $S_xTy(z)=zq_1(z)/p_1(z)$. Hence
$S_xTy=MS_xy$. \end{proof}

\begin{lemma}\label{empty2} Let $T$ be a linear operator with no non-trivial
finite dimensional invariant subspaces acting on a linear space $X$
and $A$ be a $T$-independent subset of $X$. Then for each $x\in A$,
$$
F(T,x)\cap \spann\left(\bigcup_{y\in A\setminus\{x\}}
F(T,y)\right)=\{0\},
$$
where the spaces $F(T,u)$ are defined in $(\ref{xa})$.
\end{lemma}

\begin{proof} Assume that the intersection in the above display
contains a non-zero vector $u$. Since $u\in F(T,x)$, there exist
$p\in\pp^*$ and $q\in\pp$ such that $p(T)u=q(T)x$. Since $u\neq 0$
and $p\neq 0$, according to Lemma~\ref{inj}, $p(T)u\neq 0$. It
follows that $q\neq 0$. On the other hand, since $u$ is a non-zero
element of the span of the union of $F(T,y)$ for $y\in
A\setminus\{x\}$, there exist pairwise different $x_1,\dots,x_n\in
A\setminus\{x\}$ and $u_j\in F(T,x_j)\setminus \{0\}$ such that
$u=u_1+{\dots}+u_n$. Pick $p_j\in\pp^*$ and $q_j\in\pp$ for which
$p_j(T)u_j=q_j(T)x_j$ for $1\leq j\leq n$. Since $u_j\neq 0$ and
$p_j\neq 0$, from Lemma~\ref{inj} it follows that $p_j(T)u_j\neq 0$
and therefore $q_j\neq 0$. Consider the polynomials $p_0=pp_1\dots
p_n$, $\widetilde p_0=p_0/p$, $\widetilde p_j=p_0/p_j$ for $1\leq
j\leq n$. Then
$$
p_0(T)u=(q\widetilde p_0)(T)x\ \ \text{and}\ \
p_0(T)u_j=(q_j\widetilde p_j)(T)x_j\ \ \text{for $1\leq j\leq n$}.
$$
Taking into account that $u=u_1+{\dots}+u_n$, we obtain
$$
(q\widetilde p_0)(T)x=\sum_{j=1}^n (q_j\widetilde p_j)(T)x_j.
$$
Since the polynomials $q\widetilde p_0$ and $q_j\widetilde p_j$ are
non-zero, the last display contradicts the $T$-independence of $A$,
since $x,x_1,\dots,x_n$ are pairwise different elements of $A$.
\end{proof}

\begin{proof}[Proof of Lemma~$\ref{empty1}$]
Clearly there exists in $L$ a maximal $T$-independent subset
$A=\{x_1,\dots,x_k\}$ (since $T$-independence implies linear
independence and $L$ is finite dimensional, $A$ is a finite set). It
follows from the maximality of $A$ that $L$ is contained in the sum
of $F(T,x_j)$ for $1\leq j\leq k$. The last sum is direct according
to Lemma~\ref{empty2}:
$$
L\subseteq N=\bigoplus_{j=1}^k F(T,x_j).
$$
Thus any $x\in N$ can be uniquely presented as a sum
$u^x_1+\dots+u^x_k$, where $u^x_j\in F(T,x_j)$. Using
Lemmas~\ref{empty0} and~\ref{empty00}, we can consider the linear
operator
$$
R:N\to \rr^k,\qquad R(x)=(R_1x,\dots,R_kx),\ \ \text{where}\ \
R_jx=S_{x_j}u_j^x\ \ \text{for}\ \ 1\leq j\leq k.
$$
According to Lemmas~\ref{empty00} and~\ref{empty0} we also have that
$N$ is a $T$-biinvariant subspace of $X$ and $R(Tx)_j=M(Rx)_j$ for
$1\leq j\leq k$. For each $x\in N$, let
$$
\delta(x)=\max_{1\leq j\leq k}\deg R_jx.
$$
Clearly $\delta(0)=-\infty$ and $\delta(x)\in\Z$ for each $x\in
N\setminus \{0\}$. Let also
$$
\Delta^+=\sup_{x\in L} \delta(x)\ \ \text{and}\ \
\Delta^-=\inf_{x\in L\setminus\{0\}} \delta(x).
$$
\par\smallskip
Using the fact that $L$ is finite dimensional, we will show that
$\Delta^+$ and $\Delta^-$ are finite. Indeed, assume that either
$\Delta^+=+\infty$ or $\Delta^-=-\infty$. Then there exists a
sequence $\{u_l\}_{l\in\N}$ of non-zero elements of $L$ such that
the sequence $\{\delta(u_l)\}_{l\in\N}$ is strictly monotonic. For
each $l$ we can pick $j(l)\in\{1,\dots,k\}$ such that
$\delta(u_l)=\deg R_{j(l)}u_l$. Then there is $\nu\in \{1,\dots,k\}$
such that the set $B_\nu=\{l\in\N:j(l)=\nu\}$ is infinite. It
follows that the degrees $\deg R_{\nu}u_l$ for $l\in B_\nu$ are
pairwise different. Property (d3) of the degree function implies
that the rational functions $R_{\nu}u_l$ for $l\in B_\nu$ are
linearly independent. Hence the infinite set $\{u_l:l\in B_\nu\}$ is
linearly independent in $X$, which is impossible since all $u_l$
belong to the finite dimensional space $L$. Thus $\Delta^+$ and
$\Delta^-$ are finite.

Now let $p\in\pp^*$ and $m=\deg p$. From (d1) and the equality
$R_jTx=MR_jx$, we immediately get that $\delta(p(T)x)=\delta(x)+m$
for each $x\in N$. Therefore, $\inf\bigl\{\delta(x):x\in
p(T)(L)\setminus\{0\}\bigr\}=\Delta^-+m$. In particular, if
$m>\Delta^+-\Delta^-$, then
$$
\inf_{x\in
p(T)(L)\setminus\{0\}}\delta(x)=\Delta^-+m>\Delta^+=\sup_{x\in L}
\delta(x).
$$
Thus $\delta(u)>\delta(v)$ for any non-zero $u\in P(T)(L)$ and $v\in
L$, which implies that $p(T)(L)\cap L=\{0\}$ whenever $\deg
p>\Delta^+-\Delta^-$. Thus the number $n_0=\Delta^+-\Delta^-+1$
satisfies the desired condition.
\end{proof}

\subsection{Proof of Theorem~\ref{omegA}}

The implications
$(\ref{omegA}.4)\Longrightarrow(\ref{omegA}.3)\Longrightarrow(\ref{omegA}.2)$
are trivial. Assume that $T$ is transitive and $T'$ has a
non-trivial finite dimensional invariant subspace. Then $T$ has a
non-trivial closed invariant subspace of finite codimension. Passing
to the quotient by this subspace, we obtain a transitive operator on
a finite dimensional topological vector space. Since there is only
one Hausdorff linear topology on a finite dimensional space, we
arrive to a transitive operator on a finite dimensional Banach
space. Since transitivity and hypercyclicity for operators on
separable Banach spaces are equivalent (see, for instance,
\cite{gri}), we obatin a hypercyclic operator on a finite
dimensional Banach space. On the other hand, it is well known that
such operators do not exist, see, for instance, \cite{ww}. This
proves the implication
$(\ref{omegA}.2)\Longrightarrow(\ref{omegA}.1)$. It remains to show
that (\ref{omegA}.1) implies (\ref{omegA}.4).

Assume that (\ref{omegA}.1) is satisfied and (\ref{omegA}.4) fails.
Then there exist non-empty open subsets $U$ and $V$ of $X$ and a
sequence $\{p_l\}_{l\in\Z_+}$ of polynomials such that $\deg
p_l\to\infty$ and $p_l(T)(V)\cap U=\varnothing$ for each $l\in\Z_+$.
Since $X$ carries weak topology, there exist two finite linearly
independent sets $\{f_1,\dots,f_n\}$ and $\{g_1,\dots,g_m\}$ in $X'$
and two vectors $(a_1,\dots,a_n)\in\K^n$ and
$(b_1,\dots,b_m)\in\K^m$ such that $U_0\subseteq U$ and
$V_0\subseteq V$, where
\begin{equation*}
U_0=\{u\in X:f_j(u)=a_j\ \ \text{for}\ \ 1\leq j\leq n\}\ \
\text{and}\ \ V_0=\{u\in X:g_j(u)=b_j\ \ \text{for}\ \ 1\leq j\leq
m\}.
\end{equation*}
Let $L=\spann\{f_1,\dots,f_n,g_1,\dots,g_m\}$. Since $T'$ has non
non-trivial finite dimensional invariant subspaces, by
Lemma~\ref{empty1}, for any sufficiently large $l$, $p_l(T')(L)\cap
L=\{0\}$. For such an $l$, the equality $p_l(T')(L)\cap L=\{0\}$
together with the injectivity of $p_l(T')$, provided by
Lemma~\ref{inj}, and the definition of $L$ imply that the vectors
$p_l(T')f_1,\dots,p_l(T')f_n,g_1,\dots,g_m$ are linearly
independent. Hence there exists a vector $u\in X$ such that
$$
\text{$p_l(T')f_j(u)=a_j$ \ for \ $1\leq j\leq n$\ \ and\ \
$g_j(u)=b_j$ \ for $1\leq j\leq m$.}
$$
Since $p_l(T')f_j(u)=f_j(p_l(T)u)$, the last display implies that
$u\in V_0\subseteq V$ and $p_l(T)u\in U_0\subseteq U$. Hence
$p_l(T)(V)\cap U$ contains $p_l(T)u$ and therefore is non-empty.
This contradiction completes the proof.

\section{Spaces without supercyclic semigroups $\{T_t\}_{t\in\R_+}$}

We shall prove Theorem~\ref{omeg} and show that on certain
topological vector spaces there are no strongly continuous
supercyclic semigroups $\{T_t\}_{t\in\R_+}$. In this section by the
dimension $\dim X$ of a vector space $X$ we mean its algebraic
dimension (=the cardinality of the Hamel basis). Symbol $\ccc$
stands for the cardinality of continuum: $\ccc=2^{\aleph_0}$. The
next theorem is the main result of this section.

\begin{theorem}\label{nosemi} Let $X$ be an infinite dimensional
locally convex space such that either in $X$ or in $X'_\sigma$ there
are no compact metrizable subsets whose linear span has dimension
$\ccc$. Then there are no strongly continuous supercyclic semigroups
$\{T_t\}_{t\in\R_+}$ on $X$.
\end{theorem}

The above theorem immediately implies the following corollary.

\begin{corollary}\label{nosemi1} Let $X$ be an infinite dimensional
locally convex space such that $\dim X<\ccc$ or $\dim X'<\ccc$. Then
there are no strongly continuous supercyclic semigroups
$\{T_t\}_{t\in\R_+}$ on $X$.
\end{corollary}

We prove Theorem~\ref{nosemi} at the end of this section. First, we
shall prove Theorem~\ref{omeg} by means of application of
Theorems~\ref{omegA} and~\ref{nosemi}.

\begin{proof}[Proof of Theorem~$\ref{omeg}$] Let $T\in L(\omega)$ be
such that $T'$ has no non-trivial finite dimensional invariant
subspaces and $\{p_l\}_{l\in\Z_+}$ be a sequence of polynomials such
that $\deg p_l\to\infty$ as $l\to\infty$. Since the topology of
$\omega$ is weak, Theorem~\ref{omegA} implies that for each
non-empty open subsets $U$ and $V$ of $\omega$, $p_l(U)\cap V\neq
\varnothing$ for all sufficiently large $l$. Hence
$\{(x,p_l(T)x):x\in\omega,\ l\in\Z_+\}$ is dense in
$\omega\times\omega$. By Theorem~U, $\{p_l(T):l\in \Z_+\}$ is
universal. It remains to show that there are no strongly continuous
supercyclic semigroups $\{T_t\}_{t\in\R_+}$ on $X$. Recall that
$\omega'=\phi$ and $\dim\phi=\aleph_0<\ccc$. Thus
Corollary~\ref{nosemi1} implies that there are no supercyclic
strongly continuous operator semigroups $\{T_t\}_{t\in\R_+}$ on
$\omega$.
\end{proof}

The rest of the section is devoted to the proof of
Theorem~\ref{nosemi}. We need some preparation.

\begin{lemma}\label{fds}Let $X$ be a finite dimensional topological
vector space of the real dimension $>2$. Then there is no
supercyclic strongly continuous operator semigroup
$\{T_t\}_{t\in\R_+}$ on $X$.
\end{lemma}

\begin{proof} As well-known, any strongly continuous operator
semigroup $\{T_t\}_{t\in\R_+}$ on $\K^n$ has shape
$\{e^{tA}\}_{t\in\R_+}$, where $A\in L(\K^n)$. Assume the contrary.
Then there exists $n\in\N$ and $A\in L(\K^n)$ such that the
semigroup $\{e^{tA}\}_{t\in\R_+}$ is supercyclic and $n\geq 3$ if
$\K=\R$, $n\geq 2$ if $\K=\C$. Since the operators $e^{tA}$ are
invertible and commute with each other, Proposition~G implies that
the set $W$ of universal elements for the family $\{ze^{tA}:z\in\K,\
t\in\R_+\}$ is dense in $\K^n$. On the other hand, for each $c>0$
and any $x\in\K^n$, from the restrictions on $n$ it follows that the
closed set $\{ze^{tA}x:z\in\K,\ 0\leq t\leq c\}$ is nowhere dense in
$\K^n$ (one can use smoothness of the map $(z,t)\mapsto ze^{tA}x$ to
see that the topological dimension of $\{ze^{tA}x:z\in\K,\ 0\leq
t\leq c\}$ is less than that of $\K^n$). Hence, each $x\in W$ is
universal for the family $\{ze^{tA}:z\in\K,\ t>c\}$ for any $c>0$.
Now if $(a,b)$ is a finite subinterval of $(0,\infty)$, it is easy
to see that the family $\{ze^{tkA}:z\in\K,\ a<t<b,\ k\in\Z_+\}$
contains $\{ze^{tA}:z\in\K,\ t>c\}$ for a sufficiently large $c>0$.
It follows that for each $x\in W$, the set $\{ze^{tkA}x:z\in\K,\
a<t<b,\ k\in\Z_+\}$ is dense in $\K^n$. Taking into account that
$(a,b)$ is arbitrary and $W$ is dense in $\K^n$, we see that
$\{(t,x,ze^{tkA}x:t\in\R_+,\ z\in\K,\ x\in\R^n, k\in\Z_+\}$ is dense
in $\R_+\times\K^n\times\K^n$. Applying Theorem~\ref{genericsup}, we
see that for a generic $t\in\R_+$ in the Baire category sense the
operator $e^{tA}$ is supercyclic. This contradicts the well-known
fact that there are no supercyclic operators on finite dimensional
spaces of real dimension $>2$.
\end{proof}

The following lemma appears as Lemma~2 in \cite{shka1}. It is worth
noting that under the Continuum Hypothesis its statement becomes
trivial.

\begin{lemma}\label{ch} Let $(M,d)$ be a separable complete metric space
and $X$ be a topological vector space. Then for any continuous map
$f:M\to X$, the algebraic dimension of $\spann f(M)$ is either
finite or countable or continuum.
\end{lemma}

We use the above lemma to prove the following dichotomy.

\begin{lemma}\label{dich} Let $\{T_t\}_{t\in\R_+}$ be a strongly
continuous operator semigroup on a topological vector space $X$ and
$x\in X$. Then the space $C(x)=\spann\{T_tx:t\in\R_+\}$ is either
finite dimensional or has dimension $\ccc$.
\end{lemma}

\begin{proof} From Lemma~\ref{ch} it follows that $\dim C(x)$ is
either finite or $\aleph_0$ or $\ccc$. Thus it suffices to rule out
the case $\dim C(x)=\aleph_0$.

Assume that $\dim C(x)=\aleph_0$. Restricting the operators $T_t$ to
the invariant subspace $C(x)$, we can without loss of generality
assume that $C(x)=X$. Thus $\dim X=\aleph_0$ and therefore $X$ is
the union of an increasing sequence $\{E_n\}_{n\in\Z_+}$ of finite
dimensional subspaces. For each $\epsilon>0$ let
$X_\epsilon=\spann\{T_tx:t\geq\epsilon\}$. First, we shall show that
each $X_\epsilon$ is finite dimensional. Let $\epsilon>0$,
$0<\alpha<\epsilon$ and  $A_n=\{t\in[\alpha,\epsilon]:T_t x\in
E_n\}$ for $n\in\Z_+$. Clearly $A_n$ are closed subsets of the
interval $[\alpha,\epsilon]$ and
$[\alpha,\epsilon]=\bigcup\limits_{n=0}^\infty A_n$ since $X$ is the
union of $E_n$. By the Baire category theorem there is $n\in\Z_+$
such that $A_n$ has non-empty interior in $[\alpha,\epsilon]$. Hence
we can pick $a,b\in\R$ such that $\alpha\leq a<b\leq\epsilon$ and
$T_tx\in E_n$ for any $t\in [a,b]$. We shall show that $X_a\subseteq
E_n$. Assume, it is not the case. Then the number
$c=\inf\{t\in[a,\infty):T_tx\notin E_n\}$ belongs to $[b,\infty)$.
Since the set $\{t\in\R_+:T_t\in E_n\}$ is closed, $T_cx\in E_n$.
Since $[a,b]$ is uncountable and the span of $\{T_t:t\in[a,b]\}$ is
finite dimensional, we can pick $a\leq t_0<t_1<{\dots}<t_n\leq b$
and $c_1,\dots,c_{n-1}\in\K$ such that
\begin{equation}\label{ttn}
T_{t_n}x=\sum_{j=1}^{n-1} c_jT_{t_j}x.
\end{equation}
Since $T_cx\in E_n$, by definition of $c$, we can pick $t\in
(c,c+t_n-t_{n-1})$ such that $T_tx\notin E_n$. Since $t>c\geq t_n$,
formula (\ref{ttn}) implies that
$$
T_tx=T_{t-t_n}T_{t_n}x=T_{t-t_n}\sum_{j=1}^{n-1}
c_jT_{t_j}x=\sum_{j=1}^{n-1} c_jT_{t-t_n+t_j}x\in E_n.
$$
because $a\leq t-t_n+t_j\leq c$ for $1\leq j\leq n-1$. This
contradiction proves that $X_\epsilon\subseteq X_a\subseteq E_n$.
Thus $X_\epsilon$ is finite dimensional for each $\epsilon>0$. Since
$T_t(X)=T_t(C(x))\subseteq X_t$, it follows that $T_t$ has finite
rank for any $t>0$.

Now assume that $t>0$. Since $T_t$ has finite rank, $F_t=\ker T_t$
is a closed subspace of $X$ of finite codimension. It is also clear
that $F_t$ is $T_s$-invariant for each $s\in\R_+$. Passing to
quotient operators, $S_s\in L(X/F_t)$, $S_s(u+F_t)=T_su+F_t$, we
arrive to a strongly continuous semigroup $\{S_s\}_{s\in\R_+}$ on
the finite dimensional space $X/F_t$. Hence there is $A\in L(X/F_t)$
such that $S_s=e^{sA}$ for any $s\in\R_+$. Thus each $S_s$ is
invertible. Since each $S_s$ is a quotient of $T_s$, we see that
$$
{\tt rk}\,T_s\geq {\tt rk}\,S_s=\dim X/F_t={\tt rk}\,T_t\ \ \
\text{for any $t>0$ and $s\geq 0$}.
$$
It follows that $T_t$ for $t>0$ have the same rank $k\in\N$. Passing
to the limit as $t\to 0$, we see that the identity operator $I=T_0$
is the strong operator topology limit of a sequence of rank $k$
operators. Hence ${\tt rk}\,I\leq k$. That is, $X$ is finite
dimensional. This contradiction completes the proof.
\end{proof}

\begin{lemma}\label{phi} Let $X$ be a topological vector space
in which the linear span of each metrizable compact subset has
dimension $<\ccc$. Then for any strongly continuous operator
semigroup $\{T_t\}_{t\in\R_+}$ on $X$ and any $x\in X$, the space
$C(x)=\spann\{T_tx:t\in\R_+\}$ is finite dimensional.
\end{lemma}

\begin{proof} Let $\{T_t\}_{t\in\R_+}$ be a strongly continuous
operator semigroup on $X$ and $x\in X$. Strong continuity of
$\{T_t\}_{t\in\R_+}$ implies that for any $n\in\N$, the set
$K_s=\{T_tx:0\leq t\leq n\}$ is compact and metrizable. Hence $\dim
E_n<\ccc$ for any $n\in\N$, where $E_n=\spann(K_n)$. Since the sum
of countably many cardinals strictly less than $\ccc$ is strictly
less than $\ccc$, we see that
$$
\dim C(x)=\dim \biggl(\bigcup_{n=1}^\infty E_n\biggr)\leq
\sum_{n=1}^\infty \dim E_n<\ccc.
$$
By Lemma~\ref{dich}, $C(x)$ is finite dimensional.
\end{proof}

\begin{proposition}\label{gege} Let $X$ be an infinite dimensional
locally convex space such that in $X'_\sigma$ the span of any
compact metrizable subset has dimension $<\ccc$. Then there is no
strongly continuous supercyclic operator semigroup
$\{T_t\}_{t\in\R_+}$ on $X$.
\end{proposition}

\begin{proof} Assume that there exists a supercyclic strongly continuous
operator semigroup $\{T_t\}_{t\in\R_+}$ on $X$. It is
straightforward to verify that $\{T'_t\}_{t\in\R_+}$ is a strongly
continuous semigroup on $X'_\sigma$. Pick three linearly independent
vectors $f_1$, $f_2$ and $f_3$ in $X'$. By Lemma~\ref{phi},
$E_j=\spann\{T'_t f_j:t\in\R_+\}$ is finite dimensional for $1\leq
j\leq 3$. Clearly each $E_j$ is $T'_t$-invariant for any $t\in\R_+$.
Then $E=E_1+E_2+E_3$ is finite dimensional and $T'_t$-invariant for
any $t\in\R_+$. Since $f_j\in E$ for $1\leq j\leq 3$, $\dim E\geq
3$. Since $E$ is $T'_t$ invariant, we see that its annihilator
$F=\{x\in X:f(x)=0\ \ \text{for any}\ \ f\in E\}$ if $T_t$-invariant
for each $t\in\R_+$. Thus we can consider the quotient operators
$S_t\in L(X/F)$, $S_t (x+F)=T_tx+F$. Clearly $\{S_t\}_{t\in\R_+}$ is
a strongly continuous operator semigroup on $X/F$. Moreover,
$\{S_t\}_{t\in\R_+}$ is supercyclic since $\{T_t\}_{t\in\R_+}$ is.
Now since $\dim E=\dim X/F$, $X/F$ is finite dimensional and has
dimension $\geq3$. By Lemma~\ref{fds} there are no strongly
continuous supercyclic operator semigroups on $X/F$. This
contradiction completes the proof.
\end{proof}

\begin{proof}[Proof of Theorem~$\ref{nosemi}$] If $X$ has no compact
metrizable subsets whose linear span has dimension $\ccc$,
Lemma~\ref{phi} implies that the linear span of any orbit
$\{T_tx:t\in\R_+\}$ is finite dimensional. It follows that
$\{T_t\}_{t\in\R_+}$ is not supercyclic. It remains to consider the
case when $X'_\sigma$ has no compact metrizable subsets whose linear
span has dimension $\ccc$ and apply Proposition~\ref{gege}.
\end{proof}

\section{The space $\phi$ \label{P}}

Recall that $\phi$ is a linear space of countable algebraic
dimension carrying the strongest locally convex topology. In this
section we mention certain properties of $\phi$, mainly those which
are related to continuous linear operators. It is well known
\cite{shifer} that $\phi$ is complete and all linear subspaces of
$\phi$ are closed. Moreover, infinite dimensional subspaces of
$\phi$ are isomorphic to $\phi$. It is also well-known that for any
topology $\theta$ on $\phi$ such that $(\phi,\theta)$ is a
topological vector space, $\theta$ is weaker than the original
topology of $\phi$. The latter observation immediately implies the
following lemma.

\begin{lemma}\label{phi1} For any topological vector space $X$ and
any linear map $T:\phi\to X$, $T$ is continuous.
\end{lemma}

\begin{lemma}\label{phi2} Let $X$ be a topological vector space and
$T:X\to\phi$ be a surjective continuous linear operator. Then $X$ is
isomorphic to $\phi\times\ker T$.
\end{lemma}

\begin{proof} Since $T$ is linear and surjective, there exists a
linear map $S:\phi\to X$ such that $TS=I$. By Lemma~\ref{phi1}, $S$
is continuous. Consider the linear maps
$$
A:\phi\times\ker T\to X,\quad A(u,y)=y+Su\ \ \text{and}\ \ B:X\to
\phi\times\ker T,\quad Bx=(Tx,x-STx).
$$
It is easy to see that $A$ and $B$ are continuous and that $AB=I$
and $BA=I$. Hence $B$ is a required isomorphism.
\end{proof}

\begin{corollary}\label{phi3} Let $X$ be a topological vector space.
Then the following conditions are equivalent.
\begin{itemize}\itemsep=-2pt
\item[{\rm (\ref{phi3}.1)}]$X$ is isomorphic to a space of the shape $Y\times\phi$,
where $Y$ is a topological vector space$;$
\item[{\rm (\ref{phi3}.2)}]$X$ has a quotient isomorphic to $\phi;$
\item[{\rm (\ref{phi3}.3)}]there is $T\in L(X,\phi)$ such that $T(X)$ is
infinite dimensional.
\end{itemize}
\end{corollary}

\begin{proof} The implications
$(\ref{phi3}.1)\Longrightarrow(\ref{phi3}.2)\Longrightarrow(\ref{phi3}.3)$
are trivial. Assume that (\ref{phi3}.3) is satisfied. That is, there
is $T\in L(X,\phi)$ with infinite dimensional $T(X)$. Since any
infinite dimensional linear subspace of $\phi$ is isomorphic to
$\phi$, we see that $T(X)$ is isomorphic to $\phi$. Hence there is a
surjective $S\in L(X,\phi)$. By Lemma~\ref{phi2} $X$ is isomorphic
to $Y\times\phi$, where $Y=\ker S$. Hence (\ref{phi3}.3) implies
(\ref{phi3}.1), which completes the proof.
\end{proof}

\subsection{Cyclic operators on $\phi$}

Clearly $\phi$ is isomorphic to the space $\pp$ of all polynomials
over $\K$ endowed with the strongest locally convex topology. The
shift operator on $\phi$ is obviously similar to the operator
\begin{equation}\label{M}
M:\pp\to\pp,\quad Mp(z)=zp(z).
\end{equation}
For each $n\in\N$ we denote
\begin{equation}\label{pn}
\pp_n=\{p\in\pp:\deg p<n\}.
\end{equation}
Clearly $\pp_n$ is an $n$-dimensional subspace of $\pp$.

\begin{lemma}\label{cphi1} An operator $T\in L(\phi)$ is cyclic if
and only if $T$ is similar to the operator $M$.
\end{lemma}

\begin{proof} Clearly $1$ is a cyclic vector for $M$. Hence any
operator similar to $M$ is cyclic. Now let $T\in L(\phi)$ be cyclic
and $x\in\phi$ be a cyclic vector for $T$. Then the vectors $T^nx$
for $n\in\Z_+$ are linearly independent. Indeed, otherwise their
span is finite dimensional, which contradicts cyclicity of $x$ for
$T$. Since any linear subspace of $\phi$ is closed, we see that
$\{T^nx:n\in\Z_+\}$ is an algebraic basis of $\phi$. It is easy to
see then that the linear map $J:\pp\to\phi$, $Jp=p(T)x$ is
invertible and $T=JMJ^{-1}$. By Lemma~\ref{phi1}, $J$ and $J^{-1}$
are continuous. Hence $T$ is similar to $M$.
\end{proof}

We need a multicyclic version of the above lemma.

\begin{lemma}\label{cphi2} Let $T\in L(\phi)$. Then the following
conditions are equivalent.
\begin{itemize}\itemsep=-2pt
\item[{\rm (\ref{cphi2}.1)}]$T$ is multicyclic$;$
\item[{\rm (\ref{cphi2}.2)}]there exists $k\in\N$ and a linear
subspace $Y$ of $\phi$ of finite codimension such that
$T(Y)\subseteq Y$ and the restriction $T\bigr|_Y\in L(Y)$ is similar
to $M^k$, where $M$ is the operator defined in $(\ref{M})$.
\end{itemize}
\end{lemma}

\begin{proof} First, assume that (\ref{cphi2}.2) is satisfied. Pick
a finite dimensional subspace $Z$ of $\phi$ such that $\phi=Z\oplus
Y$. Since $T\bigr|_Y$ is similar to $M^k$, we can pick an invertible
linear operator $J:\pp\to\phi$ such that $Ty=JM^kJ^{-1}y$ for any
$y\in Y$. Let $L=Z+J(\pp_k)$, where $\pp_k$ is defined in
(\ref{pn}). Clearly $L$ is a finite dimensional subspace of $\phi$.
From the equality $T\bigr|_Y=JM^kJ^{-1}$ it easily follows that for
any $n\in\N$,
$$
L+T(L)+{\dots}+T^{n-1}(L)\supseteq Z+J(\pp_{nk}).
$$
Hence the linear span of the union of $T^j(L)$ for $j\in\Z_+$
contains $Z+J(\pp)=Z+Y=\phi$. Thus $T$ is $m$-cyclic with $m=\dim
L$. The implication $(\ref{cphi2}.2)\Longrightarrow(\ref{cphi2}.1)$
has been verified.

Assume now that $T$ is $n$-cyclic for some $n\in\N$. Then there is
an $n$-dimensional subspace $L$ of $\phi$ such that
\begin{equation}\label{multi}
\spann\biggl(\bigcup_{k=0}^\infty T^n(L)\biggr)=\phi
\end{equation}
(again, we use the fact that any linear subspace of $\phi$ is closed
and therefore a dense subspace of $\phi$ must coincide with $\phi$).
We will use the concept of $T$-independence, introduced in
Section~\ref{wea}. Since $T$-independence implies linear
independence, any $T$-independent subset of $L$ has at most $n$
elements. Let $k$ be the maximum of cardinalities of $T$-independent
subsets of $L$ and $A$ be a $T$-independent subset of cardinality
$k$. Since $A$ is linearly independent, we can pick a subset
$B\subset L$ of cardinality $n-k$ such that $A\cup B$ is a basis in
$L$. From the definition of $k$ it follows that for any $b\in B$,
$A\cup\{b\}$ is not $T$ independent and therefore using
$T$-independence of $A$, we can find polynomials $p_b$ and $p_{b,a}$
for $a\in A$ such that
\begin{equation}\label{tind}
p_b\neq 0\quad\text{and}\quad p_b(T)b=\sum_{a\in A}p_{b,a}(T)a.
\end{equation}
Now let $m=0$ if $B=\varnothing$ and  $m=\max\limits_{b\in B}\deg
p_b$ otherwise. Consider the spaces
$$
Z=\spann\{T^jb:b\in B,\ 0\leq j\leq m\}\ \ \text{and}\ \
Y=\spann\{T^ja:a\in A,\ j\in\Z_+\}.
$$
Then $Z$ is finite dimensional and $T(Y)\subseteq Y$. Obviously,
\begin{equation}\label{tind0}
T^ja\in Y\subseteq Y+Z\quad\text{for any $a\in A$ and $j\in\Z_+$}.
\end{equation}
Let $b\in B$ and $j\in\Z_+$. Since $p_b\neq 0$ and $\deg p_b\leq m$,
we can find polynomials $q,r$ such that $\deg r<m$ and
$t^j=q(t)p_b(t)+r(t)$. Then $T^jb=q(T)p_b(T)b+r(T)b$. Since $\deg
r<m$, $r(T)b\in Z$. According to (\ref{tind}), $p_b(T)b\in Y$. Since
$Y$ is invariant for $T$, $q(T)p_b(T)b\in Y$. Thus
\begin{equation}\label{tind1}
T^jb\in Y+Z\quad\text{for any $b\in B$ and $j\in\Z_+$}.
\end{equation}
Since $A\cup B$ is a basis of $L$, from (\ref{tind0}) and
(\ref{tind1}) it follows that $T^j(L)\subseteq Y+Z$ for each
$j\in\Z_+$. According to (\ref{multi}), $\phi=Y+Z$. Since $Z$ is
finite dimensional, we see that $Y$ has finite codimension in
$\phi$. In particular, $Y$ is non-trivial and therefore
$A\neq\varnothing$. That is, $1\leq k\leq n$ and
$A=\{a_1,\dots,a_k\}$. Now consider the linear operator $J:\pp\to
Y$, which sends the monomial $t^{l}$ to $T^{j}a_s$, where $j\in
\Z_+$ and $s\in \{1,\dots,k\}$ are uniquely defined by the equation
$l+1=jk+s$. By definition of $Y$, $J$ is onto. From $T$-independence
of $A$ it follows that $J$ is also one-to-one. By definition of $J$,
we have $Jt^{l+k}=TJt^l$. Hence $JM^kp=TJp$ for any $p\in\pp$. That
is, $M^k$ and $T\bigr|_{Y}$ are similar.
\end{proof}

\begin{corollary}\label{dran} Let $T$ be a multicyclic operator on
$\phi$. Then $T$ is not onto.
\end{corollary}

\begin{proof} According to Lemma~\ref{cphi2}, we can decompose $\phi$ into
a direct sum $\phi=Y\oplus Z$, where $Z$ has finite dimension
$m\in\Z_+$, $T(Y)\subseteq Y$ and $T\bigr|_Y$ is similar to $M^k$
for some $k\in\N$. Since $T\bigr|_Y$ is similar to $M^k$,
$T^{m+1}(Y)$ has codimension $k(m+1)>m$ in $Y$. Hence $\dim
\phi/T^{m+1}(Y)>m$. On the other hand, $\dim T^{m+1}(Z)\leq \dim
Z=m$. Thus $T^{m+1}(\phi)=T^{m+1}(Z)+T^{m+1}(Y)$ has positive
codimension in $\phi$. Hence $T^{m+1}$ is not onto and so is $T$.
\end{proof}

\subsection{Proof of Theorem~\ref{timesphi}}

In this section $X$ is a topological vector space, which has no
quotient isomorphic to $\phi$. We have to show that there are no
cyclic operators with dense range on $X\times \phi$. Assume the
contrary and let $T\in L(X\times\phi)$ be a cyclic operator with
dense range. Consider the matrix representation of $T$:
$$
T=\begin{pmatrix}A&B\\ C&D\end{pmatrix}, \ \ \text{where}\ \ A\in
L(X),\ B\in L(\phi,X),\ C\in L(X,\phi)\ \text{and}\ D\in L(\phi)
$$
with $T$ acting according to the formula $T(x,u)=(Ax+Bu,Cx+Du)$.
Since $T$ is cyclic, we can pick a vector $(x,u)\in X\times \phi$
such that $E=\spann\{T^k(x,u):k\in\Z_+\}$ is dense in $X\times\phi$.
Since $T$ has dense range, then $T^m$ has dense range for any
$m\in\Z_+$. Thus $E_m=T^m(E)=\spann\{T^k(x,u):k\geq m\}$ is dense in
$X\times\phi$ for each $m\in\Z_+$. Let $T^k(x,u)=(x_k,u_k)$ for
$k\in\Z_+$, where $x_k\in X$ and $u_k\in\phi$. Since
$E_m=\spann\{(x_k,u_k):k\geq m\}$ is dense in $X\times\phi$, we see
that $F_m=\spann\{u_k:k\geq m\}$ is dense in $\phi$ for any
$m\in\Z_+$. Hence $F_m=\phi$ for any $m\in\Z_+$. Since $X$ has no
quotients isomorphic to $\phi$, Lemma~\ref{phi3} implies that
$L_0=C(X)$ is a finite dimensional subspace of $\phi$. Then the
space
$$
L=\spann\bigl(L_0\cup \{u\}\bigr).
$$
is also finite dimensional. Clearly $u_0=u\in L$ and
$u_{k+1}=Cx_{k}+Du_{k}\in Du_{k}+L$ for any $k\in\Z_+$. It follows
that each $u_{k}$ belongs to the space spanned by the union of
$D^m(L)$ for $m\in\Z_+$. Since $L$ is finite dimensional, $D$ is
multicyclic. By Lemma~\ref{cphi2}, we can decompose $\phi$ into a
direct sum $\phi=Y\oplus Z$, where $Z$ is finite dimensional,
$D(Y)\subseteq Y$ and $D\bigr|_Y$ is similar to $M^n$ for some
$n\in\N$. That is, there exists an invertible $J\in L(Y,\pp)$ such
that $D\bigr|_Y=J^{-1}M^nJ$. Let also $P\in L(\phi)$ be the linear
projection onto $Y$ along $Z$. We consider two cases.

{\bf Case 1. } \ The sequence $\{\deg JPu_k\}$ is bounded from
above. In this case $\spann\{JPu_k:k\in\Z_+\}$ is finite
dimensional. Since $J$ is invertible, $\spann\{Pu_k:k\in\Z_+\}$ is
finite dimensional. Since $P$ has finite dimensional kernel,
$F_0=\spann\{u_k:k\in\Z_+\}$ is finite dimensional. We have arrived
to a contradiction with the equality $F_0=\phi$.

{\bf Case 2. } \ The sequence $\{\deg JPu_k\}$ is unbounded from
above. Since $N=(L+Z+D(Z))\cap Y$ is finite dimensional,
$$
m=\sup_{w\in N\setminus\{0\}} \deg Jw\in \Z_+.
$$
We shall show that $\deg JPu_{k+1}=n+\deg JPu_{k}$ whenever $\deg
JPu_{k}>m$. Indeed, let $k\in \Z_+$ be such that $\deg JPu_{k}>m$.
By definition of $P$, $u_k-Pu_k\in Z$ and $u_{k+1}-Pu_{k+1}\in Z$.
As we know, $u_{k+1}\in Du_{k}+L$. Hence $Pu_{k+1}\in
DPu_k+L+Z+D(Z)$. Since $Pu_{k+1}$ and $DPu_k$ belong to $Y$, we have
$Pu_{k+1}\in DPu_k+N$. Thus there is $w\in N$ such that
$Pu_{k+1}=DPu_k+w$. Hence $JPu_{k+1}=JDPu_k+Jw=M^nJPu_k+Jw$. Since
$\deg M^nJPu_k=n+\deg JPu_k>m\geq \deg Jw$, we have $\deg
JPu_{k+1}=n+\deg JPu_k$ (the degree of the sum of two polynomials of
different degrees equals to the maximum of the degrees). Since the
sequence $\{\deg JPu_k\}$ is unbounded from above, there is $k\in
\Z_+$ such that $\deg JPu_k>m$ and according to the just proven
statement, we will have $\deg JPu_j=\deg JPu_k+n(j-k)$ for $j\geq
k$. Since any family of polynomials with pairwise different degrees
is linearly independent, we see that the vectors $JPu_j$ for $j\geq
k$ are linearly independent. Since $JP$ is a linear operator, the
vectors $u_j$ for $j\geq k$ are linearly independent in $\phi$.
Hence the sequence of spaces $F_j$ for $j\geq k$ is strictly
decreasing. On the other hand, we know that $F_j=\phi$ for each
$j\in\Z_+$. This contradiction completes the proof.

\section{Hypercyclicity of operators on direct sums}

We shall prove the following lemma, which is a key for the proof of
Theorem~\ref{sumfre}.

\begin{lemma}\label{sumlcs} Let $\{X_n\}_{n\in\Z_+}$ be a sequence
of infinite dimensional locally convex spaces such that
\begin{itemize}\itemsep=-2pt
\item[{\rm(\ref{sumlcs}.1)}]there exists a dense linear subspace $Y$
of $X_0$, carrying a topology, stronger than the one inherited from
$X_0$ and turning $Y$ into a separable metrizable topological vector
space$;$
\item[{\rm(\ref{sumlcs}.2)}]there exists $T_0\in L(X_0,X_0\oplus
X_1)$ with dense range$;$
\item[{\rm(\ref{sumlcs}.3)}]for each $n\in\N$, there exists $T_n\in
L(X_n,X_{n+1}\times\K)$ with dense range.
\end{itemize}
Then there is a hypercyclic operator $S$ on
$X=\bigoplus\limits_{n=0}^\infty X_n$. \end{lemma}

It is worth noting that condition (\ref{sumlcs}.1) implies that
$X_0$ is separable, condition (\ref{sumlcs}.2) implies that $X_1$ is
separable and condition (\ref{sumlcs}.3) implies that $X_n$ for
$n\geq 2$ are all separable. Thus $X$ is separable. We need the
following auxiliary lemma.

\begin{lemma}\label{aux1} Let $X$ and $Y$ be topological vector spaces
such that there exists $T\in L(X,Y\times\K)$ with dense range. Then
for any closed hyperplane $H$ of $X$, there exists $S\in L(H,Y)$
with dense range.
\end{lemma}

\begin{proof} Let $T_0\in L(H,Y\times \K)$  be the restriction
of $T$ to $H$. We can write $T_0=(S_0,g)$, where $S_0\in L(H,Y)$ and
$g\in H'$. If $T_0$ has dense range, then $S=S_0$  is a continuous
linear operator from $H$ to $Y$ with dense range. It remains to
consider the case when the range of $T_0$ is not dense. Since the
range of $T$ is dense and $T_0$ is a restriction of $T$ to a closed
hyperplane, the codimension of the closure of $T_0(H)$ in $Y\times
\K$ does not exceed 1. Thus the codimension of the closure of
$T_0(H)$ in $Y\times \K$ is exactly $1$ and there is a non-zero
$\psi\in (Y\times\K)'$ such that $T_0(H)$ is a dense subspace of
$\ker\psi$. If $\ker\psi=Y\times\{0\}$, then again we can take
$S=S_0$. If $\ker\psi\neq Y\times\{0\}$, we can pick $y\in Y$ such
that $\psi(y)=1$. Take $S\in L(H,Y)$, $Sx=S_0x+g(x)y$. It is
straightforward to verify that $S$ has dense range.
\end{proof}

\begin{proof}[Proof of Lemma~$\ref{sumlcs}$] Let $\{U_n\}_{n\in\Z_+}$
be a base of topology of $Y$ and
$$
Z_n=\{x\in X:x_j=0\ \ \text{for}\ \ j>n\}.
$$
Clearly $X$ is the union of the increasing sequence of subspaces
$Z_n$ and each $Z_n$ is naturally isomorphic to
$\bigoplus\limits_{k=0}^n X_k$. We shall construct inductively a
sequence of operators $S_k\in L(Z_k,Z_{k+1})$ and vectors $y_k\in Y$
satisfying the following conditions for any $k\in\Z_+$:
\begin{itemize}\itemsep=-2pt
\item[(a1)]$S_jx=S_kx$ for any $j<k$ and $x\in Z_j$;
\item[(a2)]$S_k(Z_{k})$ is dense in $Z_{k+1}$;
\item[(a3)]$S_k\dots S_0y_k\notin Z_k$;
\item[(a4)]$S_k\dots S_0y_{k-1}=y_k$ if $k\geq 1$;
\item[(a5)]$y_k\in U_k$.
\end{itemize}
By (\ref{sumlcs}.2) there exists $S_0\in L(Z_0,Z_1)$ with dense
range. Since $Y$ is dense in $X_0=Z_0$, $S_0$ has dense range and
$Z_0$ is nowhere dense in $Z_1$, we can pick $y_0\in U_0$ such that
$S_0y_0\notin Z_0$. The basis of induction has been constructed.
Assume now that $n\in\N$ and $y_k\in Y$, $S_k\in L(Z_k,Z_{k+1})$
satisfying (a1--a5) for $k<n$ are already constructed. According to
(a3) for $k=n-1$, we have $w=S_{n-1}\dots S_0y_{n-1}\notin Z_{n-1}$.
That is, the $n^{\rm th}$ component $w_n$ of $w$ is non-zero. Since
$X_n$ is locally convex, we can pick a closed hyperplane $H$ in
$X_n$ such that $w_n\notin H$. Let $P$ be the linear projection on
$Z_n$ onto $H$ along $Z_{n-1}\oplus \spann\{w_n\}$. From
(\ref{sumlcs}.3) and Lemma~\ref{aux1} it follows that there is $R\in
L(H,X_{n+1})$ with dense range. According to (a3) for $k=n-1$, the
operator $S_{n-1}\dots S_0$ from $Z_0$ to $Z_n$ has dense range.
Hence the operator $Q=RPS_{n-1}\dots S_0$ from $Z_0$ to $X_{n+1}$
has dense range. Since $Y$ is dense in $Z_0$, we can pick $y_n\in
U_n$ such that $Qy_n\neq 0$. Now we define the operator $S_n:Z_n\to
Z_{n+1}$. It is easy to see that $Z_n=Z_{n-1}\oplus H\oplus
\spann\{w\}$. We set
$$
S_n(x+y+sw)=S_{n-1}x+Ry+sy_n\ \ \text{for}\ \ x\in Z_{n-1},\ y\in H\
\ \text{and}\ \ s\in\K.
$$
The operator $S_n$ is continuous since $S_{n-1}$ and $R$ are
continuous. Clearly (a1) and (a5) for $k=n$ are satisfied. Next,
$S_n(Z_n)\supseteq S_{n-1}(Z_{n-1})+R(H)$. Since $S_{n-1}(Z_{n-1})$
is dense in $Z_{n}$ (condition (a2) for $k=n-1$) and $R(H)$ is dense
in $X_{n+1}$, we see that $S_n(Z_n)$ is dense in $Z_{n+1}=Z_n\oplus
X_n$,which gives us (a2) for $k=n$. From the last display and the
relation $Qy_n\neq 0$ it follows that (a3) is satisfied for $k=n$.
Finally, since $S_nw=y_n$ from the definition of $w$  we see that
(a4) for $k=n$ is also satisfied. Thus the inductive construction of
$S_k$ and $y_k$ is complete.

Condition (a2) ensures that there is a unique operator $S\in L(X)$
such that $S\bigr|_{Z_n}=S_n$ for any $n\in\Z_+$. From (a4) it now
follows that
$$
S^{k+1}y_{k-1}=y_k\ \ \ \text{for each $k\in\Z_+$}.
$$
According to the above display, the set $A=\{y_n:n\in\Z_+\}$ is
contained in the orbit $O(S,y_0)$. By (a5) $A$ is dense in $Y$.
Since $Y$ is dense in $X_0$ and carries a stronger topology, $A$ is
dense in $X_0=Z_0$. By (a2) $S^m(A)$ is dense in $Z_m$ for each
$m\in\Z_+$. Since $A\subset O(S,y_0)$, we have that $S^m(A)\subset
O(S,y_0)$ and therefore $O(S,y_0)\cap Z_m$ is dense in $Z_m$ for
each $m\in\Z_+$. Hence $O(S,y_0)$ is dense in $X$. That is, $y_0$ is
a hypercyclic vector for $S$.
\end{proof}

\begin{remark} The orbit of the hypercyclic vector constructed in
the proof of Lemma~\ref{sumlcs} is not just dense. It is
sequentially dense. The latter property is strictly stronger than
density already for countable direct sums of separable infinite
dimensional Banach spaces.
\end{remark}

\subsection{Proof of Theorem~\ref{sumM}}

Let $X_n\in \M$ for each $n\in\Z_+$ and
$X=\bigoplus\limits_{n=0}^\infty X_n$. We shall apply
Lemma~\ref{sumlcs}. Condition (\ref{sumlcs}.1) is satisfied
according Lemma~\ref{l11}. Indeed $X_1$ admits a Banach $s$-disk
with dense span. By Lemma~\ref{clM}, spaces $X_0\times X_1$ and
$X_n\times \K$ belong to $\M$. From Lemma~\ref{dera} it follows that
conditions (\ref{sumlcs}.2) and (\ref{sumlcs}.3) are also satisfied.
Thus by Lemma~\ref{sumlcs}, there is hypercyclic $T\in L(X)$.

\subsection{Proof of Theorem~\ref{sumfre}}

As we have already mentioned, s Fr\'echet space $X$ belongs to $\M$
if and only if $X$ is infinite dimensional, separable and
non-isomorphic to $\omega$. Moreover, in any separable Fr\'echet
space there is an $\ell_1$-sequence with dense span.

\begin{lemma}\label{aux2} Let $X$ and $Y$ be separable infinite
dimensional Fr\'echet spaces. Then the following conditions are
equivalent
\begin{itemize}\itemsep=-2pt
\item[\rm (\ref{aux2}.1)]there is no $T\in L(X,Y)$ with dense
range$;$
\item[\rm (\ref{aux2}.2)]$X$ is isomorphic to $\omega$ and $Y$ is not
isomorphic to $\omega$.
\end{itemize}
\end{lemma}

\begin{proof} If both $X$ and $Y$ are isomorphic to $\omega$, then
obviously there is a surjective $T\in L(X,Y)$. If $X$ is isomorphic
to $\omega$, $Y$ is not and $T\in L(X,Y)$, then $Z=T(X)$ carries
minimal locally convex topology \cite{bonet} since $\omega$ does. It
follows that $Z$ is either finite dimensional or isomorphic to
$\omega$ and therefore complete. Hence $Z$ is closed in $Y$. It
follows that $Z=\overline{Z}\neq Y$ since $Y$ is neither finite
dimensional nor isomorphic to $\omega$. Thus there is not $T\in
L(X,Y)$ with dense range. It remains to show that there is $T\in
L(X,Y)$ with dense range if $X$ is not isomorphic to $\omega$. In
this case the topology of $X$ is not weak and it remains to apply
Lemma~\ref{dera}.
\end{proof}

\begin{lemma}\label{cases} Let $X$ be the countable locally convex
direct sum of a sequence of separable Fr\'echet spaces infinitely
many of which are infinite dimensional. Then one of the following
two possibilities occurs:
\begin{itemize}\itemsep=-2pt
\item[\rm (\ref{cases}.1)]$X$ is isomorphic to $Y\oplus Z$, where $Y$ is a
separable infinite dimensional Fr\'echet space and $Z$ is the
locally convex direct sum of an infinite countable number of copies
of $\omega;$
\item[\rm (\ref{cases}.2)]$X$ is isomorphic to $\bigoplus\limits_{n=0}^\infty Y_n$,
where each $Y_n$ is a separable infinite dimensional Fr\'echet space
non-isomorphic to $\omega$.
\end{itemize}
\end{lemma}

\begin{proof} Separating the finite dimensional spaces, spaces
isomorphic to $\omega$ and infinite dimensional spaces
non-isomorphic to $\omega$, we see that
$$
X=\bigoplus_{\alpha\in A}X_\alpha\oplus\bigoplus_{\beta\in
B}X_\beta\oplus\bigoplus_{\gamma\in C}X_\gamma,
$$
where the sets $A$, $B$ and $C$ are pairwise disjoint, $X_\alpha$ is
isomorphic to $\omega$ for each $\alpha\in A$, $X_\beta$ is a
separable infinite dimensional Fr\'echet space non-isomorphic to
$\omega$ for any $\beta\in B$, $X_\gamma$ is finite dimensional for
each $\gamma\in C$, $A\cup B$ is infinite and countable and $C$ is
either finite or countable.

If $B$ and $C$ are finite, then $A$ is infinite. Pick $\alpha_0\in
A$. Then
$$
X=Y\oplus\bigoplus_{\alpha\in A\setminus\{\alpha_0\}} X_\alpha,\ \
\text{where}\ \  Y=X_{\alpha_0}\oplus\bigoplus_{\beta\in
B}X_\beta\oplus\bigoplus_{\gamma\in C}X_\gamma.
$$
Clearly $Y$ is a separable infinite dimensional Fr\'echet space.
Since each $X_\alpha$ is isomorphic to $\omega$, we fall into the
case (\ref{cases}.1). If $B$ is finite and $C$ is infinite, then $A$
is infinite and both $A$ and $C$ can be enumerated by elements of
$\Z_+$: $A=\{\alpha_n:n\in\Z_+\}$, $C=\{\gamma_n:n\in\Z_+\}$. Then
$$
X=Y\oplus\bigoplus_{n\in\Z_+} (X_{\alpha_{n+1}}\oplus
X_{\gamma_n}),\ \ \text{where}\ \
Y=X_{\alpha_0}\oplus\bigoplus_{\beta\in B}X_\beta.
$$
Again $Y$ is a separable infinite dimensional Fr\'echet space. Since
each $X_{\alpha_{n+1}}\oplus X_{\gamma_n}$ is isomorphic to
$\omega$, (\ref{cases}.1) is satisfied. If $B$ is infinite and
$A\cup C$ is infinite, then we enumerate both $B$ and $A\cup C$ by
the elements of $\Z_+$: $B=\{\beta_n:n\in\Z_+\}$, $A\cup
C=\{\rho_n:n\in\Z_+\}$. We arrive to
$$
X=\bigoplus_{n\in\Z_+} (X_{\beta_n}\oplus X_{\rho_n}).
$$
Since each $X_{\beta_n}\oplus X_{\rho_n}$ is a separable infinite
dimensional Fr\'echet space non-isomorphic to $\omega$,
(\ref{cases}.2) is satisfied. Finally, if $B$ is infinite and $A\cup
C$ is finite, we fix $\beta_0\in B$ and write
$$
X=Z\oplus \bigoplus_{\beta\in B\setminus\{\beta_0\}} X_\beta,\ \
\text{where}\ \ Z=X_{\beta_0}\oplus\bigoplus_{\rho\in A\cup
C}X_\rho.
$$
Again $Z$ and each $X_\beta$ are separable infinite dimensional
Fr\'echet spaces non-isomorphic to $\omega$ and we fall into the
case (\ref{cases}.2).
\end{proof}

We are ready to prove Theorem~\ref{sumfre}. Let $X$ be a countable
infinite direct sum of separable Fr\'echet spaces. If all the spaces
in the sum, except for finitely many, are finite dimensional, then
$X$ is isomorphic to $Y\times \phi$, where $Y$ is a Fr\'echet space.
According to Theorem~\ref{timesphi}, $X$ admits no cyclic operator
with dense range. In particular, there is no supercyclic operator on
$X$. If there are infinitely many infinite dimensional spaces in the
sum defining $X$, then according to Lemma~\ref{cases}, we see that
$X$ is isomorphic to
$$
Y=\bigoplus_{n\in\Z_+}Y_n,
$$
where $Y_n$ are all separable infinite dimensional Fr\'echet spaces
and either all $Y_n$ are non-isomorphic to $\omega$ or all $Y_n$ for
$n\geq 1$ are isomorphic to $\omega$. In any case from
Lemma~\ref{aux2} it follows that there exists $T_0\in
L(Y_0,Y_0\oplus Y_1)$ with dense range and for each $n\in\N$, there
exists $T_n\in L(Y_n,Y_{n+1}\times\K)$ with dense range. By
Lemma~\ref{sumlcs}, there is a hypercyclic operator on $X$. The
proof of Theorem~\ref{sumfre} is complete.

\section{Hypercyclic operators on countable unions of spaces}

The following lemma is a main tool in the proof of
Theorem~\ref{lbs}.

\begin{lemma}\label{limit} Let a locally convex space be the union
of an increasing sequence $\{X_n\}_{n\in\N}$ of its closed linear
subspaces. Assume also that for any $n\in\N$ there is an
$\ell_1$-sequence with dense span in $X_n$ and the topology of
$X_n/X_{n-1}$ is not weak, where $X_0=\{0\}$. Then there exists a
linear map $S:X\to X$ and $x_0\in X_1$ such that for any $n\in\N$,
$S\bigr|_{X_n}\in L(X_n,X_{n+1})$ and the orbit
$O(S,x)=\{S^kx:k\in\Z_+\}$ is dense in $X$.
\end{lemma}

Note that we do  not claim continuity of the above operator $S$ on
$X$. Although if, for instance, $X$ is the inductive limit of the
sequence $\{X_n\}$, then continuity of $S$ will immediately follow
from the continuity of the restrictions $S\bigr|_{X_n}$.

\begin{proof}[Proof of $Lemma~\ref{limit}$] For each $n\in\N$, let
$\{x_{n,k}\}_{k\in\Z_+}$ be an $\ell_1$-sequence with dense span in
$X_n$. For any $n\in\N$, we apply Lemma~\ref{equi1} with the triple
of spaces $(Y,Y_1,Y_0)$ being $(X,X_n,X_{n-1})$ to obtain a sequence
$\{f_{n,k}\}_{k\in\Z_+}$ in $X'$ such that $\{f_{n,k}:k\in\Z_+\}$ is
uniformly equicontinuous, each $f_{n,k}$ vanishes on $X_{n-1}$ and
$\phi\subseteq \bigl\{\{f_{n,k}(x)\}_{k\in\Z_+}:x\in X_{n}\bigr\}$.
According to Lemma~\ref{l11}, there exists a Banach disk $K$ in $X$
such that $X_K$ is a dense subspace of $X_1$ and the Banach space
$X_K$ is separable. Let $\{U_n\}_{n\in\N}$ be a base of topology of
$X_K$. We shall construct inductively a sequence of operators
$S_k\in L(X,X_{k+1})$ and vectors $y_k\in X_K$ satisfying the
following conditions for any $k\in\Z_+$:
\begin{itemize}\itemsep=-2pt
\item[(p1)]$S_jx=S_kx$ for any $j<k$ and $x\in X_j$;
\item[(p2)]$S_k(X_{k})$ is dense in $X_{k+1}$;
\item[(p3)]$f_{k+1,0}(S_k\dots S_1y_k)\neq0$;
\item[(p4)]$S_k\dots S_1y_{k-1}=y_k$ if $k\geq 2$;
\item[(p5)]$y_k\in U_k$.
\end{itemize}
Consider the linear map $S_1:X \to X_2$ defined by the formula
$$
S_1 x=\sum_{k=0}^\infty 2^{-k}f_{1,k}(x) x_{2,k}.
$$
Since $\{x_{2,k}\}_{k\in\Z_+}$ is an $\ell_1$-sequence in $X_2$ and
$\{f_{1,k}:k\in\Z_+\}$ is uniformly equicontinuous, the above
display defines a continuous linear operator from $X$ to $X_2$.
Since $\phi\subseteq \bigl\{\{f_{1,k}(x)\}_{k\in\Z_+}:x\in
X_1\bigr\}$, $S_1(X_1)$ contains $\spann\{x_{2,k}:k\in\Z_+\}$. Hence
$S_1(X_1)$ is dense in $X_2$. Since $X_K$ is dense in $X_1$, $S_1$
has dense range and $X_2\cap \ker f_{2,0}$ is nowhere dense in
$X_2$, we can pick $y_1\in U_1$ such that $f_{2,0}(S_1y_1)\neq0$.
The basis of induction has been constructed. Assume now that
$n\geq2$ and $y_k\in X_K$, $S_k\in L(X,X_{k+1})$, satisfying
(p1--p5) for $k\leq n-1$, are already constructed. According to (p3)
for $k=n-1$, we have $f_{n,0}(w)\neq 0$, where $w=S_{n-1}\dots
S_1y_{n-1}$. Since $H=X_n\cap \ker f_{n,0}$ is a closed hyperplane
in $X_{n}$, we have $X_n=H\oplus\spann\{w\}$. Let also $H_0=H\cap
\ker f_{n,1}$. Then $H_0$ is a closed hyperplane of $H$. By (p2) for
$k<n$, the operator $S_{n-1}\dots S_1$ from $X_1$ to $X_n$ has dense
range. Since $X_K$ is dense in $X_1$, we can pick $y_n\in U_n$ such
that $u=S_{n-1}\dots S_1y_n\notin H_0\oplus\spann\{w\}$. Thus
$X_{n}=H_0\oplus \spann\{u,w\}$. From definition of $H_0$ it now
follows that the matrix $\begin{pmatrix}f_{n,0}(w)&f_{n,1}(w)\\
f_{n,0}(u)&f_{n,1}(u)\end{pmatrix}$ is invertible. The latter
property allows us for any two vectors $x_0,x_1\in X_{n+1}$ to find
$y_0,y_1\in\spann\{x_0,x_1\}\subset X_{n+1}$ satisfying
$f_{n,0}(w)y_0+f_{n,1}(w)y_1=x_0$ and
$f_{n,0}(u)y_0+f_{n,1}(u)y_1=x_1$. Pick any vector $v\in X_{n+1}$
such that $f_{n+1,0}(v)\neq 0$ and let
$$
x_0=y_n-S_{n-1}w-\sum_{k=0}^\infty 2^{-k}f_{n,k+2}(w)x_{n+1,k},\quad
x_1=v-S_{n-1}u-\sum_{k=0}^\infty 2^{-k}f_{n,k+2}(u)x_{n+1,k}.
$$
The above series converge since $\{x_{n+1,k}\}_{k\in\Z_+}$ is an
$\ell_1$-sequence and $\{f_{n,k}:k\in\Z_+\}$ is uniformly
equicontinuous. Applying the above property to the pair $x_0,x_1\in
X_{n+1}$, we find $y_0,y_1\in X_{n+1}$ such that
$$
\text{$f_{n,0}(w)y_0+f_{n,1}(w)y_1=x_0$ and
$f_{n,0}(u)y_0+f_{n,1}(u)y_1=x_1$.}
$$
Consider now the linear map $S_n:X\to X_{n+1}$ defined by the
formula
$$
S_nx=S_{n-1}x+f_{n,0}(x)y_0+f_{n,1}(x)y_1+\sum_{k=0}^\infty
2^{-k}f_{n,k+2}(x)x_{n+1,k}.
$$
The above display defines a continuous linear operator since
$\{x_{n+1,k}\}_{k\in\Z_+}$ is an $\ell_1$-sequence and
$\{f_{n,k}:k\in\Z_+\}$ is uniformly equicontinuous. From the last
three displays it follows that $S_nw=y_n$ and $S_nu=v$. From
definition of $w$ and $u$ and the relation $f_{n+1,0}(v)\neq 0$ it
follows that (p3) and (p4) for $k=n$ are satisfied. Clearly (p5) for
$k=n$ is also satisfied. Since each $f_{n,k}$ vanishes on $X_{n-1}$,
we have from the last display that $S_{n}x=S_{n-1}x$ for any $x\in
X_{n-1}$. Hence (p1) for $k=n$ is satisfied. It remains to verify
(p2) for $k=n$. Let $U$ be a non-empty open subset of $X_{n+1}$.
Since $E=\spann\{x_{n+1,k}:k\in\Z_+\}$ is dense in $X_{n+1}$, we can
find $x\in E$ and a convex balanced neighborhood $W$ of zero in
$X_{n+1}$ such that $x+W \subseteq U$. Since $\phi\subseteq
\bigl\{\{f_{n,k}(x)\}_{k\in\Z_+}:x\in X_n\bigr\}$, for each $x\in
E=\spann\{x_{n+1,k}:k\in\Z_+\}$, we can pick $y\in X_n$ such that
$f_{n,0}(y)=f_{n,1}(y)=0$ and $x=\sum\limits_{k=0}^\infty
2^{-k}f_{n,k+2}(y)x_{n+1,k}$. It follows that $S_ny=S_{n-1}y+x$. By
(p2) for $k=n-1$, $S_{n-1}(X_{n-1})$ is dense in $X_n$. Since
$S_{n-1}y\in X_n$, we can find $r\in X_{n-1}$ such that $S_{n-1}r\in
S_{n-1}y-W$. By the already proven property (p1) for $k=n$,
$S_{n-1}r=S_nr$. Hence $S_nr\in S_{n-1}y-W$. Using the equality
$S_ny=S_{n-1}y+x$, we get $S_n(y-r)\in x+W\subseteq U$. Hence any
non-empty open subset of $X_{n+1}$ contains elements of $S_n(X_n)$,
which proves (p2) for $k=n$.  Thus the inductive construction of
$S_k$ and $y_k$ is complete.

Condition (p2) ensures that there is a unique linear map $S:X\to X$
such that $S\bigr|_{X_n}=S_n\bigr|_{X_n}$ for any $n\in\N$. From
(p4) it now follows that $S^{k+1}y_{k}=y_{k+1}$ for each $k\in\N$.
Thus the set $A=\{y_n:n\in\N\}$ is contained in the orbit
$O(S,y_1)$. By (p5) $A$ is dense in $X_K$ and therefore is dense in
$X_1$. By (p2) $S^m(A)$ is dense in $X_{m+1}$ for each $m\in\Z_+$.
Since $A\subset O(S,y_1)$, we have that $S^m(A)\subset O(S,y_1)$ and
therefore $O(S,y_1)\cap X_m$ is dense in $X_m$ for each $m\in\N$.
Hence $O(S,y_1)$ is dense in $X$.
\end{proof}

Before proving Theorem~\ref{lbs}, we need to make the following two
elementary observations.

\begin{lemma}\label{lbs0}Let $X$ be an LB-space and $Y$ be a closed
linear subspace of $X$. Then either $X/Y$ is finite dimensional or
the topology of $X/Y$ is not weak.
\end{lemma}

\begin{proof} Since $X$ is an LB-space, it is the inductive limit of
a  sequence $(X_n,\|\cdot \|_n)$ of Banach spaces. If $X/Y$ is
infinite dimensional, we can find a linearly independent sequence
$\{f_n\}_{n\in\Z_+}$ in $X'$ such that each $f_n$ vanishes on $Y$.
Next, we find a sequence $\{\epsilon_n\}_{n\in\Z_+}$ of positive
numbers converging to zero fast enough to ensure that
$\epsilon_n\bigl\|f_n\bigr|_{X_k}\bigr\|_k\to 0$ as $n\to\infty$ for
each $k\in\N$. It follows that the sequence $\epsilon_nf_n$ is
pointwise convergent to zero on $X$. Since any LB-space is barrelled
\cite{shifer,bonet}, the set $\{\epsilon_n f_n:n\in\Z_+\}$ is
uniformly equicontinuous. Hence
$p(x)=\sup\{\epsilon_n|f_n(x)|:n\in\Z_+\}$ is a continuous seminorm
on $X$. Since each $f_n$ vanishes on $Y$, $Y\subseteq \ker p$. Then
$\widetilde p(x+Y)=p(x)$ is a continuous seminorm on $X/Y$. Since
$f_n$ are linearly independent $\ker p$ has infinite codimension in
$X$ and therefore $\ker\widetilde p$ has infinite codimension in
$X/Y$. Hence the topology of $X/Y$ is not weak.
\end{proof}

\begin{lemma}\label{lbs1} Let $X$ be an inductive limit of a
sequence $\{X_n\}_{n\in\Z_+}$ of Banach spaces such that $X_0$ is
dense in $X$. Then $X$ has no quotients isomorphic to $\phi$.
\end{lemma}

\begin{proof} Assume that $X$ has a quotient isomorphic to $\phi$.
By Lemma~\ref{phi2} then $X$ is isomorphic to $Y\times \phi$ for
some closed linear subspace $Y$ of $X$. Let $J:X_0\to X$ be the
natural embedding. Since $X_0$ is dense in $X$, $J$ has dense range.
Hence $J':X'\to X'_0$ is injective. Since $X$ is isomorphic to
$Y\times\phi$, we have that $X'_\beta$ is isomorphic to
$Y'_\beta\times \omega$ ($\omega$ is naturally isomorphic to
$\phi'_\beta$). Hence, there exists an injective continuous linear
operator from $\omega$ to the Banach space $X'_0$ (with the topology
$\beta(X'_0,X_0)$). That is impossible, since any injective
continuous linear operator from $\omega$ to a locally convex space
is an isomorphism onto image and $\omega$ is non-normable.
\end{proof}

\subsection{Proof of Theorem~\ref{lbs}}

Throughout this section $X$ is the inductive limit of a sequence
$\{X_n\}_{n\in\Z_+}$ of separable Banach spaces. Let also
$\overline{X}_n$ be the closure of $X_n$ in $X$. If one of
$\overline{X}_n$ coincides with $X$, then it is easy to see that
$X\in\M$ and therefore there is a hypercyclic operator on $X$ by
Corollary~\ref{saan1}. It remains to consider the case when
$\overline{X}_n\neq X$ for each $n\in\Z_+$. In this case $X$ is the
union of the increasing sequence $\{\overline{X}_n\}$ of proper
closed linear subspaces. Passing to a subsequence, if necessary, we
can assume that one of the two following conditions is satisfied:
either $\dim \overline{X}_{n+1}/\overline{X}_{n}=\infty$ for each
$n\in\Z_+$ or $0<\dim \overline{X}_{n+1}/\overline{X}_{n}<\infty$
for each $n\in\Z_+$.

{\bf Case 1} \ $0<\dim \overline{X}_{n+1}/\overline{X}_{n}<\infty$
for each $n\in\Z_+$. In this case, for any $n\in\Z_+$, we can pick a
non-trivial finite dimensional subspace $Y_n$ of $X_{n+1}$ such that
$\overline{X_n}\oplus Y_n=\overline{X_{n+1}}$. Thus the vector space
$X$ can be written as an algebraic direct sum
\begin{equation}\label{di1}
X=\overline{X_0}\oplus\bigoplus_{n=0}^\infty Y_n.
\end{equation}
Apart from the original topology $\tau$ on $X$, we can consider the
topology $\theta$, turning the sum (\ref{di1}) into a locally convex
direct sum. Obviously $\tau\subseteq\theta$. On the other hand, if
$W$ is a balanced convex $\theta$-neighborhood of $0$ in $X$, then
$W\cap \overline{X_n}$ is a $\tau$-neighborhood of zero in
$\overline{X_n}$ for any $n\in\Z_+$. Indeed, it follows from the
fact that $\overline{X_n}=\overline{X_0}\oplus Z_n$, where
$Z_n=\bigoplus\limits_{j=0}^{n-1} Y_j$ and $Z_n$ is finite
dimensional. Since the topology of each $X_n$ is stronger than the
one inherited from $X$, we see that $W\cap X_n$ is a neighborhood of
zero in $X_n$ for each $n\in\Z_+$. By definition of the inductive
limit topology, $W$ is a $\tau$-neighborhood of zero in $X$. Hence
$\theta\subseteq\tau$. Thus $\theta=\tau$ and therefore $X$ is
isomorphic to $\overline{X_0}\times Y$, where $Y$ is the locally
convex direct sum of $Y_n$ for $n\in\Z_+$. Since $Y_n$ are finite
dimensional, $Y$ is isomorphic to $\phi$. Since $\overline{X_0}$ is
the inductive limit of the sequence $\overline{X_0}\cap X_n$ of
Banach spaces, from Lemma~\ref{lbs1} it follows that
$\overline{X_0}$ has no quotients isomorphic to $\phi$. Since $X$ is
isomorphic to $\overline{X_0}\times\phi$, Theorem~\ref{timesphi}
implies that there are no cyclic operators with dense range on $X$.

{\bf Case 2} \ $\overline{X}_{n+1}/\overline{X}_{n}$ is infinite
dimensional for each $n\in\Z_+$. By Lemma~\ref{lbs0}, the topology
of each $\overline{X}_{n+1}/\overline{X}_{n}$ is not weak. Let
$n\in\Z_+$. Since $X_n$ is a separable Banach space, there is an
$\ell_1$-sequence $\{x_{n,k}\}_{k\in\Z_+}$ in $X_n$ with dense span.
Since the topology on $X_n$ inherited from $X$ is weaker than the
Banach space topology of $X_n$, $\{x_{n,k}\}_{k\in\Z_+}$ is an
$\ell_1$-sequence with dense span in $\overline{X}_{n}$. By
Lemma~\ref{limit}, there exists a linear map $S:X\to X$ and $x_0\in
X$ such that for any $n\in\N$, $S\bigr|_{\overline{X}_n}\in
L(\overline{X}_n,\overline{X}_{n+1})$ and the orbit
$O(S,x)=\{S^kx:k\in\Z_+\}$ is dense in $X$. Since the topology of
$X_n$ is stronger than the one inherited from $X$, we have that each
restriction $S\bigr|_{X_n}$ is a continuous linear operator from
$X_n$ to $X$. Since $X$ is the inductive limit of the sequence
$\{X_n\}_{n\in\Z_+}$, $S:X\to X$ is continuous. Hence $S$ is a
hypercyclic continuous linear operator on $X$. The proof of
Theorem~\ref{lbs} is now complete.

\section{Remarks on mixing versus hereditarily hypercyclic}

We start with the following remark. As we have already mentioned,
$\phi$ supports no supercyclic operator \cite{bonper}, which follows
also from Theorem~\ref{timesphi}. On the other hand, $\phi$ supports
a transitive operator \cite{fre}. The latter statement can be easily
strengthened with the help of Corollary~\ref{main1}. Namely, take
the backward shift $T$ on $\phi$. That is $Te_0=0$ and
$Te_n=e_{n-1}$ for $n\geq 1$, where $\{e_n\}_{n\in\Z_+}$ is the
standard basis in $\phi$. Clearly $T$ is a generalized backward
shift and therefore $T$ is an extended backward shift. By
Corollary~\ref{main1}, $I+T$ is mixing. Thus we have the following
proposition.

\begin{proposition}\label{uuu} $\phi$ supports a mixing operator and supports
no supercyclic operators.
\end{proposition}

On the other hand, a topological vector space of countable algebraic
dimension can support a hypercyclic operator, as observed by several
authors, see \cite{fre}, for instance. The following proposition
formalizes and extends this observation.

\begin{proposition}\label{hmi}Let $X$ be a normed space of countable
algebraic dimension. Then there exists a hypercyclic mixing operator
$T\in L(X)$.
\end{proposition}

\begin{proof} Let $\overline{X}$ be the completion of $X$. Then
$\overline{X}$ is a separable infinite dimensional Banach space. By
Corollary~\ref{saan1}, there is a hereditarily hypercyclic operator
on $S\in L(\overline{X})$. Let $x\in \overline{X}$ be a hypercyclic
vector for $S$ and $E$ be the linear span of the orbit of $x$:
$E=\spann\{S^nx:n\in\Z_+\}$. Grivaux \cite{gri2} demonstrated that
for any two countably dimensional dense linear subspaces $E$ and $F$
of a separable infinite dimensional Banach space $Y$, there is an
isomorphism $J:Y\to Y$ such that $J(E)=F$. Hence there is an
isomorphism $J: \overline{X}\to  \overline{X}$ such that $J(X)=E$.
Let now $T_0=J^{-1}SJ$. Since $J(X)=E$ and $E$ is $S$-invariant, $X$
is $T_0$-invariant. Thus the restriction $T$ of $T_0$ to $X$ is a
continuous linear operator on $X$. Moreover, since the $S$-orbit of
$x$ is dense in $\overline{X}$, the $T_0$-orbit of $J^{-1}x$ is
dense in $\overline{X}$. Since $J^{-1}x\in X$, the latter orbit is
exactly the $T$-orbit of $J^{-1}x$ and therefore $J^{-1}x$ is
hypercyclic for $T$. Hence $T$ is hypercyclic. Next, $T_0$ is mixing
since it is similar to the mixing operator $S$. Hence $T$ is mixing
as a restriction of a mixing operator to a dense subspace.
\end{proof}

By Proposition~\ref{untr2}, if $X$ is a Baire separable and
metrizable topological vector space, then any mixing $T\in L(X)$ is
hereditarily hypercyclic. From the above proposition it follows that
there are mixing operators on countably dimensional normed spaces.
The next theorem however implies that there are no hereditarily
hypercyclic operators on countably dimensional topological vector
spaces, emphasizing the necessity of the Baire condition in
Proposition~\ref{untr2}.

\begin{theorem}\label{al1} Let $X$ be a topological vector space
such that there exists a hereditarily universal family
$\{T_n:n\in\Z_+\}\subset L(X)$. Then $\dim X>\aleph_0$.
\end{theorem}

\begin{proof} Since the topology of any topological vector space can
be defined by a family of quasinorms \cite{shifer}, we can pick a
non-zero continuous quasinorm $p$ on $X$. That is,
$p:X\to[0,\infty)$ is non-zero, continuous, $p(x+y)\leq p(x)+p(y)$
for any $x,y\in X$, $p(0)=0$, $p(zx)=p(x)$ if $x\in X$, $z\in\K$,
$|z|=1$ and $(X,\tau_p)$ is a (not necessarily Hausdorff)
topological vector space, where $\tau_p$ is the topology defined by
the pseudometric $d(x,y)=p(x-y)$. The latter property implies that
$p(tx_n)\to 0$ for any $t\in\K$ and any sequence
$\{x_n\}_{n\in\Z_+}$ in $X$ such that $p(x_n)\to 0$.

Let $\kappa$ be the first uncountable ordinal (commonly denoted
$\omega_1$). We shall construct inductively sequences
$\{x_\alpha\}_{\alpha<\kappa}$ and $\{A_\alpha\}_{\alpha<\kappa}$ of
vectors in $X$ and subsets of $\Z_+$ respectively such that for any
$\alpha<\kappa$,
\begin{itemize}\itemsep=-2pt
\item[(s1)]$A_\alpha$ is infinite and $x_\alpha$ is a universal
vector for the family $\{T_n:n\in A_\alpha\}$;
\item[(s2)]$p(T_nx_\beta)\to 0$ as $n\to\infty$, $n\in A_\alpha$ for any
$\beta<\alpha$;
\item[(s3)]$A_\alpha\setminus A_\beta$ is finite for any
$\beta<\alpha$.
\end{itemize}

For the basis of induction we take $A_0=\Z_+$ and $x_0$ being a
universal vector for the family $\{T_n:n\in\Z_+\}$. It remains to
describe the induction step. Assume that $\gamma<\kappa$ and
$x_\alpha$, $A_\alpha$ satisfying (s1--s3) for $\alpha<\gamma$ are
already constructed. We have to construct $x_\gamma$ and $A_\gamma$
satisfying (s1--s3) for $\alpha=\gamma$.

{\bf Case 1:} $\gamma$ has the immediate predecessor. That is
$\gamma=\rho+1$ for some ordinal $\rho<\kappa$. Since $x_\rho$ is
universal for $T_n:n\in A_\rho$, we can pick an infinite subset
$A_\gamma\subset A_\rho$ such that $p(T_nx_\rho)\to 0$ as
$n\to\infty$, $n\in A_\gamma$. Since $A_\gamma$ is contained in
$A_\rho$, from (s3) for $\alpha\leq\rho$ it follows that
$A_\gamma\setminus A_\beta$ is finite for any $\beta<\gamma$. Hence
(s3) for $\alpha=\gamma$ is satisfied. Now from (s3) for
$\alpha=\gamma$ and (s2) for $\alpha<\gamma$ it follows that (s2) is
satisfied for $\alpha=\gamma$. Next, since $\{T_n:n\in\Z_+\}$ is
hereditarily universal, we can pick $x_\gamma\in X$ universal for
$\{T_n:n\in A_\gamma\}$. Hence (s1) for $\alpha=\gamma$ is also
satisfied.

{\bf Case 2:} $\gamma$ is a limit ordinal. Since $\gamma$ is a
countable ordinal, we can pick a strictly increasing sequence
$\{\alpha_n\}_{n\in\Z_+}$ of ordinals such that
$\gamma=\sup\{\alpha_n:n\in\Z_+\}$. Now pick consecutively $n_0$
from $A_{\alpha_0}$, $n_1>n_0$ from $A_{\alpha_0}\cap A_{\alpha_1}$,
$n_2>n_1$ from $A_{\alpha_0}\cap A_{\alpha_1}\cap A_{\alpha_2}$ etc.
The choice is possible since by (s3) for $\alpha<\gamma$, each
$A_{\alpha_0}\cap {\dots} \cap A_{\alpha_n}$ is infinite. Now let
$A_\gamma=\{n_j:j\in\Z_+\}$. Since $A_\gamma\setminus
A_{\alpha_j}\subseteq\{n_0,\dots,n_{j-1}\}$, $A_\gamma\setminus
A_{\alpha_j}$ is finite for each $j\in\Z_+$. Now if $\beta<\gamma$,
we can pick $j\in\Z_+$ such that $\beta<\alpha_j<\gamma$. Then
$A_\gamma\setminus A_\beta\subseteq (A_\gamma\setminus
A_{\alpha_j})\cup (A_{\alpha_j}\setminus A_{\beta})$ is finite by
(s3) with $\alpha=\alpha_j$. Moreover, since $A_\gamma$ is contained
in $A_{\alpha_j}$ up to a finite set, from (s2) with
$\alpha=\alpha_j$ it follows that $p(T_nx_\beta)\to 0$ as
$n\to\infty$, $n\in A_\gamma$. Hence (s2) and (s3) for
$\alpha=\gamma$ are satisfied. Finally, since $\{T_n:n\in\Z_+\}$ is
hereditarily universal, we can pick $x_\gamma\in X$ universal for
$\{T_n:n\in A_\gamma\}$. Hence (s1) for $\alpha=\gamma$ is also
satisfied. This concludes the construction of
$\{x_\alpha\}_{\alpha<\kappa}$ and $\{A_\alpha\}_{\alpha<\kappa}$
satisfying (s1--s3).

In order to prove that $\dim X>\aleph_0$, it suffices to show that
vectors $\{x_\alpha\}_{\alpha<\kappa}$ are linearly independent.
Assume the contrary. Then there are $n\in\N$,
$z_1,\dots,z_n\in\K\setminus\{0\}$ and ordinals
$\alpha_1<{\dots}<\alpha_n<\kappa$ such that $\sum\limits_{j=1}^n
z_j x_{\alpha_j}=0$. By (s2) with $\alpha=\alpha_n$, we see that
$p(T_kx_{\alpha_j})\to 0$ as $k\to\infty$, $k\in A_{\alpha_n}$ for
$1\leq j<n$. Denoting $c_j=-z_j/z_n$ for $1\leq j<n$ and using
linearity of $T_k$, we obtain $T_kx_{\alpha_n}=\sum\limits_{1\leq
j<n} c_j T_kx_{\alpha_j}$ for any $k\in\Z_+$. Since $p$ is a
quasinorm, we have
$$
p(T_k x_{\alpha_n})\leq \sum_{1\leq j<n}p(c_j T_kx_{\alpha_j})\to 0\
\ \text{as $k\to\infty$, $k\in A_{\alpha_n}$}.
$$
The above display contradicts universality of $x_{\alpha_n}$ for
$\{T_k:k\in A_{\alpha_n}\}$, which is (s1) with $\alpha=\alpha_n$.
This contradiction completes the proof.
\end{proof}

\begin{corollary}\label{al2}A topological vector space of countable
algebraic dimension supports no hereditarily hypercyclic operators.
\end{corollary}

It is worth noting that there are infinite dimensional separable
normed spaces, which support no supercyclic or transitive operators.
We call a continuous linear operator $T$ on a topological vector
space $X$ {\it simple} if $T$ has shape $T=zI+S$, where $z\in\K$ and
$S$ has finite rank. Observe that a simple operator on an infinite
dimensional topological vector space is never transitive or
supercyclic. Indeed, let $T$ be a simple operator on an infinite
dimensional topological vector space and $\lambda\in\K$, $S\in L(X)$
be such that $T=\lambda I+S$ and $S$ has finite rank. Then $L=S(X)$
is finite dimensional. Since $X/L$ is infinite dimensional, we can
pick non-empty open subsets $U_0$ and $V_0$ of $X/L$ such that
$U_1=\{zu:z\in\K\setminus\{0\},\ u\in U_0\}$ does not intersect
$V_1=\{zv:z\in\K\setminus\{0\},\ v\in V_0\}$. Let $U=\{zx:x+L\in
U_0,\ z\in\K\setminus\{0\}\}$ and $V=\{zx:x+L\in V_0,\
z\in\K\setminus\{0\}\}$. Clearly $U$ and $V$ are non-empty open
subsets of $X$. Using the equalities $T=\lambda I+S$ and $S(X)=L$,
it is easy to see that $T^n(U)\cap V=\varnothing$ for any
$n\in\Z_+$. Hence $T$ is non-transitive. Moreover since $U$ and $V$
are stable under multiplication by non-zero scalars, the projective
orbit $\{zT^nx:n\in\Z_+,\ z\in\K\}$ of any $x\in U$ does not meet
$V$. Hence $U$ contains no supercylic vectors for $T$. Since the set
of supercyclic vectors of any continuous linear operator is either
dense or empty, $T$ is non-supercyclic.

We say that a topological vector space $X$ is {\it simple} if it is
infinite dimensional and any $T\in L(X)$ is simple. Thus simple
topological vector spaces support no supercyclic or transitive
operators. Various examples of simple separable infinite dimensional
normed spaces can be found in the literature
\cite{vald,pol,vmill,mar,shka,shka1}. Moreover, according to
Valdivia \cite{vald}, in any separable infinite dimensional Fr\'eche
spaced there is a dense simple hyperplane. All the examples of this
type existing in the literature with one exception \cite{shka} are
constructed with the help of the axiom of choice and the spaces
produced are not Borel measurable in their completions. In
\cite{shka} there is a constructive example of a simple separable
infinite dimensional pre-Hilbert space $H$ which is a countable
union of compact sets.

Finally recall that an infinite dimensional topological vector space
$X$ is called {\it rigid} if $L(X)$ consists only of the operators
of the form $\lambda I$ for $\lambda\in\K$. Of course, a rigid space
can not be locally convex. Clearly there are no transitive or cyclic
continuous linear operators on a rigid topological vector space.
Since there exist rigid separable $\F$-spaces \cite{kal}, we see
that there are separable infinite dimensional $\F$-spaces on which
there are no cyclic or transitive operators.

\section{Concluding remarks}

We start by observing that the following question remains open.

\begin{question}
Is there a hereditarily hypercyclic operator on a countable direct
sum of separable infinite dimensional Banach spaces?
\end{question}

The most of the above results rely upon the underlying space being
locally convex or at least having plenty of continuous linear
functionals and for a good reason. As mentioned in the previous
section, there are separable infinite dimensional $\F$-spaces on
which there are no cyclic or transitive operators. On the other
hand, the absence of non-zero continuous linear functionals on a
topological vector space does not guarantee the absence of
hypercyclic operators on it. It is well-known \cite{kal} that the
spaces $L_p[0,1]$ for $0\leq p<1$ are separable $\F$-spaces having
no non-zero continuous linear functionals. Ansari \cite{ansa1}
raises a question whether these spaces support hypercyclic
operators. Theorem~\ref{main} provides an easy answer to this
question. Namely, consider the operator $T\in L(L_p[0,1])$,
$Tf(x)=f(x/2)$. It is straightforward to see that $T$ is onto and
has dense generalized kernel. Thus $I+T$ is hereditarily hypercyclic
according to Corollary~\ref{main2}.

It is obvious that an extended backward shift has dense range and
dense generalized kernel. Unfortunately, the converse is not true in
general. This leads naturally to the following question.

\begin{question}
Let $T$ be a continuous linear operator on a separable Banach space,
which has dense range and dense generalized kernel. Is it true that
$I+T$ is mixing or at least hypercyclic?
\end{question}

From Corollary~\ref{main1} and Corollary~\ref{mi00} it follows that
if $T$ is an extended backward shift on a separable infinite
dimensional Banach space $X$, then both $I+T$  and $e^T$ are
hereditarily hypercyclic. This reminds of the following question
raised by Berm\'udez, Bonilla, Conejero and Peris in \cite{holo}.

\medskip
\noindent {\bf Question~B$^{\bf 2}$CP.} \ \it Let $X$ be a complex
Banach space and $T\in L(X)$ be such that its spectrum $\sigma(T)$
is connected and contains $0$. Does hypercyclicity of $I+T$ imply
hypercyclicity of $e^T$? Does hypercyclicity of $e^T$ imply
hypercyclicity of $I+T$? \rm
\medskip

We show that the answer to both parts of the above question is
negative. Before doing this we would like to raise a similar
question, which remains open.

\begin{question} \label{qqq2} Let $X$ be a Banach
space and $T\in L(X)$ be quasinilpotent. Is hypercyclicity of $I+T$
equivalent to hypercyclicity of $e^T$?
\end{question}

First, we introduce some notation. Let $\D=\{z\in\C:|z|<1\}$,
$\H^2(\D)$ be the Hardy Hilbert space on the unit disk and
$\H^\infty(\D)$ be the space of bounded holomorphic functions
$f:\D\to\C$. For any $\phi\in \H^\infty(\D)$, the multiplication
operator
$$
M_\phi f(z)=\phi(z)f(z)
$$
is a bounded linear operator on $\H^2(\D)$. It is also clear that
$\sigma(M_\phi)=\overline{\phi(\D)}$. If $M_\phi^\star$ is the
adjoint of $M_\phi$, then $\sigma(M_\phi^\star)$ is the reflection
of $\sigma(M_\phi)$ with respect to the real axis. The following
proposition is a direct consequence of a theorem by Godefroy and
Shapiro \cite{gs}[Theorem~4.9].

\begin{proposition}\label{shi} Let $\phi\in \H^\infty(\D)$. Then
$M_\phi^\star$ is hypercyclic if and only if
\begin{equation}\label{123}
\phi(\D)\cap \T\neq\varnothing.
\end{equation}
\end{proposition}

The above Proposition calls for the following comment. A bounded
linear operator $T$ on a separable infinite dimensional Banach space
$X$ is said to satisfy the {\it Kitai Criterion} \cite{kitai} (see
also \cite{kc}) if there exist dense subsets $E$ and $F$ of $X$ and
a map $S:F\to F$ such that $TSy=y$ for any $y\in F$, $T^nx\to 0$ and
$S^ny\to 0$ as $n\to\infty$ for any $x\in E$ and $y\in F$. As it is
shown in \cite{kitai,kc}, any operator satisfying the Kitai
Criterion is hypercyclic. Moreover, any operator, satisfying the
Kitai Criterion is hereditarily hypercyclic and therefore mixing
\cite{gri}. Hypercyclicity in the proof of the above result in
\cite{gs} is demonstrated via application of the Kitai Criterion.
Thus the following slightly stronger statement holds.

\begin{corollary}\label{shi1} Let $\phi\in \H^\infty(\D)$. Then
$M_\phi^\star$ is hereditarily hypercyclic if $(\ref{123})$ is
satisfied and $M_\phi^\star$ is non-hypercyclic if $(\ref{123})$ is
not satisfied.
\end{corollary}

Now we demonstrate that the answer to both parts of Question~B$^2$CP
is negative. Consider the subset $U$ of $\C$ being the interior of
the triangle with vertices $-1$, $i$ and $-i$. In other words
$U=\{a+bi:a,b\in\R,\ a<0,\ b-a<1,\ b+a>-1\}$. Next, let
$V=\{a+bi:a,b\in\R,\ 0<b<1,\ |a|<1-\sqrt{1-b^2}\}$. The boundary of
$V$ consists of the interval $[-1+i,1+i]$ and two circle arcs. It is
clear that $U$ and $V$ are bounded, open, connected and simply
connected. From the definition of the sets $U$ and $V$ it
immediately follows that the open set $1+U=\{1+z:z\in U\}$
intersects the unit circle. On the other hand, since $U$ is
contained in the left half-plane, we see that $e^U=\{e^z:z\in
U\}\subseteq \D$. Similarly, we see that $1+V\cap \D=\varnothing$
and the open set $e^V$ intersects the unit circle. According to the
Riemann Theorem \cite{mark}, there exist holomorphic homeomorphisms
$\phi:\D\to U$ and $\psi:\D\to V$. Obviously $\phi,\psi\in
\H^\infty(\D)$. Since $I+M_\phi^\star=M_{1+\phi}^\star$,
$e^{M_\psi^\star}=M^\star_{e^\psi}$ and both $(1+\phi)(\D)=1+U$ and
$e^\psi(\D)=e^V$ intersect the unit circle, Corollary~\ref{shi1}
implies that $I+M_\phi^\star$ and $e^{M_\psi^\star}$ are
hereditarily hypercyclic. Since $I+M_\psi^\star=M_{1+\psi}^\star$,
$e^{M_\phi^\star}=M^\star_{e^\phi}$, $e^\phi(\D)=e^U$ is contained
in $\D$, and $(1+\psi)(\D)=1+V$ does not intersect $\overline{\D}$,
Corollary~\ref{shi1} implies that $e^{M_\phi^\star}$ and
$I+M_\psi^\star$ are non-hypercyclic. Finally, observe that
$\sigma(M_\phi^\star)$ is the closure of $U$ and
$\sigma(M_\psi^\star)$ is the closure of $-V$ and therefore the
spectra of $M_\phi^\star$ and $M_\psi^\star$ are connected and
contain $0$. Taking into account that all separable infinite
dimensional Hilbert spaces are isomorphic, we arrive to the
following result, which answers negatively the Question~B$^2$CP.

\begin{proposition}\label{que} There exist bounded linear operators
$A$, $B$ on the complex Hilbert space $\ell_2$ such that $\sigma(A)$
and $\sigma(B)$ are connected and contain $0$, $I+A$ and $e^B$ are
hereditarily hypercyclic, while $e^A$ and $I+B$ are non-hypercyclic.
\end{proposition}

Finally, if the answer to Question~\ref{qqq2} is affirmative, then
the following interesting question naturally arises.

\begin{question}\label{q3} Let $A$ be a quasinilpotent bounded
linear operator on a complex Banach space $X$ and $f$ be an entire
function on one variable such that $f(0)=f'(0)=1$. Is it true that
hypercyclicity of $f(A)$ is equivalent to hypercyclicity of $I+A$?
\end{question}

\subsection{Spaces $C_k(M)$ and their duals}

Let $(M,d)$ be a separable metric space and $C_k(M)$ be the space of
continuous functions $f:M\to\K$ with the compact-open topology (=the
topology of uniform convergence on compact subsets of $M$). It is
easy to see that $C_k(M)$ is complete. Moreover, $C_k(M)$ is
metrizable if and only $M$ is locally compact and $C_k(M)$ carries
weak topology if and only if $M$ is discrete. On the other hand,
there always is an $\ell_1$-sequence with dense span in $C_k(M)$.
Indeed, let $\overline{M}$ be a metrizable compactification of $M$.
Since $C(\overline{M})$ is a separable Banach space, there is an
$\ell_1$-sequence $\{f_n\}_{n\in\Z_+}$ with dense span in
$C(\overline{M})$. Since $C(\overline{M})$ is  densely and
continuously embedded into $C_k(M)$, $\{f_n\}_{n\in\Z_+}$ is an
$\ell_1$-sequence with dense span in $C_k(M)$. Thus $C_k(M)\in\M$ if
and only if $M$ is non-discrete. Corollary~\ref{saan1} implies now
the following proposition.

\begin{proposition}\label{ck} If $(M,d)$ is a separable non-discrete
metric space, then there is a hereditarily hypercyclic operator on
$C_k(X)$.
\end{proposition}

The reason for inclusion of the above proposition is to demonstrate
that Theorem~\ref{saan} and Corollary~\ref{saan1} are applicable far
beyond metrizable or LB-spaces. The spaces $C_k(M)$ can have quite
ugly structure indeed. For instance, take $M$ being the set $\Q$ of
rational numbers with the metric induced from $\R$, and you have got
the space $C_k(\Q)$, which does not fall into any of the
well-understood and studied classes of locally convex spaces. We
would like to raise the following question.

\begin{question}Characterize separable metric spaces $(M,d)$ such
that $C_k(M)$ supports a dual hypercyclic operator.
\end{question}

It is worth noting that if $M$ is discrete, then either $C_k(M)$ is
finite dimensional or is isomorphic to $\omega$ and therefore does
not support a dual hypercyclic operator (there are no hypercyclic
operators on $\phi=\omega'$). In general, ${\cal M}(M)=(C_k(M))'$
can be naturally identified with the space of finite $\K$-valued
Borel $\sigma$-additive measures on $M$ with compact support. This
dual space is separable in the strong topology if and only if all
compact subsets of $M$ are finite or countable. If it is not the
case, there is no point to look for dual hypercyclic operators on
$C_k(M)$. Thus the only spaces $(M,d)$ for which $C_k(M)$ has a
chance to support a dual hypercyclic operator are non-discrete
spaces with no uncountable compact subsets. The first natural
candidate to consider is $\Q$.

Note also that although Theorem~\ref{dual} provides answers to the
questions of Petersson, mentioned in the introduction, it does not
characterize Fr\'echet spaces, supporting a dual hypercyclic
operator.

\begin{question}Characterize Fr\'echet spaces $X$ such
that $X$ supports a dual hypercyclic operator.
\end{question}

The most natural Fr\'echet space for which we do not know whether it
supports a dual hypercyclic operator is the countable power
$\ell_2^\N$ of the Hilbert space $\ell_2$.

\subsection{The Hypercyclicity Criterion}

The following universality criterion is proved by B\'es and Peris
\cite[Theorem~2.3 and Remark~2.6]{bp}. It is formulated in \cite{bp}
in the case when $X$ is an $\F$-space, but the proof works without
any changes for Baire separable metrizable topological vector
spaces.

\begin{thmbp} \it Let $\{T_n\}_{n\in\Z_+}$ be a sequence of
pairwise commuting continuous linear operators with dense range on a
Baire separable metrizable topological vector space $X$. Then the
following conditions are equivalent:
\begin{itemize}\itemsep=-2pt
\item[\rm (a)]The family $\{T_n\oplus T_n\}_{n\in\Z_+}$ is
universal$;$
\item[\rm (b)]There exists an infinite subset $A$ of $\Z_+$ such that the family
$\{T_n\}_{n\in A}$ is hereditarily universal$;$
\item[\rm (c)]There exist a strictly increasing sequence $\{n_k\}$ of non-negative
integers, dense subsets $E$ and $F$ of $X$ and maps $S_k:F\to X$ for
$k\in\Z_+$ such that $T_{n_k}x\to 0$, $S_ky\to 0$ and
$T_{n_k}S_ky\to y$ as $k\to\infty$ for any $x\in E$ and $y\in F$.
\end{itemize}
\end{thmbp}

We formulate now the so-called Hypercyclicity and Supercyclicity
Criteria, which follow easily from the above theorem.

\begin{thmhc} \it Let $X$ be a Baire separable metrizable
topological vector space and $T\in L(X)$. Then the following
conditions are equivalent:
\begin{itemize}\itemsep=-2pt
\item[\rm (a)]$T\oplus T$ is hypercyclic$;$
\item[\rm (b)]There exists an infinite subset $A$ of $\Z_+$
such that the family $\{T^n\}_{n\in A}$ is hereditarily universal$;$
\item[\rm (c)]There exist a strictly increasing sequence
$\{n_k\}$ of non-negative integers, dense subsets $E$ and $F$ of $X$
and maps $S_k:F\to X$ for $k\in\Z_+$ such that $T^{n_k}x\to 0$,
$S_ky\to 0$ and $T^{n_k}S_ky\to y$ as $k\to\infty$ for any $x\in E$
and $y\in F$.
\end{itemize}
\end{thmhc}

\begin{thmsc} \it Let $X$ be a Baire separable metrizable
topological vector space and $T\in L(X)$. Then the following
conditions are equivalent:
\begin{itemize}\itemsep=-2pt
\item[\rm (a)]$T\oplus T$ is supercyclic$;$
\item[\rm (b)]There exists an infinite subset $A$ of $\Z_+$
and a sequence $\{s_n\}_{n\in A}$ of positive numbers such that the
family $\{s_nT^n\}_{n\in A}$ is hereditarily universal$;$
\item[\rm (c)]There exist a strictly increasing sequence
$\{n_k\}$ of non-negative integers, dense subsets $E$ and $F$ of
$X$, and a sequence $\{s_k\}_{k\in \Z_+}$ of positive numbers and
maps $S_k:F\to X$ for $k\in\Z_+$ such that $s_kT^{n_k}x\to 0$,
$s_k^{-1}S_ky\to 0$ and $T^{n_k}S_ky\to y$ as $k\to\infty$ for any
$x\in E$ and $y\in F$.
\end{itemize}
\end{thmsc}

An operator satisfying the condition (c) of Theorem~HC (respectively
Theorem~SC) is said to satisfy the Hypercyclicity (respectively,
Supercyclicity) Criterion. The long standing question whether any
hypercyclic operator $T$ on a Banach space satisfies the
Hypercyclicity Criterion, was recently solved negatively by Read and
De La Rosa \cite{read}. Their result was extended by Bayart and
Matheron \cite{bama}, who demonstrated that on any separable Banach
space with an unconditional Schauder basis such that the forward
shift operator associated with this basis is bounded, there is a
hypercyclic operator $T$ such that $T\oplus T$ is not hypercyclic.
This leaves open the following question raised in \cite{bama}.

\begin{question} \label{que00}
Does there exist a separable infinite dimensional Banach space $X$
such that any hypercyclic operator on $X$ satisfies the
Hypercyclicity Criterion?
\end{question}

It is observed in \cite{bama} that any $T\in L(\omega)$ satisfies
the Hypercyclicity Criterion. It also follows from
Theorem~\ref{omeg} and Theorem~HC. Thus the above question in the
class of Fr\'echet spaces has an affirmative answer, which leads to
the following problem.

\begin{question}
Characterize separable infinite dimensional Fr\'echet spaces $X$ on
which any hypercyclic operator on $X$ satisfies the Hypercyclicity
Criterion?
\end{question}

It is worth noting that non-hypercyclicity  of $T\oplus T$ in
\cite{read,bama} is ensured by the existence of a non-zero
continuous bilinear form $b:X\times X\to\K$ with respect to which
$T$ is symmetric: $b(Tx,y)=b(x,Ty)$ for any $x,y\in X$. The
following proposition formalizes the corresponding implication.
Similar statements have been proved by many authors in various
particular cases. The proof goes along the same lines as in any of
them.

Let $X$ and $Y$ be topological vector spaces and $b:X\times X\to Y$
be separately continuous and bilinear. We say that $T\in L(X)$ is
{\it $b$-symmetric} if $b(Tx,y)=b(x,Ty)$ for any $x,y\in X$. Recall
also that $b$ is called {\it symmetric} if $b(x,y)=b(y,x)$ for any
$x,y\in X$ and $b$ is called {\it antisymmetric} if $b(x,y)=-b(y,x)$
for any $x,y\in X$.

\begin{proposition} \label{tpt} Let $X$ and $Y$ be a topological vector
spaces, $b:X\times X\to Y$ be separately continuous, non-zero and
bilinear and $T\in L(X)$ be $b$-symmetric. Then $T\oplus T$ is
non-cyclic. If additionally $b$ is non-symmetric, then $T^2$ is
non-cyclic.
\end{proposition}

\begin{proof} Consider the left and right kernels of $b$:
\begin{equation}\label{x0x1}
X_0=\{x\in X:b(x,y)=0\ \text{for any $y\in X$}\}\ \ \text{and}\ \
X_1=\{y\in X:b(x,y)=0\ \text{for any $x\in X$}\}.
\end{equation}
Separate continuity of $b$ implies that $X_0$ and $X_1$ are closed
linear subspaces of $X$. Since $b$ is non-zero, we have $X_0\neq X$
and $X_1\neq X$. From $b$-symmetry of $T$ it follows that $X_0$ and
$X_1$ are both $T$-invariant. Hence $X_0\times X$ and $X\times X_1$
are $T\oplus T$-invariant proper closed subspaces of $X\times X$,
they can not contain a cyclic vector for $T\oplus T$. Assume that
$T\oplus T$ has a cyclic vector $(x,y)\in X\times X$. Then $x\notin
X_0$ and $y\notin X_1$. Consider now a continuous linear operator
$\Phi:X\times X\to Y$ defined by the formula
$$
\Phi(u,v)=b(x,v)-b(u,y).
$$
Since $x\notin X_0$ and $y\notin X_1$, we have $\Phi\neq 0$. On the
other hand, using $b$-symmetry of $T$, we have
$$
\Phi((T\oplus T)^n(x,y))=b(x,T^ny)-b(T^nx,y)=0.
$$
Thus the orbit of $(x,y)$ with respect to $T\oplus T$ lies in the
proper closed linear subspace $\ker\Phi$ of $X$, which contradicts
cyclicity of $(x,y)$ for $T\oplus T$.

Assume now that $b$ is non-symmetric. Then $c(x,y)=b(x,y)-b(y,x)$ is
a non-zero separately continuous bilinear map from $X\times X$ to
$Y$. Moreover, $T$ is $c$-symmetric. Assume that $x$ is a cyclic
vector for $T^2$. Then $x$ can not lie in the right kernel of $c$,
which is a proper closed $T$-invariant subspace of $X$. Hence the
operator $\Psi\in L(X,Y)$, $\Psi(u)=c(u,x)$ is non-zero. On the
other hand, for any $n\in\Z_+$,
$$
\Psi(T^{2n}x)=c(T^{2n}x,x)=c(T^nx,T^nx)=b(T^nx,T^nx)-b(T^nx,T^nx)=0.
$$
Hence the orbit of $x$ with respect to $T^2$ lies in the proper
closed linear subspace $\ker\Psi$ of $X$, which contradicts
cyclicity of $x$ for $T^2$.
\end{proof}

This looks like a proper place to reproduce the following question
of Grivaux.

\begin{question} Let $X$ be a Banach space and $T\in L(X)$ be such
that $T\oplus T$ is cyclic. Does it follow that $T^2$ is cyclic?
\end{question}

As a straightforward toy illustration of the above proposition one
can consider the following fact. Let $(\Omega,{\cal F},\mu)$ be a
measure space, $g\in L_\infty(\mu)$, $0<p<\infty$ and $T\in
L(L_p(\mu))$ be the operator of multiplication by $g$: $Tf=fg$ for
$f\in L_p(\mu)$. Then $T\oplus T$ is non-cyclic. Indeed, consider
the continuous bilinear map $b:L_p(\mu)\times L_p(\mu)\to
L_{p/2}(\mu)$, $b(f,h)=fh$. Clearly $b$ is non-zero and $T$ is
$b$-symmetric. By Proposition~\ref{tpt}, $T\oplus T$ is non-cyclic.
The above mentioned result of Bayart and Matheron can now be
formulated in the following way.

\begin{thmbm} Let $X$ be a separable infinite dimensional Banach
space with an unconditional Schauder basis such that the forward
shift operator associated with this basis is bounded.  Then there
exists a hypercyclic $T\in L(X)$ and a non-zero continuous bilinear
form $b:X\times X\to\K$ such that $T$ is $b$-symmetric. In
particular, $T\oplus T$ is non-cyclic.
\end{thmbm}

The form $b$ in the above theorem must be symmetric. Indeed,
otherwise, by Proposition~\ref{tpt}, $T^2$ is non-cyclic, which
contradicts hypercyclicity of $T$ according to the Ansari theorem
\cite{ansa} on hypercyclicity of powers of hypercyclic operators. An
answer to the following question could help in better understanding
of the phenomenon of hypercyclic operators not satisfying the
Hypercyclicity Criterion.

\begin{question} Let $T$ be a hypercyclic continuous linear operator on a
Banach space $X$ such that $T\oplus T$ is non-hypercyclic. Does
there exist a non-zero symmetric continuous bilinear form $b:X\times
X\to \K$ such that $T\oplus T$ is $b$-symmetric?
\end{question}

It is worth noting that non-existence of such a form $b$ is
equivalent to the density of the range of the operator $I\otimes
T-T\otimes I$ acting on the projective tensor product
$X\widehat{\otimes}_\pi X$.

Another observation concerning Theorem~BM is that operators
constructed in \cite{bama} have huge spectrum. Namely, their
spectrum contains a disk centered at 0 of radius $>1$. On the other
hand, we know (see Theorems~\ref{cococo} and~\ref{cocococo}) that
any separable infinite dimensional complex Banach space supports
plenty of hypercyclic operators with the spectrum being the
singleton $\{1\}$. This leads to the following question.

\begin{question} Let $T$ be a hypercyclic continuous linear operator on
a complex Banach space $X$ such that $\sigma(T)=\{1\}$. Is $T\oplus
T$ hypercyclic?
\end{question}

It is worth noting that an affirmative answer to the above question
would take care of Problem~\ref{que}. Indeed, the spectrum of any
hypercyclic operator on a hereditarily indecomposable complex Banach
space \cite{hi1} is a singleton $\{z\}$ with $z\in\T$.

\subsection{$n$-supercyclic operators}

Recently Feldman \cite{fe1} has introduced the notion of an
$n$-supercyclic operator for $n\in\N$. A bounded linear operator $T$
on a Banach space $X$ is called {\it $n$-supercyclic} for $n\in\N$
if there exists an $n$-dimensional linear subspace $L$ of $X$ such
that its orbit $\{T^nx:n\in\Z_+,\ x\in L\}$ is dense in $X$. Such a
space $L$ is called an $n$-{\it supercyclic subspace} for $T$.
Clearly, $1$-supercyclicity coincides with the usual supercyclicity.
In \cite{fe1}, for any $n\in\N$, $n\geq 2$, a bounded linear
operator $T$ on $\ell_2$ is constructed, which is $n$-supercyclic
and not $(n-1)$-supercyclic. The construction is based on the
observation that if $T_k$ for $1\leq k\leq n$ are bounded linear
operators on Banach spaces $X_k$ all satisfying the Supercyclicity
Criterion with the same sequence $\{n_k\}$, then the direct sum
$T_1\oplus {\dots}\oplus T_n$ is $n$-supercyclic. This observation
leads to the natural question whether the direct sum of $n$
supercyclic operators should be $n$-supercyclic. The following
proposition provides a negative answer to this question.

\begin{proposition} \label{dirsum} There exists a hypercyclic operator
$T\in L(\ell_2)$ such that $T\oplus T$ is not $n$-supercyclic for
any $n\in \N$.
\end{proposition}

The above proposition follows immediately from Theorem~BM and the
next proposition, which implies that for the operator $T$ from
Theorem~BM, $T\oplus T$ is not $n$-supercyclic for any $n\in\N$.

\begin{proposition}\label{tptns} Let $X$ be an infinite dimensional
topological vector space, $b:X\times X\to \K$ be non-zero,
separately continuous and bilinear and $T\in L(X)$ be a
$b$-symmetric operator with no non-trivial closed invariant
subspaces of finite codimension. Then for any finite dimensional
linear subspace $L$ of $X\times X$, the set
$\{(p(T)x,p(T)y):p\in\pp,\ (x,y)\in L\}$ is nowhere dense in
$X\times X$.
\end{proposition}

In order to prove Proposition~\ref{tptns}, we need the following
lemma.

\begin{lemma} \label{bilin} Let $m\in\N$, $L$ and $X$ be topological
vector spaces, such that $\dim L=k\leq m$ and $B:L\times X\to \K^m$
be a continuous bilinear map such that for each non-zero $u\in L$,
the linear map $B(u,\cdot)$ is surjective. Then the set $A=\{x\in
X:B(u,x)=0\ \ \text{for some non-zero}\ \ u\in L\}$ is closed and
nowhere dense in $X$.
\end{lemma}

\begin{proof} First, observe that it is enough to prove
the required statement in the case, when $X$ is finite dimensional.
Indeed, let $e_1,\dots,e_k$ be a basis of $L$ and
$X_1=\bigcap\limits_{j=1}^k \ker B(e_j,\cdot)$. Then $X_1$ is a
closed linear subspace such that $B(u,x)=0$ for any $(u,x)\in
L\times X_1$. Moreover, $X_1$ has codimension at most $km$ in $X$.
Thus we can pick a finite dimensional subspace $Y$ of $X$ such that
$X=X_1\oplus Y$. Since $B(u,x)=0$ for any $(u,x)\in L\times X_1$, we
have that for each $u\in L\setminus\{0\}$, the restriction of
$B(u,\cdot)$ to $Y$ is onto. It is also clear that
$A=\pi^{-1}(A_0)$, where $A_0=\{x\in Y:B(u,x)=0\ \ \text{for some
non-zero}\ \ u\in L\}$ and $\pi$ is the projection in $X$ onto $Y$
along $X_1$. Since $\pi$ is continuous and open, $A$ is closed and
nowhere dense in $X$ if and only if $A_0$ is closed and nowhere
dense in $Y$, which is finite dimensional.

Thus without loss of generality, we can assume that $X$ is finite
dimensional. Consider the unit sphere $S$ in $L$ with respect to
some Hilbert space norm on $L$. For each $u\in S$, $B(u,\cdot)$ is
onto and we can pick an $m$-dimensional subspace $Z_u$ of $X$ such
that the restriction of $B(u,\cdot)$ to $Z_u$ is invertible.
Clearly, the set $V_u$ of those $v\in S$ for which the restriction
of $B(u,\cdot)$ to $Z_u$ is invertible is open in $S$ and contains
$u$. Thus we can pick a neighborhood $W_u$ of $u$ such that
$B(v,\cdot)$ is onto for each $v\in \overline{W_u}$. Since $u\in
W_u$, the family $\{W_u:u\in S\}$ is an open cover of the compact
space $S$ and therefore we can choose $u_1,\dots,u_r\in S$ such that
$S=\bigcup\limits_{j=1}^r W_{u_j}$. Then
$$
A=\bigcup_{j=1}^r A_j,\ \ \text{where}\ \ A_j=\{x\in X:B(u,x)=0\ \
\text{for some non-zero}\ \ u\in \overline{W_{u_j}}\}.
$$
It suffices to show that each $A_j$ is closed and nowhere dense. Let
$1\leq j\leq r$. Closeness of $A_j$ is rather easy. Indeed, let
$\{x_n\}_{n\in\Z_+}$ be a sequence of elements of $A_j$ converging
to $x\in X$. Since $x_n\in A_j$, we can pick $w_n\in
\overline{W_{u_j}}$ such that $B(w_n,x_n)=0$. Since
$\overline{W_{u_j}}$ is compact, we, passing to a subsequence, if
necessary, can assume that $w_n\to w\in \overline{W_{u_j}}$. Since
$B$ is continuous, we have $0=B(w_n,x_n)\to B(w,x)$. Thus $B(w,x)=0$
and $x\in A_j$. That is, $A_j$ is closed. It remains to show that it
is nowhere dense. Pick a linear subspace $Y_j$ of $X$ such that
$Z_{u_j}\oplus Y_j=X$. For each $u\in V_{u_j}$ let $T_u$ be the
restriction of $B(u,\cdot)$ to $Z_{u_j}$ and $S_u$ be the
restriction of $B(u,\cdot)$ to $Y_j$. Let $z\in Z_{u_j}$ and $y\in
Y_j$. Then $z+y\in A_j$ if and only if there exists $u\in
\overline{W_{u_j}}$ such that $T_uz+S_uy=0$. Since $T_u$ is
invertible, the latter is equivalent to $z=-T_u^{-1}S_u y$. Thus
$A_j=F(Y_j\times \overline{W_{u_j}})$, where $F:Y_j\times V_{u_j}\to
X$, $F(y,u)=y-T_u^{-1}S_uy$. Let $n=\dim X$. Then $\dim Y_j=n-m$.
Hence $Y_j\times  V_{u_j}$ is a manifold of dimension
$\alpha(n-m+k)-1$, where $\alpha=1$ if $\K=\R$ and $\alpha=2$ if
$\K=\C$. It is clear that $F$ is smooth and therefor is Lipschitzian
on any compact set. Since a Lipschitzian map does not increase the
Hausdorff dimension, we see that $A_j=F(Y_j\times
\overline{W_{u_j}})$ is a countable union of compact sets of
Hausdorff dimension at most $\alpha(n-m+k)-1<\alpha n$. Since
$X\simeq \K^{n}$, any compact subset of $X$ of Hausdorff dimension
$<\alpha n$ is nowhere dense. Thus $A_j$ is a Baire first category
set. Since $A_j$ is closed, it is nowhere dense.
\end{proof}

\begin{proof}[Proof of Proposition~$\ref{tptns}$] First, we consider
the case of non-degenerate $b$. That is, we assume that both the
left and the right kernels $X_0$ and $X_1$ of $b$ defined by
(\ref{x0x1}) are trivial. For each $k\in \Z_+$ and $(x,y)\in X\times
X$, consider the linear functional $\Phi_k(x,y)\in (X\times X)'$
defined by the formula
$$
\Phi_k(x,y)(u,v)=b(T^kx,v)-b(u,T^ky).
$$
First, we shall check that for any $(x,y)\neq (0,0)$, the
functionals $\Phi_k(x,y)$ for $k\in\Z_+$ are linearly independent.
Assume the contrary. Then there exists a non-zero polynomial $p$
such that $b(p(T)x,v)=b(u,p(T)y)$ for any $u,v\in X$. Since the
left-hand side of the last equality does not depend on $u$ and the
right-hand side does not depend on $v$, they both do not depend on
both $u$ and $v$. Hence $b(p(T)x,v)=b(u,p(T)y)=0$ for any $u,v\in
X$. Since $T$ is $b$-symmetric, we have $b(x,p(T)v)=b(p(T)u,y)=0$
for any $u,v\in X$. Hence $p(T')\phi=p(T')\psi=0$, where
$\phi,\psi\in X'$, $\phi(v)=b(x,v)$, $\psi(u)=b(u,y)$. Since $T$ has
no non-trivial closed invariant subspaces of finite codimension,
$T'$ has no non-trivial finite dimensional invariant subspaces. By
Lemma~\ref{inj}, $p(T')$ is injective. Hence $\phi=\psi=0$. Since
$b$ is non-degenerate, we then have $x=y=0$, which contradicts with
the assumption $(x,y)\neq(0,0)$. Thus the functionals $\Phi_k(x,y)$
for $k\in\Z_+$ are linearly independent for each $(x,y)\neq(0,0)$.

Let $L$ be a finite dimensional linear subspace of $X\times X$,
$\dim L=k\in\N$ and $m\in\N$, $m>k$. Consider the bilinear map
$$
B:L\times (X\times X)\to \K^m,\quad
B((x,y),(u,v))=\{\Phi_j(x,y)(u,v)\}_{j=0}^{m-1}.
$$
Since for any non-zero $(x,y)\in L$, the functionals
$\Phi_0(x,y),\dots,\Phi_{m-1}(x,y)$ are linearly independent, we see
that the linear map $B((x,y),\cdot)$ is onto. By Lemma~\ref{bilin},
the set
$$
A=\{(u,v)\in X\times X: \text{there is non-zero $(x,y)\in L$ such
that $B((x,y),(u,v))=0$}\}
$$
is closed and nowhere dense in $X\times X$. Let now $(0,0)\neq
(u,v)\in p(T\oplus T)(L)$ for some $p\in \pp$. Then
$(u,v)=(p(T)x,p(T)y)$ for some $(x,y)\in L\setminus\{(0,0)\}$. Then
$\Phi_j(x,y)(u,v)=(b(T^jx,p(T)y)-b(p(T)x,T^jy))=0$ for any
$j\in\Z_+$ since $T$ is $b$-symmetric. It follows that $(u,v)\in A$.
That is, the nowhere dense set $A$ contains the set
$\{(p(T)x,p(T)y):p\in\pp,\ (x,y)\in L\}$. Thus the latter set is
nowhere dense.

It remains to reduce the general case to the case of non-degenerate
$b$. Since $b$ is non-zero, at least one of the bilinear forms
$b_0(x,y)=b(x,y)+b(y,x)$ or $b_1(x,y)=b(x,y)-b(y,x)$ is non-zero.
Clearly $T$ is symmetric with respect to both $b_0$ and $b_1$. Thus
replacing $b$ by either $b_0$ or $b_1$, if necessary, we can assume
that $b$ is symmetric or antisymmetric. Then left and right kernels
$X_0$ and $X_1$ of $b$  defined by (\ref{x0x1}) coincide. Since $T$
is $b$-symmetric, $X_0$ is $T$-invariant. Let $T_0\in L(X/X_0)$ be
defined as $T_0(x+X_0)=Tx+X_0$ and $\beta$ be the bilinear form on
$X/X_0$ defined by the formula $\beta(x+X_0,y+X_0)=b(x,y)$. Since
$X_0$ is the right and the left kernel of $b$ and is $T$-invariant,
the operator $T_0$ and the form $\beta$ are well-defined and $\beta$
is non-degenerate. Moreover, it is easy to see that $T_0$ is
$\beta$-symmetric and has no non-trivial closed invariant subspaces
of finite codimension. Assume now that $L$ a finite dimensional
linear subspace $L$ of $X\times X$, $A=\{(p(T)x,p(T)y):p\in\pp,\
(x,y)\in L\}$ and $B=\{(p(T_0)(x+X_0),p(T_0)(y+X_0)):p\in\pp,\
(x,y)\in L\}$. According to the first part of the proof, $B$ is
nowhere dense in $X/X_0$. Since $A=\pi^{-1}(B)$, where
$\pi(x)=x+X_0$ is the canonical map from $X$ onto $X/X_0$, we see
that $A$ is nowhere dense in $X$.
\end{proof}

The following question is raised by Bourdon, Feldman and  Shapiro in
\cite{fe2}.

\begin{question} Let $X$ be a complex Banach space, $n\in\N$ and $T\in
L(X)$ be such that $T$ is $n$-supercyclic and
$\sigma_p(T')=\varnothing$. Is $T$ cyclic?
\end{question}

It is worth noting that the only known examples of $n$-supercyclic
operators $T$ with $\sigma_p(T')=\varnothing$ are the mentioned
direct sums of operators satisfying the Supercyclicity Criterion
with the same sequence $\{n_k\}$. Such direct sums are all cyclic,
see \cite{shshsh}.

Feldman in \cite{fe1} has also introduced the concept of an
$\infty$-supercyclic operator. A bounded linear operator $T$ on a
Banach space $X$ is called {\it $\infty$-supercyclic} if there
exists a linear subspace $L$ of $X$ such that its orbit
$\{T^nx:n\in\Z_+,\ x\in L\}$ is dense in $X$, the space $T^n(L)$ is
not dense in $X$ for any $n\in\Z_+$ and $L$ contains no non-zero
invariant subspace of $T$. Gallardo and Motes-Rodriguez \cite{moga},
answering a question of Salas, demonstrated that the Volterra
operator
$$
V:L_2[0,1]\to L_2[0,1],\quad Vf(x)=\int_0^x f(t)\,dt
$$
is not supercyclic. In \cite{shk2} it is shown that $V$ is not
$n$-supercyclic for any $n\in\N$. However, it turns out that $V$ is
$\infty$-supercyclic.

\begin{proposition}\label{volt} The Volterra operator is
$\infty$-supercyclic.
\end{proposition}

\begin{proof} For any non-zero $h\in L_2[0,1]$ we denote by $L_h$ the
orthocomplement of $h$: $L_h=\{f\in L_2[0,1]:\langle f,h\rangle
=0\}$. It is straightforward to see that for any $h\in L_2[0,1]$,
\begin{equation}\label{vst}
\overline{V(L_{V^\star h})}=L_h, \ \ \text{where}\ \ V^\star
f(x)=\int_x^1 f(t)dt
\end{equation}
is the adjoint of $V$. Consider the space ${\mathcal E}=\{f\in
C^\infty[0,1]:f^{(j)}(1)=0\ \text{for any}\ j\in\Z_+\}$. If $h$ is
non-zero element of $\mathcal E$, then according to the above
display, $V^\star h'=-h$. Thus by (\ref{vst})
\begin{equation}\label{vst0}
\overline{V^n(L_{h})}=L_{h^{(n)}}\text{\ \ for any non-zero}\ \
h\in{\mathcal E}\ \ \text{and any}\ \ n\in\Z_+.
\end{equation}
Consider now the following specific $h\in \mathcal E$: $h(1)=0$ and
$h(x)=e^{(x-1)^{-1}}$ for $0\leq x<1$. First, we shall show that
$L_h$ does not contain any non-zero invariant subspace of $V$.
Assume the contrary. Then there exists non-zero $f\in L_2[0,1]$ such
that $V^nf\in L_h$ for each $n\in\Z_+$. That is, $0=\langle
V^nf,h\rangle=\langle f,V^{\star n}h\rangle$. Hence $h$ is not a
cyclic vector for $V^\star$. On the other hand, it is well-known
that $g\in L_2[0,1]$ is non-cyclic for $V^\star$ if and only if
there is $q\in (0,1)$ such that $g$ vanishes on $[q,1]$. Thus $h$
vanishes on a neighborhood of $1$, which is obviously not the case.
Hence $L_h$ does not contain any non-zero invariant subspace of $V$.
According to (\ref{vst0}), $\overline{V^n(L_{h})}$ is a closed
hyperplane in $L_2[0,1]$ and therefore $V^n(L_{h})$ is not dense in
$L_2[0,1]$ for each $n\in\Z_+$. Now, in order to show that $V$ is
$\infty$-supercyclic it suffices to verify that
$A=\bigcup\limits_{n=0}^\infty V^n(L_h)$ is dense in $L_2[0,1]$. By
(\ref{vst0}), $A=\bigcup\limits_{n=0}^\infty L_{h^{(n)}}$. It is
easy to see that for any $f\in L_2[0,1]$ and non-zero $g\in
L_2[0,1]$, the distance from $f$ to $L_g$ is given by the formula
${\tt dist}\,(f,L_g)=\|g\|^{-1}|\langle f,g\rangle|$. Let $q\in
(0,1)$ and $f\in L_2[0,1]$ be such that $f$ vanishes on $[q,1]$.
Then
$$
{\tt dist}\,(f,L_{h^{(n)}})=\frac{|\langle
f,h^{(n)}\rangle|}{\|h^{(n)}\|}\leq
\|f\|\frac{\ssub{\|h^{(n)}\|}{L_2[0,q]}}{\ssub{\|h^{(n)}\|}{L_2[0,1]}}.
$$
On the other hand, analyticity of $h$ on $[0,q]$ and easy lower
estimates if $\ssub{\|f^{(n)}\|}{L_2[0,1]}$ imply that
$$
\slim_{n\to\infty}\bigr(n!\ssub{\|h^{(n)}\|}{L_2[0,q]}\bigr)^{1/n}<\infty\
\ \text{and}\ \
\lim_{n\to\infty}\bigr(n!\ssub{\|h^{(n)}\|}{L_2[0,q]}\bigr)^{1/n}=\infty.
$$
From the last two displays it follows that ${\tt
dist}\,(f,L_{h^{(n)}})\to 0$ as $n\to\infty$. Hence $\overline{A}$
contains the space of all functions vanishing on a neighborhood of
1. Since the latter space is dense in $L^2[0,1]$, $A$ is dense in
$L^2[0,1]$, which completes the proof.
\end{proof}

\subsection{$\R$-cyclicity and supercyclicity}

Let $T$ be a continuous linear operator on a separable complex
topological vector space $X$. We say that $T$ is $\R$-{\it cyclic}
if there exists $x\in X$ such that the linear span of the orbit
$\{T^nx:n\in\Z_+\}$ in $X$ considered as a linear space over $\R$ is
dense in $X$. Similarly, $T$ is called $\R$-{\it supercyclic} if
there is $x\in X$ such that $\{tT^nx:n\in\Z_+,\ t\in\R\}$ is dense
in $X$ and $T$ is called $\R^+$-{\it supercyclic} if there is $x\in
X$ such that $\{tT^nx:n\in\Z_+,\ t>0\}$ is dense in $X$. Clearly any
$\R^+$-supercyclic operator is $\R$-supercyclic and any
$\R$-supercyclic operator is $\R$-cyclic. The following theorem by
Le\'on-Saavedra and M\"uller is proved in \cite{leon}. It is proved
in the case when $X$ is a Banach space, but exactly the same proof
works in general.

\begin{thmleon} Let $T$ be a continuous linear operator on a complex
locally convex space $X$ with $\sigma_p(T')=\varnothing$. Then $T$
is supercyclic if and only if $T$ is $\R^+$-supercyclic.
\end{thmleon}

As we have shown, there are bilateral weighted shifts $T_w$ on
$\ell_2(\Z)$ with the weight sequence $w$ converging to zero
arbitrarily fast and such that $I+T_w$ and $I+T'_w$ are both
hypercyclic. This happens because we allow $w$ to behave irregularly
while still satisfying the condition
 $|w_{n}|\leq a_{|n|}$ for every $n\in \Z$ with $a$ being any sequence
of positive numbers. The following proposition shows that
hypercyclicity of $I+T_w$ is incompatible with the symmetry of the
weight sequence. Recall that a continuous linear operator on a
Banach space $X$ is called {\it weakly supercyclic} if it is
supercyclic on $X$ with weak topology.

\begin{proposition} \label{symmet} Let
$w\in \ell_\infty(\Z)$ be a weight such that $|w_n|=|w_{-n}|$ for
any $n\in \Z$, and let $p$ be a polynomial with real coefficients.
Then the operator $p(T_w)$ acting on complex $\ell_2(\Z)$ is not
$\R$-cyclic. In particular, by Theorem~{\rm LM}, $p(T_w)$ is not
weakly supercyclic.
\end{proposition}

\begin{corollary}\label{symm} Let $w\in \ell_\infty(\Z)$ with $|w_n|=|w_{-n}|$
for any $n\in \Z$. Then the operator $I+T_w$ acting on $\ell_2(\Z)$
is not weakly supercyclic.
\end{corollary}

\begin{proof}[Proof of Proposition~$\ref{symmet}$]
First, note that if $w,w'\in \ell_\infty(\Z)$ satisfy $|w_n|=|w'_n|$
for any $n\in\Z$, then $T_w$ and $T_{w'}$ are isometrically similar
with a diagonal unitary operator implementing the similarity. Thus
we can, without loss of generality, assume that $w_n\in \R$ for each
$n\in\Z$. Then the operators $T_w$ and $S=p(T_w)$ have real matrix
coefficients with respect to the canonical basis. Let $H$ be the
$\R$-subspace of $\ell_2(\Z)$ consisting of the sequences with real
entries. Then $T_w(H)\subseteq H$, $S(H)\subseteq H$, $T_w(iH)=iH$
and $S(iH)\subseteq iH$. Let $T_0$ be the restriction of $T_w$ to
$H$ and $S_0$ the restriction of $S$ to $H$ considered as
$\R$-linear operators. Since $T_0$ and $S_0$ are similar to the
restrictions of $T_w$ and $S$ to $iH$, we see that $T_w$ and $S$,
considered as $\R$-linear operators, are similar to $T_0\oplus T_0$
and $S_0\oplus S_0$. Now the symmetry of the weight sequence $w$
implies that $T_0$ is isometrically similar to $T'_0$. Indeed,
$T'_0=UT_0U^{-1}$, where $Ue_n=e_{1-n}$ for $n\in\Z$. Then
$S_0=p(T_0)$ is similar to $S'_0=p(T'_0)$. Thus $S$ considered as an
$\R$-linear operator is similar to $S_0\oplus S'_0$. The last
operator is non-cyclic. Indeed, if $x\oplus y$ is a non-zero vector
in $H\oplus H$, then the orbit of $x\oplus y$ with respect to
$S_0\oplus S'_0$ is orthogonal to the non-zero vector $y\oplus
(-x)\in H\oplus H$. Thus $S$ is not $\R$-cyclic.
\end{proof}

Since $\ell_p(\Z)$ for $1\leq p\leq 2$ is contained in $\ell_2(\Z)$
and carries a stronger topology than the one inherited from
$\ell_2(\Z)$, Proposition~\ref{symmet} remains true if we replace
$\ell_2(\Z)$ by $\ell_p(\Z)$ with $1\leq p\leq 2$. On the other
hand, the unweighted shift on $\ell_p(\Z)$ with $2<p<\infty$ is
weakly supercyclic \cite{58} and therefore the statement of
Proposition~\ref{symmet} becomes false if we replace $\ell_2(\Z)$ by
$\ell_p(\Z)$ with $p>2$. At this point it is interesting to remind
that, according to Theorem~B, norm hypercyclicity and supercyclicity
of a bilateral weighted shift on $\ell_p(\Z)$ do not depend on $p$.
This leads naturally to the following question.

\begin{question}
Characterize hypercyclicity and supercyclicity of the operators of
the form $I+T$, where $T$ is a bilateral weighted shift on
$\ell_p(\Z)$. In particular, does hypercyclicity or supercyclicity
of these operators depend on the choice of $p$, $1\leq p<+\infty $?
\end{question}

\small\rm

\vskip1truecm

\scshape

\noindent Stanislav Shkarin

\noindent Queens's University Belfast

\noindent Department of Pure Mathematics

\noindent University road, Belfast, BT7 1NN, UK

\noindent E-mail address: \qquad {\tt s.shkarin@qub.ac.uk}

\end{document}